\documentclass[11pt]{article}
\usepackage{xypic}
\usepackage{mathrsfs}
\usepackage{amssymb}
\usepackage{amsfonts}
\usepackage{amsmath}

\def\ben{\begin{enumerate}}
\def\een{\end{enumerate}}

\def\bit{\begin{itemize}}
\def\eit{\end{itemize}}

\def\0{\leqno}

\input xy
\xyoption{all}
\input{tcilatex}
\begin{document}

\begin{center}
{\Large ALGEBRAIC CONSTRUCTIONS }

\bigskip


{\Large IN THE CATEGORY OF LIE ALGEBROIDS}

\ \ \\[0pt]

\ \ \\[0pt]

\ \ \\[0pt]
by


\textbf{{CONSTANTIN M. ARCU\c{S} }}

\ \ \\[0pt]

\begin{tabular}{c}
SECONDARY SCHOOL \textquotedblleft CORNELIUS RADU\textquotedblright , \\
RADINESTI VILLAGE, 217196, GORJ COUNTY, ROMANIA \\
e-mail: c\_arcus@yahoo.com, c\_arcus@radinesti.ro%
\end{tabular}
\end{center}

\ \ \\[0pt]




%


\ \

\begin{abstract}
A generalized notion of a Lie algebroid is presented. Using this, the Lie
algebroid generalized tangent bundle is obtained. A new point of view over
(linear) connections theory on a fiber bundle is presented. These
connections are characterized by o horizontal distribution of the Lie
algebroid generalized tangent bundle. Some basic properties of these
generalized connections are investigated. Special attention to the class of
linear connections is paid. The recently studied Lie algebroids connections
can be recovered as special cases within this more general framework. In
particular, all results are similar with the classical results. Formulas of
Ricci and Bianchi type and linear connections of Levi-Civita type are
presented. \ \ \bigskip\newline
\end{abstract}

\tableofcontents

\section{Introduction}

The theory of connections constitutes one of the most important chapter of
differential geometry, which has been explored in the literature (see $%
[2,3,4,5,10,11,12,13,14,15,16]$). Connections theory has become an
indispensable tool in various branches of theoretical and mathematical
physics.

If $\left( E,\pi ,M\right) $ is a fiber bundle with paracompact base and $%
\left( VTE,\tau _{E},E\right) $ is the kernel vector bundle of the tangent $%
\mathbf{B}^{\mathbf{v}}$-morphism $\left( T\pi ,\pi \right) ,$ then we
obtain the short exact sequence%
\begin{equation*}
\begin{array}{c}
\begin{array}{ccccccccc}
0 & \hookrightarrow & VTE & \hookrightarrow & TE & ^{\underrightarrow{~\ \ \
\pi !~\ \ }} & \pi ^{\ast }TM & \longrightarrow & 0 \\
&  & ~\ \ \ \downarrow \tau _{E} &  & ~\ \ \ \downarrow \tau _{E} &  & ~\ \
\ \ \ \ \ \downarrow \pi ^{\ast }\tau _{M} &  &  \\
&  & E & ^{\underrightarrow{~\ \ \ Id_{E}~\ \ }} & E & ^{\underrightarrow{~\
\ \ Id_{E}~\ \ }} & E &  &
\end{array}%
\end{array}%
\leqno(1)
\end{equation*}%
where $\pi !$ is the projection of $TE$ onto $\pi ^{\ast }TM.$

We know that a split to the right in the previous short exact sequence, i.e.
a smooth map $h\in \mathbf{Man}\left( \pi ^{\ast }TM,TE\right) $ so that $%
\pi !\circ h=Id_{\pi ^{\ast }TM},$ is a \emph{connection in the Ehresmann
sense}.

If $\left( HTE,\tau _{E},E\right) $ is the image vector bundle of the $%
\mathbf{B}^{\mathbf{v}}$-morphism $\left( h,Id_{E}\right) ,$ then the
tangent vector bundle $\left( TE,\tau _{E},E\right) $ is a Whitney sum
between the \emph{horizontal vector bundle }$\left( HTE,\tau _{E},E\right) $
and the \emph{vertical vector bundle} $\left( VTE,\tau _{E},E\right) .$

From the above notion of connection, one can easily derive more specific
types of connections by imposing additional conditions. In the literature
one can find several generalizations of the concept of Ehresmann connection
obtained by relaxing the conditions on the horizontal vector bundle.

\begin{itemize}
\item First of all, we are thinking here of the so-called \emph{partial
connections,} where $\left( HTE,\tau _{E},E\right) $ does not determine a
full complement of $\left( VTE,\tau _{E},E\right) .$ More precisely, $\Gamma
\left( HTE,\tau _{E},E\right) $ has zero intersection with $\Gamma \left(
VTE,\tau _{E},E\right) $, but $\left( HTE,\tau _{E},E\right) $ projects onto
a vector subbundle of $\left( TM,\tau _{M},M\right) .$ (see $\left[ 7\right]
$)

\item Secondly, there also exists a notion of \emph{pseudo-connection},
introduced under the name of \emph{quasi-connection }in a paper by Y. C.
Wong $\left[ 16\right] .$ Linear pseudo-connections and generalization of it
have been studied by many authors. (see $\left[ 3\right] $)
\end{itemize}

P. Popescu build the \emph{relativ tangent space} and using that he obtained
a new \emph{generalized connection }on a vector bundle.$\left[ 11\right] $
(see also $\left[ 12\right] $)

In the paper $\left[ 4\right] $ by R. L. Fernandez a \emph{contravariant
connection} in the framework of Poisson Geometry there are presented. Given
a Poisson manifold $M$ with tensor $\Lambda $ which does not have to be of
constant rank, a \emph{covariant connection} on the principal bundle $\left(
P,\pi ,M\right) $ is a $G$-invariant bundle map $h\in \mathbf{Man}\left( \pi
^{\ast }\left( T^{\ast }M\right) ,TP\right) $ so that the diagram is
commutative%
\begin{equation*}
\begin{array}{c}
\begin{array}{ccc}
~\ \ \ \ \ \ \ \ \ \ \ \ \ \ \ \pi ^{\ast }\left( T^{\ast }M\right) & ^{%
\underrightarrow{~\ \ \ h~\ \ \ }} & TP \\
\pi ^{\ast }\left( \overset{\ast }{\tau }_{M}\right) \downarrow &  & ~\ \
\downarrow T\pi \\
~\ \ \ \ \ \ \ \ \ \ \ \ \ \ \ T^{\ast }M & ^{\underrightarrow{~\ \ \ \sharp
_{\Lambda }~\ \ \ }} & TM%
\end{array}%
\end{array}%
\leqno(2)
\end{equation*}%
where $\left( \sharp _{\Lambda },Id_{M}\right) $ is the natural vector
bundle morphism induced by the Poisson tensor. In the paper $\left[ 5\right]
,$ R. L. Fernandez has extending this theory by replacing the cotangent
bundle of a Poisson manifold by a Lie algebroid over an arbitrary manifold
and the map $\sharp _{\Lambda }$ by the anchor map of the Lie algebroid.
This resulted into a notion of \emph{Lie algebroid connection} which, in
particular, turns out to be appropriate for studying the geometry of
singular distributions.

B. Langerock and F. Cantrijn $\left[ 2\right] $ proposed a \emph{general
notion of connection }on a fiber bundle $\left( E,\pi ,M\right) ,$ as being
a smooth linear bundle map $h\in \mathbf{Man}\left( \pi ^{\ast }\left(
F\right) ,TE\right) $ so that the diagram is commutative%
\begin{equation*}
\begin{array}{c}
\begin{array}{ccc}
\pi ^{\ast }\left( F\right) & ^{\underrightarrow{~\ \ \ h~\ \ \ }} & TE \\
\downarrow &  & ~\ \ \downarrow T\pi \\
F & ^{\underrightarrow{~\ \ \ \rho ~\ \ \ }} & TM%
\end{array}%
\end{array}%
\leqno(3)
\end{equation*}%
where $\left( F,\nu ,M\right) $ is an arbitrary vector bundle and $\left(
\rho ,Id_{M}\right) $ is a vector bundle morphism of $\left( F,\nu ,M\right)
$ source and $\left( TM,\tau _{M},M\right) $ target.

Different equivalent definitions of a (linear) connection on a vector bundle
are known and there are in current usage$.$ We know the following

\textbf{Theorem} \emph{If we have a short exact sequence of vector bundles
over a paracompact manifold }$M$%
\begin{equation*}
\begin{array}{c}
^{%
\begin{array}[t]{ccccccccc}
0 & \hookrightarrow & E^{\prime } & ^{\underrightarrow{~\ \ \ f~\ \ }} & E &
^{\underrightarrow{~\ \ \ g~\ \ }} & E^{\prime \prime } & \longrightarrow & 0
\\
&  & ~\ \downarrow \pi ^{\prime } &  & ~\downarrow \pi &  & ~\ \ \downarrow
\pi ^{\prime \prime } &  &  \\
&  & M & ^{\underrightarrow{~\ \ \ Id_{M}~\ \ }} & M & ^{\underrightarrow{~\
\ \ Id_{M}~\ \ }} & M &  &
\end{array}%
}%
\end{array}%
\leqno(4)
\end{equation*}%
\emph{\ then there exists a right split if and only if there exists a left
split.}

So, a split to the left in the short exact sequence (1), i.e. a smooth map $%
\Gamma \in \mathbf{Man}\left( TE,VTE\right) $ so that $\Gamma \circ
i=Id_{TE},$ is an equivalent definition with the Ehresmann connection.

We remark that the secret of the Ehresmann connection is given by the
diagrams
\begin{equation*}
\begin{array}{c}
\begin{array}[b]{ccccc}
E &  & \left( TM,\left[ ,\right] _{TM}\right) & ^{\underrightarrow{~\ \ \
Id_{TM}~\ \ }} & \left( TM,\left[ ,\right] _{TM}\right) \\
~\downarrow \pi &  & ~\ \ \downarrow \tau _{M} &  & ~\ \ \downarrow \tau _{M}
\\
M & ^{\underrightarrow{~\ \ \ Id_{M}~\ \ }} & M & ^{\underrightarrow{~\ \ \
Id_{M}~\ \ }} & M%
\end{array}%
\end{array}%
\leqno(5)
\end{equation*}%
where $\left( E,\pi ,M\right) $ is a fiber bundle and $\left( \left( TM,\tau
_{M},M\right) ,\left[ ,\right] _{TM},\left( Id_{TM},Id_{M}\right) \right) $
is the standard Lie algebroid.

First time, appeared the idea to change the standard Lie algebroid with an
arbitrary Lie algebroid as in the diagrams
\begin{equation*}
\begin{array}{c}
\begin{array}[b]{ccccc}
E &  & \left( F,\left[ ,\right] _{F}\right) & ^{\underrightarrow{~\ \ \ \rho
\ \ \ \ }} & \left( TM,\left[ ,\right] _{TM}\right) \\
~\downarrow \pi &  & ~\downarrow \nu &  & ~\ \ \downarrow \tau _{M} \\
M & ^{\underrightarrow{~\ \ \ Id_{M}~\ \ }} & M & ^{\underrightarrow{~\ \ \
Id_{M}~\ \ }} & M%
\end{array}%
\end{array}%
\leqno(6)
\end{equation*}

Second time, appeared the idea to change in the previous diagrams the
identities morphisms with arbitrary $\mathbf{Man}$-isomorphisms $h$ and $%
\eta $ as in the diagrams
\begin{equation*}
\begin{array}{c}
\begin{array}[b]{ccccccc}
E &  & \left( F,\left[ ,\right] _{F,h}\right) & ^{\underrightarrow{~\ \ \
\rho \ \ \ \ }} & \left( TM,\left[ ,\right] _{TM}\right) & ^{%
\underrightarrow{~\ \ \ Th\ \ \ \ }} & \left( TN,\left[ ,\right] _{TN}\right)
\\
~\downarrow \pi &  & ~\downarrow \nu &  & ~\ \ \downarrow \tau _{M} &  & ~\
\ \downarrow \tau _{N} \\
M & ^{\underrightarrow{~\ \ \ h~\ \ }} & N & ^{\underrightarrow{~\ \ \ \eta
~\ \ }} & M & ^{\underrightarrow{~\ \ \ h~\ \ }} & N%
\end{array}%
\end{array}%
\leqno(7)
\end{equation*}%
where
\begin{equation*}
\left( \rho ,\eta \right) \in \mathbf{B}^{\mathbf{v}}\left( \left( F,\nu
,M\right) ,\left( TM,\tau _{M},M\right) \right)
\end{equation*}%
and
\begin{equation*}
\begin{array}{ccc}
\Gamma \left( F,\nu ,N\right) \times \Gamma \left( F,\nu ,N\right) & ^{%
\underrightarrow{~\ \ \left[ ,\right] _{F,h}~\ \ }} & \Gamma \left( F,\nu
,N\right) \\
\left( u,v\right) & \longmapsto & \ \left[ u,v\right] _{F,h}%
\end{array}%
\end{equation*}%
is an operation with the following properties:

$GLA_{1}$. \emph{the equality holds good }%
\begin{equation*}
\begin{array}{c}
\left[ u,f\cdot v\right] _{F,h}=f\left[ u,v\right] _{F,h}+\Gamma \left(
Th\circ \rho ,h\circ \eta \right) \left( u\right) f\cdot v,%
\end{array}%
\end{equation*}%
\qquad \quad\ \ \emph{for all }$u,v\in \Gamma \left( F,\nu ,N\right) $\emph{%
\ and} $f\in \mathcal{F}\left( N\right) .$

\medskip $GLA_{2}$. \emph{the }$4$\emph{-tuple} $\left( \Gamma \left( F,\nu
,N\right) ,+,\cdot ,\left[ ,\right] _{F,h}\right) $ \emph{is a Lie} $%
\mathcal{F}\left( N\right) $\emph{-algebra,}

$GLA_{3}$. \emph{the }$\mathbf{Mod}$\emph{-morphism }$\Gamma \left( Th\circ
\rho ,h\circ \eta \right) $\emph{\ is a }$\mathbf{LieAlg}$\emph{-morphism of
}%
\begin{equation*}
\left( \Gamma \left( F,\nu ,N\right) ,+,\cdot ,\left[ ,\right] _{F,h}\right)
\end{equation*}%
\emph{\ source and }%
\begin{equation*}
\left( \Gamma \left( TN,\tau _{N},N\right) ,+,\cdot ,\left[ ,\right]
_{TN}\right)
\end{equation*}
\emph{target.}

So, appeared the notion of \emph{generalized Lie algebroid} which is
presented in \emph{Definition 2.1.\ }Using this new notion we build the
\emph{Lie algebroid generalized tangent bundle }in the \emph{Theorem 3.1.}
Particularly, if $\left( \left( F,\nu ,N\right) ,\left[ ,\right] _{F},\left(
\rho ,Id_{N}\right) \right) $ is a Lie algebroid, $\left( E,\pi ,M\right)
=\left( F,\nu ,N\right) $ and $h=Id_{M},$ then we obtain a similar Lie
algebroid with the the \emph{prolongation Lie algebroid}. (see $\left[
6,8,9,10\right] $) Using this general framework, in Section 4, we propose
and develop a (linear) connections theory of Ehresmann type for fiber
bundles in general and for vector bundles in particular. It covers all types
of connections mentioned. In this general framework, we can define the
covariant derivatives of sections of a fiber bundle $\left( E,\pi ,M\right) $
with respect to sections of a generalized Lie algebroid
\begin{equation*}
\left( \left( F,\nu ,N\right) ,\left[ ,\right] _{F,h},\left( \rho ,\eta
\right) \right) .
\end{equation*}

In particular, if we use the generalized Lie algebroid structure
\begin{equation*}
\left( \left[ ,\right] _{TM,Id_{M}},\left( Id_{TM},Id_{M}\right) \right)
\end{equation*}%
for the tangent bundle $\left( TM,\tau _{M},M\right) $ in our theory, then
the linear connections obtained are similar with the classical linear
connections.

It is known that in Yang-Mills theory the set
\begin{equation*}
\begin{array}{c}
Cov_{\left( E,\pi ,M\right) }^{0}%
\end{array}%
\end{equation*}%
of covariant derivatives for the vector bundle $\left( E,\pi ,M\right) $
such that
\begin{equation*}
\begin{array}{c}
X\left( \left\langle u,v\right\rangle _{E}\right) =\left\langle D_{X}\left(
u\right) ,v\right\rangle _{E}+\left\langle u,D_{X}\left( v\right)
\right\rangle _{E},%
\end{array}%
\end{equation*}%
for any $X\in \mathcal{X}\left( M\right) $ and $u,v\in \Gamma \left( E,\pi
,M\right) ,$ is very important, because the Yang-Mills theory is a
variational theory which use (see $[1]$) the Yang-Mills functional
\begin{equation*}
\begin{array}{rcl}
Cov_{\left( E,\pi ,M\right) }^{0} & ^{\underrightarrow{\mathcal{YM}}} &
\mathbb{R} \\
D_{X} & \longmapsto & \displaystyle\frac{1}{2}\tint_{M}\left\Vert \mathbb{R}%
^{D_{X}}\right\Vert ^{2}v_{g}%
\end{array}%
\end{equation*}%
where $\mathbb{R}^{D_{X}}$ is the curvature.

Using our linear connections theory, we succeed to extend the set $%
Cov_{\left( E,\pi ,M\right) }^{0}$ of Yang-Mills theory, because using all
generalized Lie algebroid structures for the tangent bundle $\left( TM,\tau
_{M},M\right) $, we obtain all possible linear connections for the vector
bundle $\left( E,\pi ,M\right) $.

More importantly, it may bring within the reach of connection theory certain
geometric structures which have not yet been considered from such a point of
view. Finally, using our theory of linear connections, the formulas of Ricci
and Bianchi type and linear connections of Levi-Civita type are presented.

\section{Preliminaries}

In general, if $\mathcal{C}$ is a category, then we denote $\left\vert
\mathcal{C}\right\vert $ the class of objects and for any $A,B{\in }%
\left\vert \mathcal{C}\right\vert $, we denote $\mathcal{C}\left( A,B\right)
$ the set of morphisms of $A$ source and $B$ target. Let$\mathbf{~Vect},$ $%
\mathbf{Liealg},~\mathbf{Mod}$\textbf{,} $\mathbf{Man,~B}$ and $\mathbf{B}^{%
\mathbf{v}}$ be the category of real vector spaces, Lie algebras, modules,
manifolds, fiber bundles and vector bundles respectively.

We know that if $\left( E,\pi ,M\right) \in \left\vert \mathbf{B}^{\mathbf{v}%
}\right\vert ,$ $\Gamma \left( E,\pi ,M\right) =\left\{ u\in \mathbf{Man}%
\left( M,E\right) :u\circ \pi =Id_{M}\right\} $ and $\mathcal{F}\left(
M\right) =\mathbf{Man}\left( M,\mathbb{R}\right) ,$ then $\left( \Gamma
\left( E,\pi ,M\right) ,+,\cdot \right) $ is a $\mathcal{F}\left( M\right) $%
-module. If \ $\left( \varphi ,\varphi _{0}\right) \in \mathbf{B}^{\mathbf{v}%
}\left( \left( E,\pi ,M\right) ,\left( E^{\prime },\pi ^{\prime },M^{\prime
}\right) \right) $ such that $\varphi _{0}\in Iso_{\mathbf{Man}}\left(
M,M^{\prime }\right) ,$ then, using the operation
\begin{equation*}
\begin{array}{ccc}
\mathcal{F}\left( M\right) \times \Gamma \left( E^{\prime },\pi ^{\prime
},M^{\prime }\right) & ^{\underrightarrow{~\ \ \cdot ~\ \ }} & \Gamma \left(
E^{\prime },\pi ^{\prime },M^{\prime }\right) \\
\left( f,u^{\prime }\right) & \longmapsto & f\circ \varphi _{0}^{-1}\cdot
u^{\prime }%
\end{array}%
\end{equation*}%
it results that $\left( \Gamma \left( E^{\prime },\pi ^{\prime },M^{\prime
}\right) ,+,\cdot \right) $ is a $\mathcal{F}\left( M\right) $-module and we
obtain the $\mathbf{Mod}$-morphism%
\begin{equation*}
\begin{array}{ccc}
\Gamma \left( E,\pi ,M\right) & ^{\underrightarrow{~\ \ \Gamma \left(
\varphi ,\varphi _{0}\right) ~\ \ }} & \Gamma \left( E^{\prime },\pi
^{\prime },M^{\prime }\right) \\
u & \longmapsto & \Gamma \left( \varphi ,\varphi _{0}\right) u%
\end{array}%
\end{equation*}%
defined by
\begin{equation*}
\begin{array}{c}
\Gamma \left( \varphi ,\varphi _{0}\right) u\left( y\right) =\varphi \left(
u_{\varphi _{0}^{-1}\left( y\right) }\right) ,%
\end{array}%
\end{equation*}%
for any $y\in M^{\prime }.$

We know that a Lie algebroid is a vector bundle $\left( F,\nu ,N\right) \in
\left\vert \mathbf{B}^{\mathbf{v}}\right\vert $ such that there exists
\begin{equation*}
\begin{array}{c}
\left( \rho ,Id_{N}\right) \in \mathbf{B}^{\mathbf{v}}\left( \left( F,\nu
,N\right) ,\left( TN,\tau _{N},N\right) \right)%
\end{array}%
\end{equation*}%
and an operation
\begin{equation*}
\begin{array}{ccc}
\Gamma \left( F,\nu ,N\right) \times \Gamma \left( F,\nu ,N\right) & ^{%
\underrightarrow{\,\left[ ,\right] _{F}\,}} & \Gamma \left( F,\nu ,N\right)
\\
\left( u,v\right) & \longmapsto & \left[ u,v\right] _{F}%
\end{array}%
\end{equation*}%
with the following properties:

\begin{itemize}
\item[$LA_{1}$.] the equality holds good
\begin{equation*}
\begin{array}{c}
\left[ u,f\cdot v\right] _{F}=f\left[ u,v\right] _{F}+\Gamma \left( \rho
,Id_{N}\right) \left( u\right) f\cdot v,%
\end{array}%
\end{equation*}%
for all $u,v\in \Gamma \left( F,\nu ,N\right) $ and $f\in \mathcal{F}\left(
N\right) ,$

\item[$LA_{2}$.] the $4$-tuple $\left( \Gamma \left( F,\nu ,N\right)
,+,\cdot ,\left[ ,\right] _{F}\right) $ is a Lie $\mathcal{F}\left( N\right)
$-algebra$,$

\item[$LA_{3}$.] the $\mathbf{Mod}$-morphism $\Gamma \left( \rho
,Id_{N}\right) $ is a $\mathbf{LieAlg}$-morphism of
\begin{equation*}
\left( \Gamma \left( F,\nu ,N\right) ,+,\cdot ,\left[ ,\right] _{F}\right)
\end{equation*}
source and
\begin{equation*}
\left( \Gamma \left( TN,\tau _{N},N\right) ,+,\cdot ,\left[ ,\right]
_{TN}\right)
\end{equation*}
target.
\end{itemize}

\textbf{Definition 2.1 }Let $M,N\in \left\vert \mathbf{Man}\right\vert ,$ $%
h\in Iso_{\mathbf{Man}}\left( M,N\right) $ and $\eta \in Iso_{\mathbf{Man}%
}\left( N,M\right) $.

If $\left( F,\nu ,N\right) \in \left\vert \mathbf{B}^{\mathbf{v}}\right\vert
$ so that there exists
\begin{equation*}
\begin{array}{c}
\left( \rho ,\eta \right) \in \mathbf{B}^{\mathbf{v}}\left( \left( F,\nu
,N\right) ,\left( TM,\tau _{M},M\right) \right)%
\end{array}%
\end{equation*}%
and an operation
\begin{equation*}
\begin{array}{ccc}
\Gamma \left( F,\nu ,N\right) \times \Gamma \left( F,\nu ,N\right) & ^{%
\underrightarrow{\left[ ,\right] _{F,h}}} & \Gamma \left( F,\nu ,N\right) \\
\left( u,v\right) & \longmapsto & \left[ u,v\right] _{F,h}%
\end{array}%
\end{equation*}%
with the following properties:\bigskip

\noindent $\qquad GLA_{1}$. the equality holds good
\begin{equation*}
\begin{array}{c}
\left[ u,f\cdot v\right] _{F,h}=f\left[ u,v\right] _{F,h}+\Gamma \left(
Th\circ \rho ,h\circ \eta \right) \left( u\right) f\cdot v,%
\end{array}%
\end{equation*}%
\qquad \quad\ \ for all $u,v\in \Gamma \left( F,\nu ,N\right) $ and $f\in
\mathcal{F}\left( N\right) .$

\medskip $GLA_{2}$. the $4$-tuple $\left( \Gamma \left( F,\nu ,N\right)
,+,\cdot ,\left[ ,\right] _{F,h}\right) $ is a Lie $\mathcal{F}\left(
N\right) $-algebra,

$GLA_{3}$. the $\mathbf{Mod}$-morphism $\Gamma \left( Th\circ \rho ,h\circ
\eta \right) $ is a $\mathbf{LieAlg}$-morphism of
\begin{equation*}
\left( \Gamma \left( F,\nu ,N\right) ,+,\cdot ,\left[ ,\right] _{F,h}\right)
\end{equation*}
source and
\begin{equation*}
\left( \Gamma \left( TN,\tau _{N},N\right) ,+,\cdot ,\left[ ,\right]
_{TN}\right)
\end{equation*}
target, \medskip \noindent then we will say that \emph{the triple }%
\begin{equation*}
\begin{array}{c}
\left( \left( F,\nu ,N\right) ,\left[ ,\right] _{F,h},\left( \rho ,\eta
\right) \right)%
\end{array}%
\leqno(2.1)
\end{equation*}%
\emph{is a generalized Lie algebroid. }The couple $\left( \left[ ,\right]
_{F,h},\left( \rho ,\eta \right) \right) $will be called \emph{generalized
Lie algebroid structure.}

\emph{Remark 2.1 In the particular case, }$\left( \eta ,h\right) =\left(
Id_{M},Id_{M}\right) ,$\emph{\ we obtain the definition of Lie algebroid.}

Let $\left( \left( F,\nu ,N\right) ,\left[ ,\right] _{F,h},\left( \rho ,\eta
\right) \right) $ be a generalized Lie algebroid.

\begin{itemize}
\item Locally, for any $\alpha ,\beta \in \overline{1,p},$ we set $\left[
t_{\alpha },t_{\beta }\right] _{F,h}=L_{\alpha \beta }^{\gamma }t_{\gamma }.$
We easily obtain that $L_{\alpha \beta }^{\gamma }=-L_{\beta \alpha
}^{\gamma },~$for any $\alpha ,\beta ,\gamma \in \overline{1,p}.$
\end{itemize}

The real local functions $L_{\alpha \beta }^{\gamma },~\alpha ,\beta ,\gamma
\in \overline{1,p}$ will be called the \emph{structure functions of the
generalized Lie algebroid }$\left( \left( F,\nu ,N\right) ,\left[ ,\right]
_{F,h},\left( \rho ,\eta \right) \right) .$

\begin{itemize}
\item We assume the following diagrams:%
\begin{equation*}
\begin{array}[b]{ccccc}
F & ^{\underrightarrow{~\ \ \ \rho ~\ \ }} & TM & ^{\underrightarrow{~\ \ \
Th~\ \ }} & TN \\
~\downarrow \nu &  & ~\ \ \ \downarrow \tau _{M} &  & ~\ \ \ \downarrow \tau
_{N} \\
N & ^{\underrightarrow{~\ \ \ \eta ~\ \ }} & M & ^{\underrightarrow{~\ \ \
h~\ \ }} & N \\
&  &  &  &  \\
\left( \chi ^{\tilde{\imath}},z^{\alpha }\right) &  & \left(
x^{i},y^{i}\right) &  & \left( \chi ^{\tilde{\imath}},z^{\tilde{\imath}%
}\right)%
\end{array}%
\end{equation*}

where $i,\tilde{\imath}\in \overline{1,m}$ and $\alpha \in \overline{1,p}.$

If%
\begin{equation*}
\left( \chi ^{\tilde{\imath}},z^{\alpha }\right) \longrightarrow \left( \chi
^{\tilde{\imath}\prime }\left( \chi ^{\tilde{\imath}}\right) ,z^{\alpha
\prime }\left( \chi ^{\tilde{\imath}},z^{\alpha }\right) \right) ,
\end{equation*}%
\begin{equation*}
\left( x^{i},y^{i}\right) \longrightarrow \left( x^{i%
{\acute{}}%
}\left( x^{i}\right) ,y^{i%
{\acute{}}%
}\left( x^{i},y^{i}\right) \right)
\end{equation*}%
and
\begin{equation*}
\left( \chi ^{\tilde{\imath}},z^{\tilde{\imath}}\right) \longrightarrow
\left( \chi ^{\tilde{\imath}\prime }\left( \chi ^{\tilde{\imath}}\right) ,z^{%
\tilde{\imath}\prime }\left( \chi ^{\tilde{\imath}},z^{\tilde{\imath}%
}\right) \right) ,
\end{equation*}%
then
\begin{equation*}
\begin{array}[b]{c}
z^{\alpha
{\acute{}}%
}=\Lambda _{\alpha }^{\alpha
{\acute{}}%
}z^{\alpha }%
\end{array}%
,
\end{equation*}%
\begin{equation*}
\begin{array}[b]{c}
y^{i%
{\acute{}}%
}=\frac{\partial x^{i%
{\acute{}}%
}}{\partial x^{i}}y^{i}%
\end{array}%
\end{equation*}%
and
\begin{equation*}
\begin{array}{c}
z^{\tilde{\imath}\prime }=\frac{\partial \chi ^{\tilde{\imath}\prime }}{%
\partial \chi ^{\tilde{\imath}}}z^{\tilde{\imath}}.%
\end{array}%
\end{equation*}

\item We assume that $\left( \theta ,\mu \right) \overset{put}{=}\left(
Th\circ \rho ,h\circ \eta \right) $. If $z^{\alpha }t_{\alpha }\in \Gamma
\left( F,\nu ,N\right) $ is arbitrary, then
\begin{equation*}
\begin{array}[t]{l}
\displaystyle%
\begin{array}{c}
\Gamma \left( Th\circ \rho ,h\circ \eta \right) \left( z^{\alpha }t_{\alpha
}\right) f\left( h\circ \eta \left( \varkappa \right) \right) =\vspace*{1mm}
\\
=\left( \theta _{\alpha }^{\tilde{\imath}}z^{\alpha }\frac{\partial f}{%
\partial \varkappa ^{\tilde{\imath}}}\right) \left( h\circ \eta \left(
\varkappa \right) \right) =\left( \left( \rho _{\alpha }^{i}\circ h\right)
\left( z^{\alpha }\circ h\right) \frac{\partial f\circ h}{\partial x^{i}}%
\right) \left( \eta \left( \varkappa \right) \right) ,%
\end{array}%
\end{array}%
\leqno(2.2)
\end{equation*}%
for any $f\in \mathcal{F}\left( N\right) $ and $\varkappa \in N.$
\end{itemize}

The coefficients $\rho _{\alpha }^{i}$ respectively $\theta _{\alpha }^{%
\tilde{\imath}}$ change to $\rho _{\alpha
{\acute{}}%
}^{i%
{\acute{}}%
}$ respectively $\theta _{\alpha
{\acute{}}%
}^{\tilde{\imath}%
{\acute{}}%
}$ according to the rule:
\begin{equation*}
\begin{array}{c}
\rho _{\alpha
{\acute{}}%
}^{i%
{\acute{}}%
}=\Lambda _{\alpha
{\acute{}}%
}^{\alpha }\rho _{\alpha }^{i}\displaystyle\frac{\partial x^{i%
{\acute{}}%
}}{\partial x^{i}},%
\end{array}%
\leqno(2.3)
\end{equation*}%
respectively%
\begin{equation*}
\begin{array}{c}
\theta _{\alpha
{\acute{}}%
}^{\tilde{\imath}%
{\acute{}}%
}=\Lambda _{\alpha
{\acute{}}%
}^{\alpha }\theta _{\alpha }^{\tilde{\imath}}\displaystyle\frac{\partial
\varkappa ^{\tilde{\imath}%
{\acute{}}%
}}{\partial \varkappa ^{\tilde{\imath}}},%
\end{array}%
\leqno(2.4)
\end{equation*}%
where
\begin{equation*}
\left\Vert \Lambda _{\alpha
{\acute{}}%
}^{\alpha }\right\Vert =\left\Vert \Lambda _{\alpha }^{\alpha
{\acute{}}%
}\right\Vert ^{-1}.
\end{equation*}

\emph{Remark 2.2 }\emph{The following equalities hold good:}%
\begin{equation*}
\begin{array}{c}
\displaystyle\rho _{\alpha }^{i}\circ h\frac{\partial f\circ h}{\partial
x^{i}}=\left( \theta _{\alpha }^{\tilde{\imath}}\frac{\partial f}{\partial
\varkappa ^{\tilde{\imath}}}\right) \circ h,\forall f\in \mathcal{F}\left(
N\right) .%
\end{array}%
\leqno(2.5)
\end{equation*}%
\emph{and }%
\begin{equation*}
\begin{array}{c}
\displaystyle\left( L_{\alpha \beta }^{\gamma }\circ h\right) \left( \rho
_{\gamma }^{k}\circ h\right) =\left( \rho _{\alpha }^{i}\circ h\right) \frac{%
\partial \left( \rho _{\beta }^{k}\circ h\right) }{\partial x^{i}}-\left(
\rho _{\beta }^{j}\circ h\right) \frac{\partial \left( \rho _{\alpha
}^{k}\circ h\right) }{\partial x^{j}}.%
\end{array}%
\leqno(2.6)
\end{equation*}

\textbf{Theorem 2.1} \emph{Let} $M,N\in \left\vert \mathbf{Man}\right\vert ,$
$h\in Iso_{\mathbf{Man}}\left( M,N\right) $ \emph{and} $\eta \in Iso_{%
\mathbf{Man}}\left( N,M\right) $\emph{\ be. Using the tangent }$\mathbf{B}^{%
\mathbf{v}}$\emph{-morphism }$\left( T\eta ,\eta \right) $\emph{\ and the
operation }%
\begin{equation*}
\begin{array}{ccc}
\Gamma \left( TN,\tau _{N},N\right) \times \Gamma \left( TN,\tau
_{N},N\right) & ^{\underrightarrow{~\ \ \left[ ,\right] _{TN,h}~\ \ }} &
\Gamma \left( TN,\tau _{N},N\right) \\
\left( u,v\right) & \longmapsto & \ \left[ u,v\right] _{TN,h}%
\end{array}%
\end{equation*}%
\emph{where }%
\begin{equation*}
\left[ u,v\right] _{TN,h}=\Gamma \left( T\left( h\circ \eta \right)
^{-1},\left( h\circ \eta \right) ^{-1}\right) \left( \left[ \Gamma \left(
T\left( h\circ \eta \right) ,h\circ \eta \right) u,\Gamma \left( T\left(
h\circ \eta \right) ,h\circ \eta \right) v\right] _{TN}\right) ,
\end{equation*}%
\emph{for any }$u,v\in \Gamma \left( TN,\tau _{N},N\right) $\emph{, we
obtain that}%
\begin{equation*}
\begin{array}{c}
\left( \left( TN,\tau _{N},N\right) ,\left( T\eta ,\eta \right) ,\left[ ,%
\right] _{TN,h}\right)%
\end{array}%
\end{equation*}%
\emph{is a generalized Lie algebroid.}

For any $\mathbf{Man}$-isomorphisms $\eta $ and $h,$ new and interesting
generalized Lie algebroid structures for the tangent vector bundle $\left(
TN,\tau _{N},N\right) $ are obtained$.$ For any base $\left\{ t_{\alpha
},~\alpha \in \overline{1,m}\right\} $ of the module of sections $\left(
\Gamma \left( TN,\tau _{N},N\right) ,+,\cdot \right) $ we obtain the
structure functions%
\begin{equation*}
\begin{array}{c}
L_{\alpha \beta }^{\gamma }=\left( \theta _{\alpha }^{i}\frac{\partial
\theta _{\beta }^{j}}{\partial x^{i}}-\theta _{\beta }^{i}\frac{\partial
\theta _{\alpha }^{j}}{\partial x^{i}}\right) \tilde{\theta}_{j}^{\gamma
},~\alpha ,\beta ,\gamma \in \overline{1,m}%
\end{array}%
\end{equation*}%
where
\begin{equation*}
\theta _{\alpha }^{i},~i,\alpha \in \overline{1,m}
\end{equation*}%
are real local functions so that
\begin{equation*}
\begin{array}{c}
\Gamma \left( T\left( h\circ \eta \right) ,h\circ \eta \right) \left(
t_{\alpha }\right) =\theta _{\alpha }^{i}\frac{\partial }{\partial x^{i}}%
\end{array}%
\end{equation*}%
and
\begin{equation*}
\tilde{\theta}_{j}^{\gamma },~i,\gamma \in \overline{1,m}
\end{equation*}%
are real local functions so that
\begin{equation*}
\begin{array}{c}
\Gamma \left( T\left( h\circ \eta \right) ^{-1},\left( h\circ \eta \right)
^{-1}\right) \left( \frac{\partial }{\partial x^{j}}\right) =\tilde{\theta}%
_{j}^{\gamma }t_{\gamma }.%
\end{array}%
\end{equation*}

In particular, using arbitrary isometries (symmetries, translations,
rotations,...) for the Euclidean $3$-dimensional space $\Sigma ,$ and
arbitrary basis for the module of sections we obtain a lot of generalized
Lie algebroid structures for the tangent vector bundle $\left( T\Sigma ,\tau
_{\Sigma },\Sigma \right) $.

Let $\left( \left( F,\nu ,M\right) ,\left[ ,\right] _{F},\left( \rho
,Id_{M}\right) \right) $ be a Lie algebroid and let $h\in Iso_{\mathbf{Man}%
}\left( M,M\right) $ be. Let $\mathcal{AF}_{F}$ be a vector fibred $\left(
m+p\right) $-atlas for the vector bundle $\left( F,\nu ,M\right) $ and let $%
\mathcal{AF}_{TM}$ be a vector fibred $\left( m+m\right) $-atlas for the
vector bundle $\left( TM,\tau _{M},M\right) $.

If $\left( U,\xi _{U}\right) \in \mathcal{AF}_{TM}$ and $\left(
V,s_{V}\right) \in \mathcal{AF}_{F}$ such that $U\cap h^{-1}\left( V\right)
\neq \phi $, then we define the application%
\begin{equation*}
\begin{array}{ccc}
\tau _{N}^{-1}(U{\cap }h^{-1}(V))) & {}^{\underrightarrow{\bar{\xi}_{U{\cap }%
h^{-1}(V)}}} & \left( U{\cap }h^{-1}(V)\right) {\times }\mathbb{R}^{m} \\
\left( \varkappa ,u\left( \varkappa \right) \right) & \longmapsto & \left(
\varkappa ,\xi _{U,\varkappa }^{-1}u\left( \varkappa \right) \right) .%
\end{array}%
\end{equation*}

\textbf{Proposition 2.1 }\emph{The set}%
\begin{equation*}
\begin{array}{c}
\overline{\mathcal{AF}}_{TM}\overset{put}{=}\underset{U\cap h^{-1}\left(
V\right) \neq \phi }{\underset{\left( U,\xi _{U}\right) \in \mathcal{AF}%
_{TM},~\left( V,s_{V}\right) \in \mathcal{AF}_{F}}{\tbigcup }}\left\{ \left(
U\cap h^{-1}\left( V\right) ,\bar{\xi}_{U{\cap }h^{-1}(V)}\right) \right\}%
\end{array}%
\end{equation*}%
\emph{is a vector fibred }$m+m$\emph{-atlas of the vector bundle }$\left(
TM,\tau _{M},M\right) .$

\emph{If }$X=X^{\tilde{\imath}}\frac{\partial }{\partial \varkappa ^{\tilde{%
\imath}}}\in \Gamma \left( TM,\tau _{M},M\right) ,$ \emph{then\ we obtain
the section }%
\begin{equation*}
\begin{array}{c}
\bar{X}=\bar{X}^{\tilde{\imath}}\circ h\frac{\partial }{\partial \bar{%
\varkappa}^{\tilde{\imath}}}\in \Gamma \left( TM,\tau _{M},M\right) ,%
\end{array}%
\end{equation*}%
\emph{such that }$\bar{X}\left( \bar{\varkappa}\right) =X\left( h\left( \bar{%
\varkappa}\right) \right) ,$ \emph{for any }$\bar{\varkappa}\in U\cap
h^{-1}\left( V\right) .$

The set $\left\{ \frac{\partial }{\partial \bar{\varkappa}^{\tilde{\imath}}}%
,~\tilde{\imath}\in \overline{1,m}\right\} $ is the natural base of the $%
\mathcal{F}\left( M\right) $-module $\left( \Gamma \left( TM,\tau
_{M},M\right) ,+,\cdot \right) .$

\emph{Remark 2.3 }If $\left( \left( F,\nu ,N\right) ,\left[ ,\right]
_{F,h},\left( \rho ,\eta \right) \right) $ is a generalized Lie algebroid,
then we obtain the inclusion $\mathbf{B}^{\mathbf{v}}$-morphism%
\begin{equation*}
\begin{array}{ccc}
~\ \ \ \ \ \ \pi ^{\ast }\ \left( h^{\ast }F\right) & \hookrightarrow & F \\
h^{\ast }\nu \downarrow &  & ~\downarrow \nu \\
~\ \ \ \ \ E & ^{\underrightarrow{~\ \ h\circ \pi ~\ \ }} & M%
\end{array}%
\leqno(2.7)
\end{equation*}

\section{The Lie algebroid generalized tangent bundle}

We consider the following diagram:
\begin{equation*}
\begin{array}{c}
\xymatrix{E\ar[d]_\pi&\left( F,\left[ ,\right] _{F,h},\left( \rho ,\eta
\right) \right)\ar[d]^\nu\\ M\ar[r]^h&N}%
\end{array}%
\leqno(3.1)
\end{equation*}%
where $\left( E,\pi ,M\right) $ is a fiber bundle and $\left( \left( F,\nu
,M\right) ,\left[ ,\right] _{F,h},\left( \rho ,\eta \right) \right) $ is a
generalized Lie algebroid.

We take $\left( x^{i},y^{a}\right) $ as canonical local coordinates on $%
\left( E,\pi ,M\right) ,$ where $i\in \overline{1,m}$ and $a\in \overline{1,r%
}.$ Let
\begin{equation*}
\left( x^{i},y^{a}\right) \longrightarrow \left( x^{i%
{\acute{}}%
}\left( x^{i}\right) ,y^{a%
{\acute{}}%
}\left( x^{i},y^{a}\right) \right)
\end{equation*}%
be a change of coordinates on $\left( E,\pi ,M\right) $. Then the
coordinates $y^{a}$ change to $y^{a%
{\acute{}}%
}$ according to the rule:
\begin{equation*}
\begin{array}{c}
y^{a%
{\acute{}}%
}=\displaystyle\frac{\partial y^{a%
{\acute{}}%
}}{\partial y^{a}}y^{a}.%
\end{array}%
\leqno(3.2)
\end{equation*}

In particular, if $\left( E,\pi ,M\right) $ is vector bundle, then the
coordinates $y^{a}$ change to $y^{a%
{\acute{}}%
}$ according to the rule:
\begin{equation*}
\begin{array}{c}
y^{a%
{\acute{}}%
}=M_{a}^{a%
{\acute{}}%
}y^{a}.%
\end{array}%
\leqno(3.2^{\prime })
\end{equation*}

Easily we obtain the following

\textbf{Theorem 3.1 }\emph{Let} $\Big({\overset{\pi ^{\ast }\ \left( h^{\ast
}F\right) }{\rho }},Id_{E}\Big)$ \emph{be the }$\mathbf{B}^{\mathbf{v}}$%
\emph{-morphism of }$\left( \pi ^{\ast }\ \left( h^{\ast }F\right) ,\pi
^{\ast }\left( h^{\ast }\nu \right) ,M\right) $\ \emph{source and} $\left(
TM,\tau _{M},M\right) $\ \emph{target, where}%
\begin{equation*}
\begin{array}{rcl}
\ \pi ^{\ast }\ \left( h^{\ast }F\right) & ^{\underrightarrow{\overset{\pi
^{\ast }\ \left( h^{\ast }F\right) }{\rho }}} & TE \\
\displaystyle Z^{\alpha }T_{\alpha }\left( u_{x}\right) & \longmapsto & %
\displaystyle\left( Z^{\alpha }\cdot \rho _{\alpha }^{i}\circ h\circ \pi
\right) \frac{\partial }{\partial x^{i}}\left( u_{x}\right)%
\end{array}%
\leqno(3.3)
\end{equation*}

\emph{Using the operation}
\begin{equation*}
\begin{array}{ccc}
\Gamma \left( \pi ^{\ast }\ \left( h^{\ast }F\right) ,\pi ^{\ast }\left(
h^{\ast }\nu \right) ,M\right) ^{2} & ^{\underrightarrow{~\ \ \left[ ,\right]
_{\pi ^{\ast }\ \left( h^{\ast }F\right) }~\ \ }} & \Gamma \left( \pi ^{\ast
}\ \left( h^{\ast }F\right) ,\pi ^{\ast }\left( h^{\ast }\nu \right)
,M\right)%
\end{array}%
\end{equation*}%
\emph{defined by}%
\begin{equation*}
\begin{array}{ll}
\left[ T_{\alpha },T_{\beta }\right] _{\pi ^{\ast }\ \left( h^{\ast
}F\right) } & =L_{\alpha \beta }^{\gamma }\circ h\circ \pi \cdot T_{\gamma },%
\vspace*{1mm} \\
\left[ T_{\alpha },fT_{\beta }\right] _{\pi ^{\ast }\ \left( h^{\ast
}F\right) } & \displaystyle=fL_{\alpha \beta }^{\gamma }\circ h\circ \pi
T_{\gamma }+\rho _{\alpha }^{i}\circ h\circ \pi \frac{\partial f}{\partial
x^{i}}T_{\beta },\vspace*{1mm} \\
\left[ fT_{\alpha },T_{\beta }\right] _{\pi ^{\ast }\ \left( h^{\ast
}F\right) } & =-\left[ T_{\beta },fT_{\alpha }\right] _{\pi ^{\ast }\ \left(
h^{\ast }F\right) },%
\end{array}%
\leqno(3.4)
\end{equation*}%
\emph{for any} $f\in \mathcal{F}\left( E\right) ,$ \emph{it results that}
\begin{equation*}
\begin{array}{c}
\left( \left( \pi ^{\ast }\ \left( h^{\ast }F\right) ,\pi ^{\ast }\left(
h^{\ast }\nu \right) ,M\right) ,\left[ ,\right] _{\pi ^{\ast }\ \left(
h^{\ast }F\right) },\left( \overset{\pi ^{\ast }\ \left( h^{\ast }F\right) }{%
\rho },Id_{E}\right) \right)%
\end{array}%
\end{equation*}%
\emph{is a Lie algebroid which is called the pull-back Lie algebroid of the
generalized Lie algebroid }$\left( \left( F,\nu ,M\right) ,\left[ ,\right]
_{F,h},\left( \rho ,\eta \right) \right) $

If $z=z^{\alpha }t_{\alpha }\in \Gamma \left( F,\nu ,M\right) ,$ then we
obtain the section%
\begin{equation*}
Z=\left( z^{\alpha }\circ h\circ \pi \right) T_{\alpha }\in \Gamma \left(
\pi ^{\ast }\left( h^{\ast }F\right) ,\pi ^{\ast }\left( h^{\ast }\nu
\right) ,E\right)
\end{equation*}%
so that $Z\left( u_{x}\right) =z\left( h\left( x\right) \right) ,$ for any $%
u_{x}\in \pi ^{-1}\left( U{\cap h}^{-1}V\right) .$

Let
\begin{equation*}
\begin{array}[b]{c}
\left( \partial _{i},\dot{\partial}_{a}\right) \overset{put}{=}\left( \frac{%
\partial }{\partial x^{i}},\frac{\partial }{\partial y^{a}}\right)%
\end{array}%
\end{equation*}
be the base sections for the Lie $\mathcal{F}\left( E\right) $-algebra
\begin{equation*}
\left( \Gamma \left( TE,\tau _{E},E\right) ,+,\cdot ,\left[ ,\right]
_{TE}\right) .
\end{equation*}

For any sections%
\begin{equation*}
\begin{array}{c}
Z^{\alpha }T_{\alpha }\in \Gamma \left( \pi ^{\ast }\left( h^{\ast }F\right)
,\pi ^{\ast }\left( h^{\ast }F\right) ,E\right)%
\end{array}%
\end{equation*}%
and%
\begin{equation*}
\begin{array}{c}
Y^{a}\dot{\partial}_{a}\in \Gamma \left( VTE,\tau _{E},E\right)%
\end{array}%
\end{equation*}%
we obtain the section
\begin{equation*}
\begin{array}{c}
Z^{\alpha }\tilde{\partial}_{\alpha }+Y^{a}\overset{\cdot }{\tilde{\partial}}%
_{a}=:Z^{\alpha }\left( T_{\alpha }\oplus \left( \rho _{\alpha }^{i}\circ
h\circ \pi \right) \partial _{i}\right) +Y^{a}\left( 0_{\pi ^{\ast }\left(
h^{\ast }F\right) }\oplus \dot{\partial}_{a}\right) \vspace*{1mm} \\
=Z^{\alpha }T_{\alpha }\oplus \left( Z^{\alpha }\left( \rho _{\alpha
}^{i}\circ h\circ \pi \right) \partial _{i}+Y^{a}\dot{\partial}_{a}\right)
\in \Gamma \left( \pi ^{\ast }\left( h^{\ast }F\right) \oplus TE,\overset{%
\oplus }{\pi },E\right) .%
\end{array}%
\end{equation*}

Since we have
\begin{equation*}
\begin{array}{c}
Z^{\alpha }\displaystyle\tilde{\partial}_{\alpha }+Y^{a}\overset{\cdot }{%
\tilde{\partial}}_{a}=0 \\
\Updownarrow \\
Z^{\alpha }T_{\alpha }=0~\wedge Z^{\alpha }\left( \rho _{\alpha }^{i}\circ
h\circ \pi \right) \partial _{i}+Y^{a}\dot{\partial}_{a}=0,%
\end{array}%
\end{equation*}%
it implies $Z^{\alpha }=0,~\alpha \in \overline{1,p}$ and $Y^{a}=0,~a\in
\overline{1,r}.$

Therefore, the sections $\tilde{\partial}_{1},...,\tilde{\partial}_{p},%
\overset{\cdot }{\tilde{\partial}}_{1},...,\overset{\cdot }{\tilde{\partial}}%
_{r}$ are linearly independent.\smallskip

We consider the vector subbundle $\left( \left( \rho ,\eta \right) TE,\left(
\rho ,\eta \right) \tau _{E},E\right) $ of the vector bundle\break $\left(
\pi ^{\ast }\left( h^{\ast }F\right) \oplus TE,\overset{\oplus }{\pi }%
,E\right) ,$ for which the $\mathcal{F}\left( E\right) $-module of sections
is the $\mathcal{F}\left( E\right) $-submodule of $\left( \Gamma \left( \pi
^{\ast }\left( h^{\ast }F\right) \oplus TE,\overset{\oplus }{\pi },E\right)
,+,\cdot \right) ,$ generated by the set of sections $\left( \tilde{\partial}%
_{\alpha },\overset{\cdot }{\tilde{\partial}}_{a}\right) .$

The base sections $\left( \tilde{\partial}_{\alpha },\overset{\cdot }{\tilde{%
\partial}}_{a}\right) $ will be called the \emph{natural }$\left( \rho ,\eta
\right) $\emph{-base.}

The matrix of coordinate transformation on $\left( \left( \rho ,\eta \right)
TE,\left( \rho ,\eta \right) \tau _{E},E\right) $ at a change of fibred
charts is
\begin{equation*}
\left\Vert
\begin{array}{cc}
\Lambda _{\alpha }^{\alpha
{\acute{}}%
}\circ h\circ \pi & 0\vspace*{1mm} \\
\left( \rho _{\alpha }^{i}\circ h\circ \pi \right) \displaystyle\frac{%
\partial y^{a%
{\acute{}}%
}}{\partial x^{i}} & \displaystyle\frac{\partial y^{a%
{\acute{}}%
}}{\partial y^{a}}%
\end{array}%
\right\Vert .\leqno(3.5)
\end{equation*}

In particular, if $\left( E,\pi ,M\right) $ is a vector bundle, then the
matrix of coordinate transformation on $\left( \left( \rho ,\eta \right)
TE,\left( \rho ,\eta \right) \tau _{E},E\right) $ at a change of fibred
charts is
\begin{equation*}
\left\Vert
\begin{array}{cc}
\Lambda _{\alpha }^{\alpha
{\acute{}}%
}\circ h\circ \pi & 0\vspace*{1mm} \\
\left( \rho _{a}^{i}\circ h\circ \pi \right) \displaystyle\frac{\partial
M_{b}^{a%
{\acute{}}%
}\circ \pi }{\partial x_{i}}y^{b} & M_{a}^{a%
{\acute{}}%
}\circ \pi%
\end{array}%
\right\Vert .\leqno(3.6)
\end{equation*}

Easily we obtain

\textbf{Theorem 3.1 }\emph{Let} $\left( \tilde{\rho},Id_{E}\right) $\ \emph{%
be the} $\mathbf{B}^{\mathbf{v}}$\emph{-morphism of }$\left( \left( \rho
,\eta \right) TE,\left( \rho ,\eta \right) \tau _{E},E\right) $\ \emph{%
source and }$\left( TE,\tau _{E},E\right) $\ \emph{target, where}
\begin{equation*}
\begin{array}{rcl}
\left( \rho ,\eta \right) TE\!\!\! & \!\!^{\underrightarrow{\tilde{\ \ \rho
\ \ }}}\!\!\! & \!\!TE\vspace*{2mm} \\
\left( Z^{\alpha }\tilde{\partial}_{\alpha }+Y^{a}\overset{\cdot }{\tilde{%
\partial}}_{a}\right) \!(u_{x})\!\!\!\! & \!\!\longmapsto \!\!\! &
\!\!\left( \!Z^{\alpha }\!\left( \rho _{\alpha }^{i}{\circ }h{\circ }\pi
\!\right) \!\partial _{i}{+}Y^{a}\dot{\partial}_{a}\right) \!(u_{x})\!\!.%
\end{array}%
\leqno(3.7)
\end{equation*}

\emph{Using the operation}
\begin{equation*}
\begin{array}{ccc}
\Gamma \left( \left( \rho ,\eta \right) TE,\left( \rho ,\eta \right) \tau
_{E},E\right) ^{2} & ^{\underrightarrow{~\ \ \left[ ,\right] _{\left( \rho
,\eta \right) TE}~\ \ }} & \Gamma \left( \left( \rho ,\eta \right) TE,\left(
\rho ,\eta \right) \tau _{E},E\right)%
\end{array}%
\end{equation*}%
\emph{defined by}%
\begin{equation*}
\begin{array}{l}
\left[ Z_{1}^{\alpha }\tilde{\partial}_{\alpha }+Y_{1}^{a}\overset{\cdot }{%
\tilde{\partial}}_{a},Z_{2}^{\beta }\tilde{\partial}_{\beta }+Y_{2}^{b}%
\overset{\cdot }{\tilde{\partial}}_{b}\right] _{\left( \rho ,\eta \right) TE}%
\vspace*{1mm} \\
\displaystyle=\left[ Z_{1}^{\alpha }T_{\alpha },Z_{2}^{\beta }T_{\beta }%
\right] _{\pi ^{\ast }\left( h^{\ast }F\right) }\oplus \left[ Z_{1}^{\alpha
}\left( \rho _{\alpha }^{i}\circ h\circ \pi \right) \partial _{i}+Y_{1}^{a}%
\dot{\partial}_{a},\right. \vspace*{1mm} \\
\hfill \displaystyle\left. Z_{2}^{\beta }\left( \rho _{\beta }^{j}\circ
h\circ \pi \right) \partial _{j}+Y_{2}^{b}\dot{\partial}_{b}\right] _{TE},%
\end{array}%
\leqno(3.8)
\end{equation*}%
\emph{for any} $Z_{1}^{\alpha }\tilde{\partial}_{\alpha }+Y_{1}^{a}\overset{%
\cdot }{\tilde{\partial}}_{a}$\emph{\ and }$Z_{2}^{\beta }\tilde{\partial}%
_{\beta }+Y_{2}^{b}\overset{\cdot }{\tilde{\partial}}_{b},$ \emph{we obtain
that the couple }$\left( \left[ ,\right] _{\left( \rho ,\eta \right)
TE},\left( \tilde{\rho},Id_{E}\right) \right) $\emph{\ is a Lie algebroid
structure for the vector bundle }$\left( \left( \rho ,\eta \right) TE,\left(
\rho ,\eta \right) \tau _{E},E\right) .$

\emph{Remark 3.2}\textbf{\ }In particular, if $h=Id_{M},$ then the Lie
algebroid
\begin{equation*}
\begin{array}{c}
\left( \left( \left( Id_{TM},Id_{M}\right) TE,\left( Id_{TM},Id_{M}\right)
\tau _{E},E\right) ,\left[ ,\right] _{\left( Id_{TM},Id_{M}\right)
TE},\left( \widetilde{Id_{TM}},Id_{E}\right) \right)%
\end{array}%
\end{equation*}%
is isomorphic with the usual Lie algebroid
\begin{equation*}
\begin{array}{c}
\left( \left( TE,\tau _{E},E\right) ,\left[ ,\right] _{TE},\left(
Id_{TE},Id_{E}\right) \right) .%
\end{array}%
\end{equation*}

This is a reason for which the Lie algebroid
\begin{equation*}
\begin{array}{c}
\left( \left( \left( \rho ,\eta \right) TE,\left( \rho ,\eta \right) \tau
_{E},E\right) ,\left[ ,\right] _{\left( \rho ,\eta \right) TE},\left( \tilde{%
\rho},Id_{E}\right) \right)%
\end{array}%
,
\end{equation*}%
will be called the \emph{Lie algebroid generalized tangent bundle.}

The vector bundle $\left( \left( \rho ,\eta \right) TE,\left( \rho ,\eta
\right) \tau _{E},E\right) $ will be called the \emph{generalized tangent
bundle.}

\subsection{The generalized tangent bundle of dual vector bundle}

Let $\left( E,\pi ,M\right) $ be a vector bundle. We build the generalized
tangent bundle of dual vector bundle $\left( \overset{\ast }{E},\overset{%
\ast }{\pi },M\right) $ using the diagram:%
\begin{equation*}
\begin{array}{rcl}
\overset{\ast }{E} &  & \left( F,\left[ ,\right] _{F,h},\left( \rho ,\eta
\right) \right) \\
\overset{\ast }{\pi }\downarrow &  & ~\downarrow \nu \\
M & ^{\underrightarrow{~\ \ \ \ h~\ \ \ \ }} & ~\ N%
\end{array}%
,\leqno(3.1.1)
\end{equation*}%
where $\left( \left( F,\nu ,N\right) ,\left[ ,\right] _{F,h},\left( \rho
,\eta \right) \right) $ is a generalized Lie algebroid.

We take $\left( x^{i},p_{a}\right) $ as canonical local coordinates on $%
\left( \overset{\ast }{E},\overset{\ast }{\pi },M\right) ,$ where $i\in
\overline{1,m}$ and $a\in \overline{1,r}.$

Consider
\begin{equation*}
\left( x^{i},p_{a}\right) \longrightarrow \left( x^{i%
{\acute{}}%
}\left( x^{i}\right) ,p_{a%
{\acute{}}%
}\left( x^{i},p_{a}\right) \right)
\end{equation*}%
a change of coordinates on $\left( \overset{\ast }{E},\overset{\ast }{\pi }%
,M\right) $. Then the coordinates $p_{a}$ change to $p_{a%
{\acute{}}%
}$ according to the rule:
\begin{equation*}
\begin{array}{c}
p_{a%
{\acute{}}%
}=M_{a%
{\acute{}}%
}^{a}p_{a}.%
\end{array}%
\leqno(3.1.2)
\end{equation*}

Let
\begin{equation*}
\begin{array}[b]{c}
\left( \overset{\ast }{\partial }_{i},\dot{\partial}^{a}\right) \overset{put}%
{=}\left( \frac{\partial }{\partial x^{i}},\frac{\partial }{\partial p_{a}}%
\right)%
\end{array}%
\end{equation*}
be the base sections for the Lie $\mathcal{F}\left( \overset{\ast }{E}%
\right) $-algebra
\begin{equation*}
\left( \Gamma \left( T\overset{\ast }{E},\tau _{\overset{\ast }{E}},\overset{%
\ast }{E}\right) ,+,\cdot ,\left[ ,\right] _{T\overset{\ast }{E}}\right) .
\end{equation*}

For any sections%
\begin{equation*}
\begin{array}[b]{c}
Z^{\alpha }T_{\alpha }\in \Gamma \left( \overset{\ast }{\pi }^{\ast }\left(
h^{\ast }F\right) ,\overset{\ast }{\pi }^{\ast }\left( h^{\ast }\nu \right) ,%
\overset{\ast }{E}\right)%
\end{array}%
\end{equation*}%
and%
\begin{equation*}
\begin{array}[b]{c}
Y_{a}\dot{\partial}^{a}\in \Gamma \left( VT\overset{\ast }{E},\tau _{\overset%
{\ast }{E}},\overset{\ast }{E}\right) ,%
\end{array}%
\end{equation*}%
we obtain the section%
\begin{equation*}
\begin{array}{l}
Z^{\alpha }\overset{\ast }{\tilde{\partial}}_{\alpha }+Y_{a}\overset{\cdot }{%
\tilde{\partial}}^{a}=:Z^{\alpha }\left( T_{\alpha }\oplus \left( \rho
_{\alpha }^{i}\circ h\circ \overset{\ast }{\pi }\right) \overset{\ast }{%
\partial }_{i}\right) +Y_{a}\left( 0_{\overset{\ast }{\pi }^{^{\ast }}\left(
h^{\ast }F\right) }\oplus \dot{\partial}^{a}\right) \vspace*{1mm} \\
\qquad \displaystyle=Z^{\alpha }T_{\alpha }\oplus \left( Z^{\alpha }\left(
\rho _{\alpha }^{i}\circ h\circ \overset{\ast }{\pi }\right) \overset{\ast }{%
\partial }_{i}+Y_{a}\dot{\partial}^{a}\right) \in \Gamma \left( \overset{%
\ast }{\pi }^{\ast }\left( h^{\ast }F\right) \oplus T\overset{\ast }{E},%
\overset{\oplus }{\overset{\ast }{\pi }},\overset{\ast }{E}\right) .%
\end{array}%
\end{equation*}

Since we have
\begin{equation*}
\begin{array}{c}
Z^{\alpha }\overset{\ast }{\tilde{\partial}}_{\alpha }+Y_{a}\overset{\cdot }{%
\tilde{\partial}}^{a}=0_{\overset{\ast }{\pi }^{\ast }\left( h^{\ast
}F\right) \oplus T\overset{\ast }{E}}\vspace*{1mm} \\
\Updownarrow \\
Z^{\alpha }T_{\alpha }=0_{\overset{\ast }{\pi }^{\ast }\left( h^{\ast
}F\right) }~\wedge Z^{\alpha }\left( \rho _{\alpha }^{i}\circ h\circ \overset%
{\ast }{\pi }\right) \overset{\ast }{\partial }_{i}+Y_{a}\dot{\partial}%
^{a}=0_{T\overset{\ast }{E}},%
\end{array}%
\end{equation*}%
it implies $Z^{\alpha }=0,~\alpha \in \overline{1,p}$ and $Y_{a}=0,~a\in
\overline{1,r}.$

Therefore, the sections
\begin{equation*}
\begin{array}{c}
\displaystyle\overset{\ast }{\tilde{\partial}}_{1},...,\overset{\ast }{%
\tilde{\partial}}_{p},\overset{\cdot }{\tilde{\partial}}^{1},...,\overset{%
\cdot }{\tilde{\partial}}^{r}%
\end{array}%
\end{equation*}%
are linearly independent.

We consider the vector subbundle
\begin{equation*}
\begin{array}{c}
\left( \left( \rho ,\eta \right) T\overset{\ast }{E},\left( \rho ,\eta
\right) \tau _{\overset{\ast }{E}},\overset{\ast }{E}\right)%
\end{array}%
\end{equation*}%
of vector bundle
\begin{equation*}
\begin{array}[b]{c}
\left( \overset{\ast }{\pi }^{\ast }\left( h^{\ast }F\right) \oplus T\overset%
{\ast }{E},\overset{\oplus }{\overset{\ast }{\pi }},\overset{\ast }{E}%
\right) ,%
\end{array}%
\end{equation*}%
for which the $\mathcal{F}\left( \overset{\ast }{E}\right) $-module of
sections is the $\mathcal{F}\left( \overset{\ast }{E}\right) $-submodule of
\begin{equation*}
\left( \Gamma \left( \overset{\ast }{\pi }^{\ast }\left( h^{\ast }F\right)
\oplus T\overset{\ast }{E},\overset{\oplus }{\overset{\ast }{\pi }},\overset{%
\ast }{E}\right) ,+,\cdot \right) ,
\end{equation*}%
generated by the family of sections $\left( \overset{\ast }{\tilde{\partial}}%
_{\alpha },\overset{\cdot }{\tilde{\partial}}^{a}\right) $ which is called
the \emph{natural }$\left( \rho ,\eta \right) $\emph{-base.}

The matrix of coordinate transformation on $\left( \left( \rho ,\eta \right)
T\overset{\ast }{E},\left( \rho ,\eta \right) \tau _{\overset{\ast }{E}},%
\overset{\ast }{E}\right) $ at a change of fibred charts is
\begin{equation*}
\left\Vert
\begin{array}{cc}
\Lambda _{\alpha }^{\alpha
{\acute{}}%
}\circ h\circ \overset{\ast }{\pi } & 0 \\
\left( \rho _{\alpha }^{i}\circ h\circ \overset{\ast }{\pi }\right) %
\displaystyle\frac{\partial M_{a%
{\acute{}}%
}^{b}\circ \overset{\ast }{\pi }}{\partial x_{i}}p_{b} & M_{a%
{\acute{}}%
}^{a}\circ \overset{\ast }{\pi }%
\end{array}%
\right\Vert .\leqno(3.1.3)
\end{equation*}

We consider the operation $\left[ ,\right] _{\left( \rho ,\eta \right) T%
\overset{\ast }{E}}$ defined by%
\begin{equation*}
\begin{array}{l}
\left[ Z_{1}^{\alpha }\overset{\ast }{\tilde{\partial}}_{\alpha }+Y_{a}^{1}%
\overset{\cdot }{\tilde{\partial}}^{a},Z_{2}^{\beta }\overset{\ast }{\tilde{%
\partial}}_{\beta }+Y_{b}^{2}\overset{\cdot }{\tilde{\partial}}^{b}\right]
_{\left( \rho ,\eta \right) T\overset{\ast }{E}}=\vspace*{1mm} \\
=\left[ Z_{1}^{\alpha }T_{\alpha },Z_{2}^{\beta }T_{\beta }\right] _{\overset%
{\ast }{\pi }^{\ast }\left( h^{\ast }F\right) }\oplus \left[ Z_{1}^{\alpha
}\left( \rho _{\alpha }^{i}\circ h\circ \overset{\ast }{\pi }\right) \overset%
{\ast }{\partial }_{i}+Y_{a}^{1}\dot{\partial}^{a},\right. \\
\hfill \left. Z_{2}^{\beta }\left( \rho _{\beta }^{j}\circ h\circ \overset{%
\ast }{\pi }\right) \overset{\ast }{\partial }_{j}+Y_{b}^{2}\dot{\partial}%
^{b}\right] _{T\overset{\ast }{E}},%
\end{array}%
\leqno(3.1.4)
\end{equation*}%
for any sections $Z_{1}^{\alpha }\overset{\ast }{\tilde{\partial}}_{\alpha
}+Y_{a}^{1}\overset{\cdot }{\tilde{\partial}}^{a}$ and $Z_{2}^{\beta }%
\overset{\ast }{\tilde{\partial}}_{\beta }+Y_{b}^{2}\overset{\cdot }{\tilde{%
\partial}}^{b}.$

Let $\left( \overset{\ast }{\tilde{\rho}},Id_{\overset{\ast }{E}}\right) $\
be the $\mathbf{B}^{\mathbf{v}}$-morphism of $\left( \left( \rho ,\eta
\right) T\overset{\ast }{E},\left( \rho ,\eta \right) \tau _{\overset{\ast }{%
E}},\overset{\ast }{E}\right) $\ source and $\left( T\overset{\ast }{E},\tau
_{\overset{\ast }{E}},\overset{\ast }{E}\right) $\ target, where
\begin{equation*}
\begin{array}{rcl}
\left( \rho ,\eta \right) T\overset{\ast }{E}\!\!\! & \!\!^{\underrightarrow{%
\ \,\overset{\ast }{\tilde{\rho}}\ \,}}\!\!\! & \!\!T\overset{\ast }{E} \\
\left( Z^{\alpha }\overset{\ast }{\tilde{\partial}}_{\alpha }+Y_{a}\overset{%
\cdot }{\tilde{\partial}}^{a}\right) \!(\overset{\ast }{u}_{x})\!\!\!\! &
\!\!\longmapsto \!\!\! & \!\!\left( \!Z^{\alpha }\!\left( \rho _{\alpha }^{i}%
{\circ }h{\circ }\overset{\ast }{\pi }\!\right) \!\overset{\ast }{\partial }%
_{i}{+Y}_{a}\dot{\partial}^{a}\right) \!(\overset{\ast }{u}_{x})\!\!.%
\end{array}%
\bigskip \leqno(3.1.5)
\end{equation*}

The Lie algebroid generalized tangent bundle of the dual vector bundle $%
\left( \overset{\ast }{E},\overset{\ast }{\pi },M\right) $ will be denoted
\begin{equation*}
\begin{array}{c}
\left( \left( \left( \rho ,\eta \right) T\overset{\ast }{E},\left( \rho
,\eta \right) \tau _{\overset{\ast }{E}},\overset{\ast }{E}\right) ,\left[ ,%
\right] _{\left( \rho ,\eta \right) T\overset{\ast }{E}},\left( \overset{%
\ast }{\tilde{\rho}},Id_{\overset{\ast }{E}}\right) \right) .%
\end{array}%
\end{equation*}

\section{(Linear) $\left( \protect\rho ,\protect\eta \right) $-connections}

We consider the diagram:
\begin{equation*}
\begin{array}{c}
\xymatrix{E\ar[d]_\pi&\left( F,\left[ ,\right] _{F,h},\left( \rho ,\eta
\right) \right)\ar[d]^\nu\\ M\ar[r]^h&N}%
\end{array}%
\end{equation*}%
where $\left( E,\pi ,M\right) \in \left\vert \mathbf{B}\right\vert $ and $%
\left( \left( F,\nu ,N\right) ,\left[ ,\right] _{F,h},\left( \rho ,\eta
\right) \right) $ is a generalized Lie algebroid.

Let
\begin{equation*}
\left( \left( \left( \rho ,\eta \right) TE,\left( \rho ,\eta \right) \tau
_{E},E\right) ,\left[ ,\right] _{\left( \rho ,\eta \right) TE},\left( \tilde{%
\rho},Id_{E}\right) \right)
\end{equation*}%
\ be the Lie algebroid generalized tangent bundle of the fiber bundle $%
\left( E,\pi ,M\right) $.

We consider the $\mathbf{B}^{\mathbf{v}}$-morphism $\left( \left( \rho ,\eta
\right) \pi !,Id_{E}\right) $ given by the commutative diagram%
\begin{equation*}
\begin{array}{c}
\xymatrix{\left( \rho ,\eta \right) TE\ar[r]^{( \rho ,\eta ) \pi
!}\ar[d]_{(\rho,\eta)\tau_E}&\pi ^{\ast }\left( h^{\ast
}F\right)\ar[d]^{pr_1} \\ E\ar[r]^{id_{E}}& E}%
\end{array}%
\leqno(4.1)
\end{equation*}

This is defined as:%
\begin{equation*}
\begin{array}{c}
\left( \rho ,\eta \right) \pi !\left( \left( Z^{\alpha }\overset{\ast }{%
\tilde{\partial}}_{\alpha }+Y^{a}\overset{\cdot }{\tilde{\partial}}%
_{a}\right) \left( u_{x}\right) \right) =\left( Z^{\alpha }T_{\alpha
}\right) \left( u_{x}\right) ,%
\end{array}%
\leqno(4.2)
\end{equation*}%
for any $Z^{\alpha }\overset{\ast }{\tilde{\partial}}_{\alpha }+Y^{a}\overset%
{\cdot }{\tilde{\partial}}_{a}\in \Gamma \left( \left( \rho ,\eta \right)
TE,\left( \rho ,\eta \right) \tau _{E},E\right) .$\medskip

Using the $\mathbf{B}^{\mathbf{v}}$-morphism $\left( \left( \rho ,\eta
\right) \pi !,Id_{E}\right) ,$ and the the $\mathbf{B}^{\mathbf{v}}$%
-morphism $\left( 2.7\right) $ we obtain the \emph{tangent }$\left( \rho
,\eta \right) $\emph{-application }$\left( \left( \rho ,\eta \right) T\pi
,h\circ \pi \right) $ of $\left( \left( \rho ,\eta \right) TE,\left( \rho
,\eta \right) \tau _{E},E\right) $ source and $\left( F,\nu ,N\right) $
target.

\textbf{Definition 4.1} The kernel of the tangent $\left( \rho ,\eta \right)
$-application\ is writen
\begin{equation*}
\left( V\left( \rho ,\eta \right) TE,\left( \rho ,\eta \right) \tau
_{E},E\right)
\end{equation*}
and it is called \emph{the vertical subbundle}.\bigskip

We remark that the set $\left\{ \overset{\cdot }{\tilde{\partial}}_{a},~a\in
\overline{1,r}\right\} $ is a base of the $\mathcal{F}\left( E\right) $%
-module
\begin{equation*}
\left( \Gamma \left( V\left( \rho ,\eta \right) TE,\left( \rho ,\eta \right)
\tau _{E},E\right) ,+,\cdot \right) .
\end{equation*}

\textbf{Proposition 4.1} \emph{The short sequence of vector bundles}%
\begin{equation*}
\begin{array}{c}
\xymatrix{0\ar@{^(->}[r]^i\ar[d]&V(\rho,\eta)TE\ar[d]\ar@{^(->}[r]^i&(\rho,%
\eta)TE\ar[r]^{(\rho,\eta)\pi!}\ar[d]&\pi ^{\ast }\left( h^{\ast
}F\right)\ar[r]\ar[d]&0\ar[d]\\
E\ar[r]^{Id_E}&E\ar[r]^{Id_E}&E\ar[r]^{Id_E}&E\ar[r]^{Id_E}&E}%
\end{array}%
\leqno(4.3)
\end{equation*}%
\emph{is exact.}

\textbf{Definition 4.2} \textit{A }$\mathbf{Man}$-morphism $\left( \rho
,\eta \right) \Gamma $ of $\left( \rho ,\eta \right) TE$ source and $V\left(
\rho ,\eta \right) TE$ target defined by%
\begin{equation*}
\begin{array}{c}
\left( \rho ,\eta \right) \Gamma \left( Z^{\gamma }\overset{\ast }{\tilde{%
\partial}}_{\gamma }+Y^{a}\overset{\cdot }{\tilde{\partial}}_{a}\right)
\left( u_{x}\right) =\left( Y^{a}+\left( \rho ,\eta \right) \Gamma _{\gamma
}^{a}Z^{\gamma }\right) \overset{\cdot }{\tilde{\partial}}_{a}\left(
u_{x}\right) ,%
\end{array}%
\leqno(4.4)
\end{equation*}%
so that the $\mathbf{B}^{\mathbf{v}}$-morphism $\left( \left( \rho ,\eta
\right) \Gamma ,Id_{E}\right) $ is a split to the left in the previous exact
sequence, will be called $\left( \rho ,\eta \right) $\emph{-connection for
the fiber bundle }$\left( E,\pi ,M\right) $.

The $\left( \rho ,Id_{M}\right) $-connection will be called $\rho $\emph{%
-connection }and will be denoted $\rho \Gamma $\emph{\ }and the $\left(
Id_{TM},Id_{M}\right) $-connection will be called \emph{connection }and will
be denoted $\Gamma $\emph{.}

\textbf{Definition 4.3 }If $\left( \rho ,\eta \right) \Gamma $ is a $\left(
\rho ,\eta \right) $-connection for the fiber bundle $\left( E,\pi ,M\right)
$, then the kernel of the $\mathbf{B}^{\mathbf{v}}$-morphism $\left( \left(
\rho ,\eta \right) \Gamma ,Id_{E}\right) $\ is written $\left( H\left( \rho
,\eta \right) TE,\left( \rho ,\eta \right) \tau _{E},E\right) $ and will be
called the \emph{horizontal vector subbundle}.

\textbf{Definition 4.4} If $\left( E,\pi ,M\right) \in \left\vert \mathbf{B}%
\right\vert $, then the $\mathbf{B}$-morphism $\left( \Pi ,\pi \right) $
defined by the commutative diagram%
\begin{equation*}
\begin{array}{c}
\xymatrix{ {V\left( \rho ,\eta \right)T}{E\ar[r]^{\qquad \Pi}}
\ar[d]_{(\rho,\eta)\tau_E} &E \ar[d]^\pi \\ E\ar[r]^\pi&M}%
\end{array}%
\leqno(4.5)
\end{equation*}%
so that the components of the image of the vector $Y^{a}\overset{\cdot }{%
\tilde{\partial}}_{a}\left( u_{x}\right) $ are the real numbers $Y^{1}\left(
u_{x}\right) ,...,Y^{r}\left( u_{x}\right) $ will be called the \emph{%
canonical projection }$\mathbf{B}$\emph{-morphism.}\medskip

In particular, if $\left( E,\pi ,M\right) \in \left\vert \mathbf{B}^{\mathbf{%
v}}\right\vert $ and $\left\{ s_{a},a\in \overline{1,r}\right\} $ is a base
of the $\mathcal{F}\left( M\right) $-module of sections $\left( \Gamma
\left( E,\pi ,M\right) ,+,\cdot \right) $, then $\Pi $\textit{\ }is defined
by\textit{\ }%
\begin{equation*}
\begin{array}{c}
\Pi \left( Y^{a}\overset{\cdot }{\tilde{\partial}}_{a}\left( u_{x}\right)
\right) =Y^{a}\left( u_{x}\right) s_{a}\left( \pi \left( u_{x}\right)
\right) =Y^{a}\left( u_{x}\right) s_{a}\left( x\right) .%
\end{array}%
\leqno(4.6)
\end{equation*}

\textbf{Theorem 4.1 }\emph{If }$\left( \rho ,\eta \right) \Gamma $\emph{\ is
a }$\left( \rho ,\eta \right) $\emph{-connection for the fiber bundle }$%
\left( E,\pi ,M\right) ,$\emph{\ then its components satisfy the law of
transformation }%
\begin{equation*}
\begin{array}{c}
\left( \rho ,\eta \right) \Gamma _{\gamma
{\acute{}}%
}^{a%
{\acute{}}%
}=\frac{\partial y^{a%
{\acute{}}%
}}{\partial y^{a}}\left[ \rho _{\gamma }^{k}\circ \left( h\circ \pi \right)
\frac{\partial y^{a}}{\partial x^{k}}+\left( \rho ,\eta \right) \Gamma
_{\gamma }^{a}\right] \Lambda _{\gamma
{\acute{}}%
}^{\gamma }\circ \left( h\circ \pi \right) .%
\end{array}%
\leqno(4.7)
\end{equation*}

\emph{If }$\left( \rho ,\eta \right) \Gamma $\emph{\ is a }$\left( \rho
,\eta \right) $\emph{-connection for the vector bundle }$\left( E,\pi
,M\right) ,$\emph{\ then its components satisfy the law of transformation }%
\begin{equation*}
\begin{array}{c}
\left( \rho ,\eta \right) \Gamma _{\gamma
{\acute{}}%
}^{a%
{\acute{}}%
}{=}M_{a}^{a%
{\acute{}}%
}{\circ }\pi \!\!\left[ \rho _{\gamma }^{k}{\circ }\left( h{\circ }\pi
\right) \!\frac{\partial M_{b%
{\acute{}}%
}^{a}\circ \pi }{\partial x^{k}}y^{b%
{\acute{}}%
}{+}\left( \rho ,\eta \right) \!\Gamma _{\gamma }^{a}\right] \!\!\Lambda
_{\gamma
{\acute{}}%
}^{\gamma }{\circ }\left( h{\circ }\pi \right) .%
\end{array}%
\leqno(4.7^{\prime })
\end{equation*}

\emph{In the particular case of Lie algebroids, }$\left( \eta ,h\right)
=\left( Id_{M},Id_{M}\right) ,$\emph{\ the relations }$\left( 4.7^{\prime
}\right) $\emph{\ become}%
\begin{equation*}
\begin{array}{c}
\rho \Gamma _{\gamma
{\acute{}}%
}^{a%
{\acute{}}%
}=M_{a}^{a%
{\acute{}}%
}\circ \pi \left[ \rho _{\gamma }^{k}\circ \pi \frac{\partial M_{b%
{\acute{}}%
}^{a}\circ \pi }{\partial x^{k}}y^{b%
{\acute{}}%
}+\rho \Gamma _{\gamma }^{a}\right] \Lambda _{\gamma
{\acute{}}%
}^{\gamma }\circ \pi .%
\end{array}%
\leqno(4.7^{\prime \prime })
\end{equation*}

\emph{In the classical case, }$\left( \rho ,\eta ,h\right) =\left(
Id_{TM},Id_{M},Id_{M}\right) ,$\emph{\ the relations }$\left( 4.7^{\prime
\prime }\right) $\emph{\ become}%
\begin{equation*}
\begin{array}{c}
\Gamma _{k%
{\acute{}}%
}^{i%
{\acute{}}%
}=\frac{\partial x^{i%
{\acute{}}%
}}{\partial x^{i}}\circ \tau _{M}\left[ \frac{\partial }{\partial x^{k}}%
\left( \frac{\partial x^{i}}{\partial x^{j%
{\acute{}}%
}}\circ \tau _{M}\right) y^{j%
{\acute{}}%
}+\Gamma _{k}^{i}\right] \frac{\partial x^{k}}{\partial x^{k%
{\acute{}}%
}}\circ \tau _{M}.%
\end{array}%
\leqno(4.7^{\prime \prime \prime })
\end{equation*}

\emph{Proof.} Let $\left( \Pi ,\pi \right) $ be the canonical projection $%
\mathbf{B}$-morphism.

Obviously, the components of
\begin{equation*}
\begin{array}{c}
\Pi \circ \left( \rho ,\eta \right) \Gamma \left( Z^{\gamma ^{\prime }}%
\overset{\ast }{\tilde{\partial}}_{\gamma ^{\prime }}+Y^{a^{\prime }}\overset%
{\cdot }{\tilde{\partial}}_{a^{\prime }}\right) \left( u_{x}\right)%
\end{array}%
\end{equation*}%
are the real numbers%
\begin{equation*}
\left( Y^{a%
{\acute{}}%
}+\left( \rho ,\eta \right) \Gamma _{\gamma
{\acute{}}%
}^{a%
{\acute{}}%
}Z^{\gamma
{\acute{}}%
}\right) \left( u_{x}\right) .
\end{equation*}

Since
\begin{equation*}
\begin{array}{l}
\displaystyle\left( Z^{\gamma ^{\prime }}\overset{\ast }{\tilde{\partial}}%
_{\gamma ^{\prime }}+Y^{a^{\prime }}\overset{\cdot }{\tilde{\partial}}%
_{a^{\prime }}\right) \left( u_{x}\right) =Z^{\gamma
{\acute{}}%
}\Lambda _{\gamma
{\acute{}}%
}^{\gamma }\circ h\circ \pi \overset{\ast }{\tilde{\partial}}_{\gamma
}\left( u_{x}\right) \vspace*{1mm} \\
\qquad \displaystyle+\left( Z^{\gamma
{\acute{}}%
}\rho _{\gamma
{\acute{}}%
}^{i%
{\acute{}}%
}\circ h\circ \pi \frac{\partial y^{a}}{\partial x^{i%
{\acute{}}%
}}+\frac{\partial y^{a}}{\partial y^{a%
{\acute{}}%
}}Y^{a%
{\acute{}}%
}\right) \overset{\cdot }{\tilde{\partial}}_{a}\left( u_{x}\right) ,%
\end{array}%
\end{equation*}%
it results that the components of%
\begin{equation*}
\begin{array}{c}
\Pi \circ \left( \rho ,\eta \right) \Gamma \left( Z^{\gamma ^{\prime }}%
\overset{\ast }{\tilde{\partial}}_{\gamma ^{\prime }}+Y^{a^{\prime }}\overset%
{\cdot }{\tilde{\partial}}_{a^{\prime }}\right) \left( u_{x}\right)%
\end{array}%
\end{equation*}%
are the real numbers
\begin{equation*}
\begin{array}{c}
\left( Z^{\gamma
{\acute{}}%
}\rho _{\gamma
{\acute{}}%
}^{i%
{\acute{}}%
}\circ h\circ \pi \frac{\partial y^{a}}{\partial x^{i%
{\acute{}}%
}}+\frac{\partial y^{a}}{\partial y^{a%
{\acute{}}%
}}Y^{a%
{\acute{}}%
}+\left( \rho ,\eta \right) \Gamma _{\gamma }^{a}Z^{\gamma
{\acute{}}%
}\Lambda _{\gamma
{\acute{}}%
}^{\gamma }\circ h\circ \pi \right) \left( u_{x}\right) \frac{\partial y^{a%
{\acute{}}%
}}{\partial y^{a}},%
\end{array}%
\end{equation*}%
where
\begin{equation*}
\begin{array}{c}
\left\Vert \frac{\partial y^{a}}{\partial y^{a%
{\acute{}}%
}}\right\Vert =\left\Vert \frac{\partial y^{a%
{\acute{}}%
}}{\partial y^{a}}\right\Vert ^{-1}.%
\end{array}%
\end{equation*}

Therefore, we have:%
\begin{equation*}
\begin{array}{c}
\left( Z^{\gamma
{\acute{}}%
}\rho _{\gamma
{\acute{}}%
}^{i%
{\acute{}}%
}\circ h\circ \pi \frac{\partial y^{a}}{\partial x^{i%
{\acute{}}%
}}+\frac{\partial y^{a}}{\partial y^{a%
{\acute{}}%
}}Y^{a%
{\acute{}}%
}+\left( \rho ,\eta \right) \Gamma _{\gamma }^{a}Z^{\gamma
{\acute{}}%
}\Lambda _{v%
{\acute{}}%
}^{\gamma }\circ h\circ \pi \right) \frac{\partial y^{a%
{\acute{}}%
}}{\partial y^{a}}=Y^{a%
{\acute{}}%
}+\left( \rho ,\eta \right) \Gamma _{\gamma
{\acute{}}%
}^{a%
{\acute{}}%
}Z^{\gamma
{\acute{}}%
}.%
\end{array}%
\end{equation*}

After some calculations we obtain:
\begin{equation*}
\begin{array}{c}
\left( \rho ,\eta \right) \Gamma _{\gamma
{\acute{}}%
}^{a%
{\acute{}}%
}=\frac{\partial y^{a%
{\acute{}}%
}}{\partial y^{a}}\left( \rho _{\gamma }^{i}\circ \left( h\circ \pi \right)
\frac{\partial y^{a}}{\partial x^{i}}+\left( \rho ,\eta \right) \Gamma
_{\gamma }^{a}\right) \Lambda _{\gamma
{\acute{}}%
}^{\gamma }\circ h\circ \pi .%
\end{array}%
\eqno{{q.e.d.}}
\end{equation*}

\emph{Remark 4.1} If $\Gamma $ is a Ehresmann connection for the vector
bundle $\left( E,\pi ,M\right) $ on components $\Gamma _{k}^{a},$ then the
differentiable real local functions $\left( \rho ,\eta \right) \Gamma
_{\gamma }^{a}=\left( \rho _{\gamma }^{k}\circ h\circ \pi \right) \Gamma
_{k}^{a}$ are the components of a $\left( \rho ,\eta \right) $-connection $%
\left( \rho ,\eta \right) \Gamma $ for the vector bundle $\left( E,\pi
,M\right) .$ This $\left( \rho ,\eta \right) $-connection will be called the
$\left( \rho ,\eta \right) $\emph{-connection associated to the connection }$%
\Gamma .$

\textbf{Definition 4.5 }If $\left( \rho ,\eta \right) \Gamma $ is a $\left(
\rho ,\eta \right) $-connection for the vector bundle $\left( E,\pi
,M\right) $ and $z=z^{\gamma }t_{\gamma }\in \Gamma \left( F,\nu ,M\right) ,$
then the application%
\begin{equation*}
\begin{array}{rcl}
\Gamma \left( E,\pi ,M\right) & ^{\underrightarrow{\left( \rho ,\eta \right)
D_{z}}} & \Gamma \left( E,\pi ,M\right) \\
u=u^{a}s_{a} & \longmapsto & \left( \rho ,\eta \right) D_{z}u%
\end{array}%
\end{equation*}%
where
\begin{equation*}
\begin{array}[b]{c}
\left( \rho ,\eta \right) D_{z}u=z^{\gamma }\circ h\left( \rho _{\gamma
}^{k}\circ h\frac{\partial u^{a}}{\partial x^{k}}+\left( \rho ,\eta \right)
\Gamma _{\gamma }^{a}\circ u\right) s_{a}%
\end{array}%
\leqno(4.8)
\end{equation*}%
will be called the \emph{covariant }$\left( \rho ,\eta \right) $\emph{%
-derivative associated to }$\left( \rho ,\eta \right) $\emph{-connection }$%
\left( \rho ,\eta \right) \Gamma $\emph{\ with respect to the section }$z$.

In the particular case of Lie algebroids, $\left( \eta ,h\right) =\left(
Id_{M},Id_{M}\right) ,$\ the relations $\left( 4.8\right) $\ become%
\begin{equation*}
\!\!%
\begin{array}[b]{c}
\rho D_{z}u=z^{\gamma }\left( \rho _{\gamma }^{k}\frac{\partial u^{a}}{%
\partial x^{k}}+\rho \Gamma _{\gamma }^{a}\circ u\right) s_{a}.%
\end{array}%
\leqno(4.8^{\prime })
\end{equation*}

In the classical case, $\left( \rho ,\eta ,h\right) =\left(
Id_{TM},Id_{M},Id_{M}\right) ,$\ the relations $\left( 4.8^{\prime }\right) $%
\ become%
\begin{equation*}
\!\!%
\begin{array}[b]{c}
D_{X}Y{=}X^{k}\left( \frac{\partial Y^{i}}{\partial x^{k}}{+}\Gamma
_{k}^{i}\circ Y\right) \partial _{i}.%
\end{array}%
\leqno(4.8^{\prime \prime })
\end{equation*}

\textbf{Definition 4.6 }Let $\left( \rho ,\eta \right) \Gamma $ be a $\left(
\rho ,\eta \right) $-connection for the fiber bundle $\left( E,\pi ,M\right)
.$ If for each local vector $\left( m+r\right) $-chart $\left(
U,s_{U}\right) $ and for each local vector $\left( n+p\right) $-chart $%
\left( V,t_{V}\right) $ so that $U\cap h^{-1}\left( V\right) \neq \phi $, it
exists the differentiable real functions $\left( \rho ,\eta \right) \Gamma
_{b\gamma }^{a}$ defined on $U\cap h^{-1}\left( V\right) $ such that
\begin{equation*}
\begin{array}{c}
\left( \rho ,\eta \right) \Gamma _{\gamma }^{a}\circ u=\left( \rho ,\eta
\right) \Gamma _{b\gamma }^{a}\cdot u^{b},\forall u=u^{b}s_{b}\in \Gamma
\left( E,\pi ,M\right) ,%
\end{array}%
\leqno(4.9)
\end{equation*}%
then we say that $\left( \rho ,\eta \right) \Gamma $\emph{\ is linear.}

The differentiable real local functions $\left( \rho ,\eta \right) \Gamma
_{b\gamma }^{a}$ will be called the \emph{Christoffel coefficients of linear
}$\left( \rho ,\eta \right) $\emph{-connection }$\left( \rho ,\eta \right)
\Gamma .$

\textbf{Proposition 4.2} \emph{If }$\left( \rho ,\eta \right) \Gamma $\emph{%
\ is a linear }$\left( \rho ,\eta \right) $\emph{-connection for the fiber
bundle }$\left( E,\pi ,M\right) ,$\emph{\ then its components satisfy the
law of transformation }%
\begin{equation*}
\begin{array}{c}
(\rho ,\eta )\Gamma _{b%
{\acute{}}%
\gamma
{\acute{}}%
}^{a%
{\acute{}}%
}{=}\frac{\partial y^{a%
{\acute{}}%
}}{\partial y^{a}}\left[ \rho _{\gamma }^{k}{\circ }h\!\frac{\partial }{%
\partial x^{k}}\!\left( \frac{\partial y^{a}}{\partial y^{b%
{\acute{}}%
}}\right) {+}\left( \rho ,\eta \right) \Gamma _{b\gamma }^{a}\frac{\partial
y^{b}}{\partial y^{b%
{\acute{}}%
}}\right] \!\Lambda _{\gamma
{\acute{}}%
}^{\gamma }{\circ }h.%
\end{array}%
\leqno(4.10)
\end{equation*}

\emph{If }$\left( \rho ,\eta \right) \Gamma $\emph{\ is a linear }$\left(
\rho ,\eta \right) $\emph{-connection for the vector bundle }$\left( E,\pi
,M\right) ,$\emph{\ then its components satisfy the law of transformation }%
\begin{equation*}
\begin{array}{c}
(\rho ,\eta )\Gamma _{b%
{\acute{}}%
\gamma
{\acute{}}%
}^{a%
{\acute{}}%
}{=}M_{a}^{a%
{\acute{}}%
}\!\left[ \rho _{\gamma }^{k}{\circ }h\frac{\partial M_{b%
{\acute{}}%
}^{a}}{\partial x^{k}}{+}(\rho ,\eta )\Gamma _{b\gamma }^{a}M_{b%
{\acute{}}%
}^{b}\right] \!\!\Lambda _{\gamma
{\acute{}}%
}^{\gamma }{\circ }h.%
\end{array}%
\leqno(4.10^{\prime })
\end{equation*}

\emph{In the particular case of Lie algebroids, }$\left( \eta ,h\right)
=\left( Id_{M},Id_{M}\right) ,$\emph{\ the relations }$\left( 4.10^{\prime
}\right) $\emph{\ become}%
\begin{equation*}
\begin{array}{c}
\rho \Gamma _{b%
{\acute{}}%
\gamma
{\acute{}}%
}^{a%
{\acute{}}%
}=M_{a}^{a%
{\acute{}}%
}\left[ \rho _{\gamma }^{k}\frac{\partial M_{b%
{\acute{}}%
}^{a}}{\partial x^{k}}+\rho \Gamma _{b\gamma }^{a}M_{b%
{\acute{}}%
}^{b}\right] \Lambda _{\gamma
{\acute{}}%
}^{\gamma }.%
\end{array}%
\leqno(4.10^{\prime \prime })
\end{equation*}

\emph{In the classical case, }$\left( \rho ,\eta ,h\right) =\left(
Id_{TM},Id_{M},Id_{M}\right) ,$\emph{\ the relations }$\left( 4.10^{\prime
\prime }\right) $\emph{\ become}%
\begin{equation*}
\begin{array}{c}
\Gamma _{j%
{\acute{}}%
k%
{\acute{}}%
}^{i%
{\acute{}}%
}=\frac{\partial x^{i%
{\acute{}}%
}}{\partial x^{i}}\left[ \frac{\partial }{\partial x^{k}}\left( \frac{%
\partial x^{i}}{\partial x^{j%
{\acute{}}%
}}\right) +\Gamma _{jk}^{i}\frac{\partial x^{j}}{\partial x^{j%
{\acute{}}%
}}\right] \frac{\partial x^{k}}{\partial x^{k%
{\acute{}}%
}}.%
\end{array}%
\leqno(4.10^{\prime \prime \prime })
\end{equation*}

\emph{Remark 4.2} If $\left( \rho ,\eta \right) \Gamma $ is a linear $\left(
\rho ,\eta \right) $-connection for the vector bundle $\left( E,\pi
,M\right) $, then, for any $z=z^{\gamma }t_{\gamma }\in \Gamma \left( F,\nu
,M\right) ,$ we obtain the \emph{covariant }$\left( \rho ,\eta \right) $%
\emph{-derivative associated to the linear }$\left( \rho ,\eta \right) $%
\emph{-connection }$\left( \rho ,\eta \right) \Gamma $\emph{\ with respect
to the section }$z$%
\begin{equation*}
\begin{array}{rcl}
\Gamma \left( E,\pi ,M\right) \!\! & \!\!^{\underrightarrow{\left( \rho
,\eta \right) D_{z}}}\!\! & \!\!\Gamma \left( E,\pi ,M\right) \\
u{=}u^{a}s_{a}\!\! & \!\!\longmapsto \!\! & \!\!(\rho ,\eta )D_{z}u%
\end{array}%
\end{equation*}%
defined by
\begin{equation*}
\!\!%
\begin{array}[b]{c}
(\rho ,\eta )D_{z}u{=}z^{\gamma }{\circ }h\left( \rho _{\gamma }^{k}{\circ }h%
\frac{\partial u^{a}}{\partial x^{k}}{+}(\rho ,\eta )\Gamma _{b\gamma
}^{a}\cdot u^{b}\right) s_{a}.%
\end{array}%
\leqno(4.11)
\end{equation*}

In the particular case of Lie algebroids, $\left( \eta ,h\right) =\left(
Id_{M},Id_{M}\right) ,$\ the relations $\left( 4.11\right) $\ become%
\begin{equation*}
\!\!%
\begin{array}[b]{c}
\rho D_{z}u{=}z^{\gamma }\left( \rho _{\gamma }^{k}\frac{\partial u^{a}}{%
\partial x^{k}}{+}\rho \Gamma _{b\gamma }^{a}\cdot u^{b}\right) s_{a}.%
\end{array}%
\leqno(4.11^{\prime })
\end{equation*}

In the classical case, $\left( \rho ,\eta ,h\right) =\left(
Id_{TM},Id_{M},Id_{M}\right) ,$\ the relations $\left( 4.11^{\prime }\right)
$\ become%
\begin{equation*}
\!\!%
\begin{array}[b]{c}
D_{X}Y{=}X^{k}\left( \frac{\partial Y^{i}}{\partial x^{k}}{+}\Gamma
_{jk}^{i}\cdot Y^{j}\right) \partial _{i}.%
\end{array}%
\leqno(4.11^{\prime \prime })
\end{equation*}

\subsection{(Linear) $\left( \protect\rho ,\protect\eta \right) $%
-connections for dual vector bundle}

Let $\left( E,\pi ,M\right) $ be a vector bundle.

We consider the following diagram:%
\begin{equation*}
\begin{array}{rcl}
\overset{\ast }{E} &  & \left( F,\left[ ,\right] _{F,h},\left( \rho ,\eta
\right) \right) \\
\overset{\ast }{\pi }\downarrow &  & ~\downarrow \nu \\
M & ^{\underrightarrow{~\ \ \ \ h~\ \ \ \ }} & ~\ N%
\end{array}%
,\leqno(4.1.1)
\end{equation*}%
where $\left( \left( F,\nu ,N\right) ,\left[ ,\right] _{F,h},\left( \rho
,\eta \right) \right) $ is a generalized Lie algebroid.

Let
\begin{equation*}
\left( \left( \left( \rho ,\eta \right) T\overset{\ast }{E},\left( \rho
,\eta \right) \tau _{\overset{\ast }{E}},\overset{\ast }{E}\right) ,\left[ ,%
\right] _{\left( \rho ,\eta \right) T\overset{\ast }{E}},\left( \overset{%
\ast }{\tilde{\rho}},Id_{\overset{\ast }{E}}\right) \right)
\end{equation*}%
\ be the Lie algebroid generalized tangent bundle of the vector bundle $%
\left( \overset{\ast }{E},\overset{\ast }{\pi },M\right) $.

We consider the $\mathbf{B}^{\mathbf{v}}$-morphism $\left( \left( \rho ,\eta
\right) \overset{\ast }{\pi }!,Id_{\overset{\ast }{E}}\right) $ given by the
commutative diagram%
\begin{equation*}
\begin{array}{c}
\xymatrix{\left( \rho ,\eta \right) T\overset{\ast }{E}\ar[r]^{( \rho ,\eta
) \overset{\ast }{\pi!} }\ar[d]_{(\rho,\eta)\tau_{\overset{\ast }{E}}}&
\overset{\ast }{\pi }^{\ast }\left( h^{\ast }F\right) \ar[d]^{pr_1} \\
\overset{\ast }{E}\ar[r]^{id_{\overset{\ast }{E}}}& \overset{\ast }{E}}%
\end{array}%
\leqno(4.1.2)
\end{equation*}

Using the components, this is defined as:%
\begin{equation*}
\begin{array}[b]{c}
\left( \rho ,\eta \right) \overset{\ast }{\pi }!\left( Z^{\alpha }\overset{%
\ast }{\tilde{\partial}}_{\alpha }+Y_{a}\overset{\cdot }{\tilde{\partial}}%
^{a}\right) \left( \overset{\ast }{u}_{x}\right) =\left( Z^{\alpha
}T_{\alpha }\right) \left( \overset{\ast }{u}_{x}\right) ,%
\end{array}%
\leqno(4.1.3)
\end{equation*}%
for any $\displaystyle Z^{\alpha }\overset{\ast }{\tilde{\partial}}_{\alpha
}+Y_{a}\overset{\cdot }{\tilde{\partial}}^{a}\in \left( \left( \rho ,\eta
\right) T\overset{\ast }{E},\left( \rho ,\eta \right) \tau _{\overset{\ast }{%
E}},\overset{\ast }{E}\right) .$\medskip

Using the $\mathbf{B}^{\mathbf{v}}$-morphism $\left( \left( \rho ,\eta
\right) \overset{\ast }{\pi }!,Id_{\overset{\ast }{E}}\right) $ and the $%
\mathbf{B}^{\mathbf{v}}$-morphism $\left( 2.7\right) $ we obtain the \emph{%
tangent }$\left( \rho ,\eta \right) $\emph{-application }$\left( \left( \rho
,\eta \right) T\overset{\ast }{\pi },h\circ \overset{\ast }{\pi }\right) $%
\emph{\ }of $\left( \left( \rho ,\eta \right) T\overset{\ast }{E},\left(
\rho ,\eta \right) \tau _{\overset{\ast }{E}},\overset{\ast }{E}\right) $
source and $\left( F,\nu ,N\right) $ target.

\textbf{Definition~4.1.1} The kernel of the tangent $\left( \rho ,\eta
\right) $-application\break
\begin{equation*}
\left( \left( \rho ,\eta \right) T\overset{\ast }{\pi },h\circ \overset{\ast
}{\pi }\right)
\end{equation*}%
\ is written as%
\begin{equation*}
\left( V\left( \rho ,\eta \right) T\overset{\ast }{E},\left( \rho ,\eta
\right) \tau _{\overset{\ast }{E}},\overset{\ast }{E}\right)
\end{equation*}%
and will be called the \emph{vertical subbundle}.\bigskip

The set $\left\{ \overset{\cdot }{\tilde{\partial}}^{a},~a\in \overline{1,r}%
\right\} $ is a base for the $\mathcal{F}\left( \overset{\ast }{E}\right) $%
-module
\begin{equation*}
\left( \Gamma \left( V\left( \rho ,\eta \right) T\overset{\ast }{E},\left(
\rho ,\eta \right) \tau _{\overset{\ast }{E}},\overset{\ast }{E}\right)
,+,\cdot \right) .
\end{equation*}

\textbf{Proposition 4.1.1} \emph{The short sequence of vector bundles }%
\begin{equation*}
\begin{array}{c}
\xymatrix{0\ar@{^(->}[r]^i\ar[d]&V(\rho,\eta)T\overset{\ast
}{E}\ar[d]\ar@{^(->}[r]^i&(\rho,\eta)T\overset{\ast
}{E}\ar[r]^{(\rho,\eta)\overset{\ast }{\pi }!}\ar[d]&\overset{\ast }{\pi
}^{\ast }\left( h^{\ast }F\right)\ar[r]\ar[d]&0\ar[d]\\ \overset{\ast
}{E}\ar[r]^{Id_{\overset{\ast }{E}}}&\overset{\ast
}{E}\ar[r]^{Id_{\overset{\ast }{E}}}&\overset{\ast
}{E}\ar[r]^{Id_{\overset{\ast }{E}}}&\overset{\ast }
{E}\ar[r]^{Id_{\overset{\ast }{E}}}&\overset{\ast }{E}}%
\end{array}%
\leqno(4.1.4)
\end{equation*}%
\emph{is exact.}

\textbf{Definition 4.1.2} \textit{A }$\mathbf{Man}$-morphism $\left( \rho
,\eta \right) \Gamma $ of $\left( \rho ,\eta \right) T\overset{\ast }{E}$
source and $V\left( \rho ,\eta \right) T\overset{\ast }{E}$ target defined by%
\begin{equation*}
\begin{array}[b]{c}
\left( \rho ,\eta \right) \Gamma \left( Z^{\gamma }\overset{\ast }{\tilde{%
\partial}}_{\gamma }+Y_{a}\overset{\cdot }{\tilde{\partial}}^{a}\right)
\left( \overset{\ast }{u}_{x}\right) =\left( Y_{b}-\left( \rho ,\eta \right)
\Gamma _{b\gamma }Z^{\gamma }\right) \overset{\cdot }{\tilde{\partial}}%
^{b}\left( \overset{\ast }{u}_{x}\right) ,%
\end{array}%
\leqno(4.1.5)
\end{equation*}%
such that the $\mathbf{B}^{\mathbf{v}}$-morphism $\left( \left( \rho ,\eta
\right) \Gamma ,Id_{\overset{\ast }{E}}\right) $ is a split to the left in
the previous exact sequence, will be called $\left( \rho ,\eta \right) $%
\emph{-connection for the dual vector bundle }$\left( \overset{\ast }{E},%
\overset{\ast }{\pi },M\right) $.

The differentiable real local functions $\left( \rho ,\eta \right) \Gamma
_{b\gamma }$ will be called the \emph{components of }$\left( \rho ,\eta
\right) $\emph{-connection }$\left( \rho ,\eta \right) \Gamma .$

The $\left( \rho ,Id_{M}\right) $-connection for the dual vector bundle $%
\left( \overset{\ast }{E},\overset{\ast }{\pi },M\right) $ will be called $%
\rho $\emph{-connection for the dual vector bundle }$\left( \overset{\ast }{E%
},\overset{\ast }{\pi },M\right) $ and will be denoted $\rho \Gamma $.

The $\left( Id_{TM},Id_{M}\right) $-connection for the dual vector bundle $%
\left( \overset{\ast }{E},\overset{\ast }{\pi },M\right) $ will be called
\emph{connection for the dual vector bundle }$\left( \overset{\ast }{E},%
\overset{\ast }{\pi },M\right) $ and will be denoted $\Gamma $.

Let $\left\{ s^{a},~a\in \overline{1,r}\right\} $ be the dual base of the
base $\left\{ s_{a},a\in \overline{1,r}\right\} .$

The $\mathbf{B}^{\mathbf{v}}$-morphism $\left( \overset{\ast }{\Pi },\overset%
{\ast }{\pi }\right) $ defined by the commutative diagram%
\begin{equation*}
\begin{array}{c}
\xymatrix{V\left( \rho ,\eta \right) T\overset{\ast }{E}\ar[r]^{\qquad
\overset{\ast }{\Pi } }\ar[d]_{(\rho,\eta)\tau_{\overset{\ast
}{E}}}&\overset{\ast }{E} \ar[d]^{\overset{\ast }{\pi }} \\ \overset{\ast
}{E}\ar[r]^{\overset{\ast }{\pi }}&M}%
\end{array}%
,\leqno(4.1.6)
\end{equation*}%
where, $\overset{\ast }{\Pi }$\textit{\ }is defined by\textit{\ }%
\begin{equation*}
\begin{array}[b]{c}
\overset{\ast }{\Pi }\left( Y_{a}\overset{\cdot }{\tilde{\partial}}%
^{a}\left( \overset{\ast }{u}_{x}\right) \right) =Y_{a}\left( \overset{\ast }%
{u}_{x}\right) s^{a}\left( \overset{\ast }{\pi }\left( \overset{\ast }{u}%
_{x}\right) \right) ,%
\end{array}%
\leqno(4.1.7)
\end{equation*}%
is canonical projection $\mathbf{B}^{\mathbf{v}}$-morphism$.$

\textbf{Theorem 4.1.1 }\emph{If }$\left( \rho ,\eta \right) \Gamma $\emph{\
is a }$\left( \rho ,\eta \right) $\emph{-connection for the vector bundle }$%
\left( \overset{\ast }{E},\overset{\ast }{\pi },M\right) ,$\emph{\ then its
components satisfy the law of transformation }%
\begin{equation*}
\begin{array}[b]{c}
\left( \rho ,\eta \right) \Gamma _{b%
{\acute{}}%
\gamma
{\acute{}}%
}=M_{b%
{\acute{}}%
}^{b}{\circ }\overset{\ast }{\pi }\left[ -\rho _{\gamma }^{k}\circ h\circ
\overset{\ast }{\pi }\cdot \frac{\partial M_{b}^{a%
{\acute{}}%
}{\circ }\overset{\ast }{\pi }}{\partial x^{k}}p_{a%
{\acute{}}%
}+\left( \rho ,\eta \right) \Gamma _{b\gamma }\right] \Lambda _{\gamma
{\acute{}}%
}^{\gamma }\circ \left( h\circ \overset{\ast }{\pi }\right) .%
\end{array}%
\leqno(4.1.8)
\end{equation*}

\emph{In the particular case of Lie algebroids, }$\left( \eta ,h\right)
=\left( Id_{M},Id_{M}\right) ,$\emph{\ the relations }$\left( 4.1.8\right) $%
\emph{\ become}%
\begin{equation*}
\begin{array}[b]{c}
\rho \Gamma _{b%
{\acute{}}%
\gamma
{\acute{}}%
}=M_{b%
{\acute{}}%
}^{b}{\circ }\overset{\ast }{\pi }\!\!\left[ -\rho _{\gamma }^{k}\circ
\overset{\ast }{\pi }\cdot \frac{\partial M_{b}^{a%
{\acute{}}%
}{\circ }\overset{\ast }{\pi }}{\partial x^{k}}p_{a%
{\acute{}}%
}+\rho \Gamma _{b\gamma }\right] \Lambda _{\gamma
{\acute{}}%
}^{\gamma }\circ \overset{\ast }{\pi }.%
\end{array}%
\leqno(4.1.8^{\prime })
\end{equation*}

\emph{In the classical case, }$\left( \rho ,\eta ,h\right) =\left(
Id_{TM},Id_{M},Id_{M}\right) ,$\emph{\ the relations }$\left( 4.1.8^{\prime
}\right) $\emph{\ become}%
\begin{equation*}
\begin{array}[b]{c}
\Gamma _{j%
{\acute{}}%
k%
{\acute{}}%
}=\frac{\partial x^{j}}{\partial x^{j%
{\acute{}}%
}}{\circ }\overset{\ast }{\tau }\!\!_{M}\left[ -\frac{\partial }{\partial
x^{k}}\left( \frac{\partial x^{i%
{\acute{}}%
}}{\partial x^{j}}{\circ }\overset{\ast }{\tau }\!\!_{M}\right) p_{i%
{\acute{}}%
}+\Gamma _{jk}\right] \frac{\partial x^{k}}{\partial x^{k%
{\acute{}}%
}}\circ \overset{\ast }{\tau }\!\!_{M}.%
\end{array}%
\leqno(4.1.8^{\prime \prime })
\end{equation*}

\emph{Proof.} Let $\left( \overset{\ast }{\Pi },\overset{\ast }{\pi }\right)
$ be the canonical projection $\mathbf{B}$-morphism.

Obviously, the components of
\begin{equation*}
\begin{array}[b]{c}
\overset{\ast }{\Pi }\circ \left( \rho ,\eta \right) \Gamma \left( Z^{\gamma
}\overset{\ast }{\tilde{\partial}}_{\gamma }+Y_{a}\overset{\cdot }{\tilde{%
\partial}}^{a}\right) \left( \overset{\ast }{u}_{x}\right)%
\end{array}%
\end{equation*}%
are the real numbers%
\begin{equation*}
\begin{array}[b]{c}
\left( Y_{b%
{\acute{}}%
}-\left( \rho ,\eta \right) \Gamma _{b%
{\acute{}}%
\gamma
{\acute{}}%
}Z^{\gamma
{\acute{}}%
}\right) \left( \overset{\ast }{u}_{x}\right) .%
\end{array}%
\end{equation*}

Since
\begin{equation*}
\begin{array}{l}
\displaystyle\left( Z^{\gamma ^{\prime }}\overset{\ast }{\tilde{\partial}}%
_{\gamma ^{\prime }}+Y_{b^{\prime }}\overset{\cdot }{\tilde{\partial}}%
^{b^{\prime }}\right) \left( \overset{\ast }{u}_{x}\right) =Z^{\gamma
{\acute{}}%
}\Lambda _{\gamma
{\acute{}}%
}^{\gamma }\circ h\circ \overset{\ast }{\pi }\cdot \overset{\ast }{\tilde{%
\partial}}_{\alpha }\left( \overset{\ast }{u}_{x}\right) \vspace*{1mm} \\
\qquad \displaystyle+\left( Z^{\gamma
{\acute{}}%
}\rho _{\gamma
{\acute{}}%
}^{i%
{\acute{}}%
}\circ h\circ \overset{\ast }{\pi }\frac{\partial M_{b}^{a%
{\acute{}}%
}{\circ }\pi }{\partial x^{i%
{\acute{}}%
}}p_{a%
{\acute{}}%
}+M_{b}^{b%
{\acute{}}%
}Y_{b%
{\acute{}}%
}\right) \overset{\cdot }{\tilde{\partial}}^{b}\left( \overset{\ast }{u}%
_{x}\right) ,%
\end{array}%
\end{equation*}%
it results that the components of%
\begin{equation*}
\begin{array}[b]{c}
\overset{\ast }{\Pi }\circ \left( \rho ,\eta \right) \Gamma \left( Z^{\gamma
^{\prime }}\overset{\ast }{\tilde{\partial}}_{\gamma ^{\prime
}}+Y_{b^{\prime }}\overset{\cdot }{\tilde{\partial}}^{b^{\prime }}\right)
\left( \overset{\ast }{u}_{x}\right)%
\end{array}%
\end{equation*}%
are the real numbers
\begin{equation*}
\begin{array}[b]{c}
\left( Z^{\gamma
{\acute{}}%
}\rho _{\gamma
{\acute{}}%
}^{k%
{\acute{}}%
}\circ h\circ \overset{\ast }{\pi }\frac{\partial M_{b}^{a%
{\acute{}}%
}{\circ }\overset{\ast }{\pi }}{\partial x^{k%
{\acute{}}%
}}p_{a%
{\acute{}}%
}+M_{b}^{b%
{\acute{}}%
}\circ \overset{\ast }{\pi }Y_{b%
{\acute{}}%
}-\left( \rho ,\eta \right) \Gamma _{b\gamma }Z^{\gamma
{\acute{}}%
}\Lambda _{\gamma
{\acute{}}%
}^{\gamma }\circ h\circ \overset{\ast }{\pi }\right) M_{b%
{\acute{}}%
}^{b}\circ \overset{\ast }{\pi }\left( \overset{\ast }{u}_{x}\right) ,%
\end{array}%
\end{equation*}%
where $\left\Vert M_{b%
{\acute{}}%
}^{b}\right\Vert =\left\Vert M_{b}^{b%
{\acute{}}%
}\right\Vert ^{-1}.$\medskip

Therefore, we have:%
\begin{equation*}
\begin{array}{l}
\displaystyle\left( Z^{\gamma
{\acute{}}%
}\rho _{\gamma
{\acute{}}%
}^{k%
{\acute{}}%
}\circ h\circ \overset{\ast }{\pi }\frac{\partial M_{b}^{a%
{\acute{}}%
}{\circ }\overset{\ast }{\pi }}{\partial x^{k%
{\acute{}}%
}}p_{a%
{\acute{}}%
}+M_{b}^{b%
{\acute{}}%
}\circ \overset{\ast }{\pi }Y_{b%
{\acute{}}%
}-\left( \rho ,\eta \right) \Gamma _{b\gamma }Z^{\gamma
{\acute{}}%
}\Lambda _{\gamma
{\acute{}}%
}^{\gamma }\circ h\circ \overset{\ast }{\pi }\right) M_{b%
{\acute{}}%
}^{b}\circ \overset{\ast }{\pi }\vspace*{1mm} \\
\qquad \qquad \displaystyle=Y_{b%
{\acute{}}%
}-\left( \rho ,\eta \right) \Gamma _{b%
{\acute{}}%
\gamma
{\acute{}}%
}Z^{\gamma
{\acute{}}%
}.%
\end{array}%
\end{equation*}

After some calculations we obtain:
\begin{equation*}
\begin{array}[b]{c}
\left( \rho ,\eta \right) \Gamma _{b%
{\acute{}}%
\gamma
{\acute{}}%
}=M_{b%
{\acute{}}%
}^{b}\circ \overset{\ast }{\pi }\left( -\rho _{\gamma }^{k}\circ h\circ
\overset{\ast }{\pi }\cdot \frac{\partial M_{b}^{a%
{\acute{}}%
}{\circ }\overset{\ast }{\pi }}{\partial x^{k}}p_{a%
{\acute{}}%
}+\left( \rho ,\eta \right) \Gamma _{b\gamma }\right) \Lambda _{\gamma
{\acute{}}%
}^{\gamma }\circ h\circ \overset{\ast }{\pi }.%
\end{array}%
\eqno{q.e.d.}
\end{equation*}

\emph{Remark 4.1.1 }If we have a set of real local functions $\left( \rho
,\eta \right) \Gamma _{b\gamma }$ which satisfies the relations of passing $%
\left( 4.1.8\right) ,$ then we have a $\left( \rho ,\eta \right) $%
-connection $\left( \rho ,\eta \right) \Gamma $ for the dual vector bundle $%
\left( \overset{\ast }{E},\overset{\ast }{\pi },M\right) .$

\textbf{Example 4.1.1} If $\Gamma $ is a Ehresmann connection for the vector
bundle $\left( \overset{\ast }{E},\overset{\ast }{\pi },M\right) $ on
components $\Gamma _{bk},$ then the differentiable real local functions
\begin{equation*}
\left( \rho ,\eta \right) \Gamma _{b\gamma }=\left( \rho _{\gamma }^{k}\circ
h\circ \overset{\ast }{\pi }\right) \Gamma _{bk}
\end{equation*}%
are the components of a $\left( \rho ,\eta \right) $-connection $\left( \rho
,\eta \right) \Gamma $ for the vector bundle $\left( \overset{\ast }{E},%
\overset{\ast }{\pi },M\right) $ which will be called the $\left( \rho ,\eta
\right) $\emph{-connection associated to the connection }$\Gamma .$

\textbf{Definition 4.1.3 }If $\left( \rho ,\eta \right) \Gamma $ is a $%
\left( \rho ,\eta \right) $-connection for the vector bundle $\left( \overset%
{\ast }{E},\overset{\ast }{\pi },M\right) $, then for any
\begin{equation*}
z=z^{\gamma }t_{\gamma }\in \Gamma \left( F,\nu ,N\right)
\end{equation*}%
the application%
\begin{equation*}
\begin{array}{rcl}
\Gamma \left( \overset{\ast }{E},\overset{\ast }{\pi },M\right) & ^{%
\underrightarrow{\left( \rho ,\eta \right) D_{z}}} & \Gamma \left( \overset{%
\ast }{E},\overset{\ast }{\pi },M\right) \\
\overset{\ast }{u}=u_{a}s^{a} & \longmapsto & \left( \rho ,\eta \right) D_{z}%
\overset{\ast }{u}%
\end{array}%
\end{equation*}%
defined by
\begin{equation*}
\begin{array}[b]{c}
\left( \rho ,\eta \right) D_{z}\overset{\ast }{u}=z^{\gamma }\circ h\left(
\rho _{\gamma }^{k}\circ h\frac{\partial u_{b}}{\partial x^{k}}-\left( \rho
,\eta \right) \Gamma _{b\gamma }\circ \overset{\ast }{u}\right) s^{b},%
\end{array}%
\leqno(4.1.9)
\end{equation*}%
will be called the \emph{covariant }$\left( \rho ,\eta \right) $\emph{%
-derivative associated to }$\left( \rho ,\eta \right) $\emph{-con\-nec\-tion
}$\left( \rho ,\eta \right) \Gamma $\emph{\ with respect to section }$z$%
\emph{.}

In the particular case of Lie algebroids, $\left( \eta ,h\right) =\left(
Id_{M},Id_{M}\right) ,$\ the relations $\left( 4.1.9\right) $\ become%
\begin{equation*}
\begin{array}[b]{c}
\rho D_{z}\overset{\ast }{u}=z^{\gamma }\left( \rho _{\gamma }^{k}\frac{%
\partial u_{b}}{\partial x^{k}}-\rho \Gamma _{b\gamma }\circ \overset{\ast }{%
u}\right) s^{b}.%
\end{array}%
\leqno(4.1.9^{\prime })
\end{equation*}

In the classical case, $\left( \rho ,\eta ,h\right) =\left(
Id_{TM},Id_{M},Id_{M}\right) ,$\ the relations $\left( 4.1.9^{\prime
}\right) $\ become%
\begin{equation*}
\begin{array}[b]{c}
D_{X}\omega =X^{k}\left( \frac{\partial \omega _{j}}{\partial x^{k}}-\Gamma
_{jk}\circ \omega \right) dx^{j}.%
\end{array}%
\leqno(4.1.9^{\prime \prime })
\end{equation*}%
\textbf{Definition 4.1.4 }We will say that the $\left( \rho ,\eta \right) $%
\emph{-connection }$\left( \rho ,\eta \right) \Gamma $\emph{\ is homogeneous}
or \emph{linear }if the local real functions $\left( \rho ,\eta \right)
\Gamma _{b\gamma }$ are homogeneous or linear on the fibre of vector bundle $%
\left( \overset{\ast }{E},\overset{\ast }{\pi },M\right) $ respectively.

\emph{Remark 4.1.2} If $\left( \rho ,\eta \right) \Gamma $ is a linear $%
\left( \rho ,\eta \right) $-connection for the vector bundle $\left( \overset%
{\ast }{E},\overset{\ast }{\pi },M\right) $, then for each local vector $%
\left( m+r\right) $-chart $\left( U,\overset{\ast }{s}_{U}\right) $ and for
each local vector $\left( n+p\right) $-chart $\left( V,t_{V}\right) $ such
that $U\cap h^{-1}\left( V\right) \neq \phi $, there exists the
differentiable real functions $\rho \Gamma _{b\gamma }^{a}$ defined on $%
U\cap h^{-1}\left( V\right) $ such that
\begin{equation*}
\begin{array}{c}
\left( \rho ,\eta \right) \Gamma _{b\gamma }\circ \overset{\ast }{u}=\left(
\rho ,\eta \right) \Gamma _{b\gamma }^{a}\cdot u_{a},%
\end{array}%
\leqno(4.1.10)
\end{equation*}%
for any $\overset{\ast }{u}=u_{a}s^{a}\in \Gamma \left( \overset{\ast }{E},%
\overset{\ast }{\pi },M\right) .$

The differentiable real local functions $\left( \rho ,\eta \right) \Gamma
_{b\gamma }^{a}$ will be called the \emph{Christoffel coefficients of linear
}$\left( \rho ,\eta \right) $\emph{-connection }$\left( \rho ,\eta \right)
\Gamma .$

\textbf{Theorem 4.1.2} \emph{If }$\left( \rho ,\eta \right) \Gamma $\emph{\
is a linear }$\left( \rho ,\eta \right) $\emph{-connection for the vector
bundle }$\left( \overset{\ast }{E},\overset{\ast }{\pi },M\right) ,$\emph{\
then its components satisfy the law of transformation }%
\begin{equation*}
\begin{array}[b]{c}
\left( \rho ,\eta \right) \Gamma _{b%
{\acute{}}%
\gamma
{\acute{}}%
}^{a%
{\acute{}}%
}=M_{b%
{\acute{}}%
}^{b}\left[ -\rho _{\gamma }^{k}\circ h\frac{\partial M_{b}^{a%
{\acute{}}%
}}{\partial x^{k}}+\left( \rho ,\eta \right) \Gamma _{b\gamma }^{a}M_{a}^{a%
{\acute{}}%
}\right] \Lambda _{\gamma
{\acute{}}%
}^{\gamma }\circ h.%
\end{array}%
\leqno(4.1.11)
\end{equation*}

\emph{In the particular case of Lie algebroids, }$\left( \eta ,h\right)
=\left( Id_{M},Id_{M}\right) ,$\emph{\ the relations }$\left( 4.1.11\right) $%
\emph{\ become}%
\begin{equation*}
\begin{array}[b]{c}
\rho \Gamma _{b%
{\acute{}}%
\gamma
{\acute{}}%
}^{a%
{\acute{}}%
}=M_{b%
{\acute{}}%
}^{b}\left[ -\rho _{\gamma }^{k}\frac{\partial M_{b}^{a%
{\acute{}}%
}}{\partial x^{k}}+\rho \Gamma _{b\gamma }^{a}M_{a}^{a%
{\acute{}}%
}\right] \Lambda _{\gamma
{\acute{}}%
}^{\gamma }.%
\end{array}%
\leqno(4.1.11^{\prime })
\end{equation*}

\emph{In the classical case, }$\left( \rho ,\eta ,h\right) =\left(
Id_{TM},Id_{M},Id_{M}\right) ,$\emph{\ the relations }$\left( 4.1.11^{\prime
}\right) $\emph{\ become}%
\begin{equation*}
\begin{array}[b]{c}
\Gamma _{j%
{\acute{}}%
k%
{\acute{}}%
}^{i%
{\acute{}}%
}=\frac{\partial x^{j}}{\partial x^{j%
{\acute{}}%
}}\left[ -\frac{\partial }{\partial x^{k}}\left( \frac{\partial x^{i%
{\acute{}}%
}}{\partial x^{j}}\right) +\Gamma _{jk}^{i}\frac{\partial x^{i%
{\acute{}}%
}}{\partial x^{i}}\right] \frac{\partial x^{k}}{\partial x^{k%
{\acute{}}%
}}.%
\end{array}%
\leqno(4.1.11^{\prime \prime })
\end{equation*}

\emph{Remark 4.1.3 }Since
\begin{equation*}
\begin{array}[b]{c}
\frac{\partial M_{b}^{a%
{\acute{}}%
}}{\partial x^{i}}M_{b%
{\acute{}}%
}^{b}+\frac{\partial M_{b%
{\acute{}}%
}^{b}}{\partial x^{i}}M_{b}^{a%
{\acute{}}%
}=0,%
\end{array}%
\end{equation*}%
it results that the relations $(4.1.11)$ are equivalent with the relations $%
(4.10^{\prime }).$

\textbf{Definition 4.1.5 }If $\left( \rho ,\eta \right) \Gamma $ is a linear
$\left( \rho ,\eta \right) $-connection for the vector bundle $\left(
\overset{\ast }{E},\overset{\ast }{\pi },M\right) $, then for any
\begin{equation*}
z=z^{\gamma }t_{\gamma }\in \Gamma \left( F,\nu ,N\right)
\end{equation*}%
the application%
\begin{equation*}
\begin{array}{rcl}
\Gamma \left( \overset{\ast }{E},\overset{\ast }{\pi },M\right) \!\! & \!\!^{%
\underrightarrow{\left( \rho ,\eta \right) D_{z}}}\!\! & \!\!\Gamma \left(
\overset{\ast }{E},\overset{\ast }{\pi },M\right) \\
\overset{\ast }{u}{=}u_{a}s^{a}\!\! & \!\!\longmapsto \!\! & \!\!(\rho ,\eta
)D_{z}\overset{\ast }{u}%
\end{array}%
\end{equation*}%
defined by
\begin{equation*}
\begin{array}[b]{c}
(\rho ,\eta )D_{z}\overset{\ast }{u}{=}z^{\gamma }{\circ }h\left( \rho
_{\gamma }^{k}{\circ }h\frac{\partial u_{b}}{\partial x^{k}}{-}(\rho ,\eta
)\Gamma _{b\gamma }^{a}\cdot u_{a}\right) s^{b}%
\end{array}%
\leqno(4.1.12)
\end{equation*}%
will be called the \emph{covariant }$\left( \rho ,\eta \right) $\emph{%
-derivative associated to linear }$\left( \rho ,\eta \right) $\emph{%
-connection }$\left( \rho ,\eta \right) \Gamma $\emph{\ with respect to
section }$z$.

In the particular case of Lie algebroids, $\left( \eta ,h\right) =\left(
Id_{M},Id_{M}\right) ,$\ the relations $\left( 4.1.12\right) $\ become%
\begin{equation*}
\begin{array}[b]{c}
\rho D_{z}\overset{\ast }{u}{=}z^{\gamma }\left( \rho _{\gamma }^{k}\frac{%
\partial u_{b}}{\partial x^{k}}{-}\rho \Gamma _{b\gamma }^{a}\cdot
u_{a}\right) s^{b}%
\end{array}%
\leqno(4.1.12^{\prime })
\end{equation*}

In the classical case, $\left( \rho ,\eta ,h\right) =\left(
Id_{TM},Id_{M},Id_{M}\right) ,$\ the relations $\left( 4.1.12^{\prime
}\right) $\ become%
\begin{equation*}
\begin{array}[b]{c}
D_{X}\omega {=}X^{k}\left( \frac{\partial \omega _{j}}{\partial x^{k}}{-}%
\Gamma _{jk}^{i}\cdot \omega _{i}\right) dx^{j}%
\end{array}%
\leqno(4.1.12^{\prime \prime })
\end{equation*}

In the next we use the same notation $\left( \rho ,\eta \right) \Gamma $ for
the linear $\left( \rho ,\eta \right) $-connection for the vector bundle $%
\left( E,\pi ,M\right) $ or for its dual $\left( \overset{\ast }{E},\overset{%
\ast }{\pi },M\right) $

\emph{Remark 4.1.4 }If $\left( \rho ,\eta \right) \Gamma $ is a linear $%
\left( \rho ,\eta \right) $-connection for the vector bundle $\left( E,\pi
,M\right) $ or for its dual $\left( \overset{\ast }{E},\overset{\ast }{\pi }%
,M\right) $ then, the tensor fields algebra
\begin{equation*}
\left( \mathcal{T}\left( E,\pi ,M\right) ,+,\cdot ,\otimes \right)
\end{equation*}%
is endowed with the $\left( \rho ,\eta \right) $-derivative
\begin{equation*}
\begin{array}{rcl}
\Gamma \left( F,\nu ,N\right) \times \mathcal{T}\left( E,\pi ,M\right) & ^{%
\underrightarrow{\left( \rho ,\eta \right) D}} & \mathcal{T}\left( E,\pi
,M\right) \\
\left( z,T\right) & \longmapsto & \left( \rho ,\eta \right) D_{z}T%
\end{array}%
\end{equation*}%
defined for a tensor field $T\in \mathcal{T}_{q}^{p}\left( E,\pi ,M\right) $
by the relation:
\begin{equation*}
\begin{array}{l}
\left( \rho ,\eta \right) D_{z}T\left( \overset{\ast }{u}_{1},...,\overset{%
\ast }{u}_{p},u_{1},...,u_{q}\right) =\Gamma \left( \rho ,\eta \right)
\left( z\right) \left( T\left( \overset{\ast }{u}_{1},...,\overset{\ast }{u}%
_{p},u_{1},...,u_{q}\right) \right) \vspace*{1mm} \\
-T\left( \left( \rho ,\eta \right) D_{z}\overset{\ast }{u}_{1},...,\overset{%
\ast }{u}_{p},u_{1},...,u_{q}\right) -...-T\left( \overset{\ast }{u}%
_{1},...,\left( \rho ,\eta \right) D_{z}\overset{\ast }{u}%
_{p},u_{1},...,u_{q}\right) \vspace*{1mm} \\
-T\left( \overset{\ast }{u}_{1},...,\overset{\ast }{u}_{p},\left( \rho ,\eta
\right) D_{z}u_{1},...,u_{q}\right) -...-T\left( \overset{\ast }{u}_{1},...,%
\overset{\ast }{u}_{p},u_{1},...,\left( \rho ,\eta \right) D_{z}u_{q}\right)
.%
\end{array}%
\leqno(4.1.13)
\end{equation*}

Moreover, it satisfies the condition
\begin{equation*}
\begin{array}{c}
\left( \rho ,\eta \right) D_{f_{1}z_{1}+f_{2}z_{2}}T=f_{1}\left( \rho ,\eta
\right) D_{z_{1}}T+f_{2}\left( \rho ,\eta \right) D_{z_{2}}T.%
\end{array}%
\leqno(4.1.14)
\end{equation*}

Consequently, if the tensor algebra $\left( \mathcal{T}\left( E,\pi
,M\right) ,+,\cdot ,\otimes \right) $ is endowed with a $\left( \rho ,\eta
\right) $-derivative defined for a tensor field $T\in \mathcal{T}%
_{q}^{p}\left( E,\pi ,M\right) $ by $\left( 4.1.13\right) $ which satisfies
the condition $\left( 4.1.14\right) $, then we can endowed $\left( E,\pi
,M\right) $ with a linear $\left( \rho ,\eta \right) $-connection $\left(
\rho ,\eta \right) \Gamma $ such that its components are defined by the
equality:
\begin{equation*}
\begin{array}{c}
\left( \rho ,\eta \right) D_{t_{\gamma }}s_{b}=\left( \rho ,\eta \right)
\Gamma _{b\gamma }^{a}s_{a}%
\end{array}%
\end{equation*}%
or
\begin{equation*}
\begin{array}{c}
\left( \rho ,\eta \right) D_{t_{\gamma }}s^{a}=-\left( \rho ,\eta \right)
\Gamma _{b\gamma }^{a}s^{b}.%
\end{array}%
\end{equation*}

The $\left( \rho ,\eta \right) $-derivative defined by $\left( 4.1.13\right)
$ will be called the \emph{covariant }$\left( \rho ,\eta \right) $\emph{%
-derivative.}

After some calculations, we obtain:
\begin{equation*}
\begin{array}{l}
\left( \rho ,\eta \right) D_{z}\left(
T_{b_{1},...,b_{q}}^{a_{1},...,a_{p}}s_{a_{1}}\otimes ...\otimes
s_{a_{p}}\otimes s^{b_{1}}\otimes ...\otimes s^{b_{q}}\right) \vspace*{1mm}
\\
=z^{\gamma }\circ h\left( \rho _{\gamma }^{k}\circ h\displaystyle\frac{%
\partial T_{b_{1},...,b_{q}}^{a_{1},...,a_{p}}}{\partial x^{k}}+\left( \rho
,\eta \right) \Gamma _{a\gamma
}^{a_{1}}T_{b_{1},...,b_{q}}^{a,a_{2}...,a_{p}}\right. \vspace*{1mm} \\
+\left( \rho ,\eta \right) \Gamma _{a\gamma
}^{a_{2}}T_{b_{1},...,b_{q}}^{a_{1},a...,a_{p}}+...+\left( \rho ,\eta
\right) \Gamma _{a\gamma }^{a_{p}}T_{b_{1},...,b_{q}}^{a_{1},a_{2}...,a}-...%
\vspace*{2mm} \\
-\left( \rho ,\eta \right) \Gamma _{b_{1}\gamma
}^{b}T_{b,b_{2},...,b_{q}}^{a_{1},a_{2}...,a_{p}}-\left( \rho ,\eta \right)
\Gamma _{b_{2}\gamma }^{b}T_{b_{1},b,...,b_{q}}^{a_{1},a_{2}...,a_{p}}-...%
\vspace*{2mm} \\
\left. -\left( \rho ,\eta \right) \Gamma _{b_{q}\gamma
}^{b}T_{b_{1},b_{2},...,b}^{a_{1},a_{2}...,a_{p}}\right) s_{a_{1}}\otimes
...\otimes s_{a_{p}}\otimes s^{b_{1}}\otimes ...\otimes s^{b_{q}}\vspace*{1mm%
}\vspace*{2mm} \\
\overset{put}{=}z^{\gamma }\circ hT_{b_{1},...,b_{q}\mid \gamma
}^{a_{1},...,a_{p}}s_{a_{1}}\otimes ...\otimes s_{a_{p}}\otimes
s^{b_{1}}\otimes ...\otimes s^{b_{q}}.%
\end{array}%
\leqno(4.1.15)
\end{equation*}

We remark that if $\left( \rho ,\eta \right) \Gamma $ is the linear $\left(
\rho ,\eta \right) $\textit{-}connection associated to the Ehresmann linear
connection $\Gamma ,$ namely $\left( \rho ,\eta \right) \Gamma _{b\alpha
}^{a}=\left( \rho _{\alpha }^{k}\circ h\right) \Gamma _{bk}^{a},$ then
\begin{equation*}
\begin{array}{c}
T_{b_{1},...,b_{q}\mid \gamma }^{a_{1},...,a_{p}}=\left( \rho _{\gamma
}^{k}\circ h\right) T_{b_{1},...,b_{q}\mid k}^{a_{1},...,a_{p}}.%
\end{array}%
\end{equation*}

In the particular case of Lie algebroids, $\left( \eta ,h\right) =\left(
Id_{M},Id_{M}\right) ,$\ the relations $(4.1.15)$ become
\begin{equation*}
\begin{array}{l}
\rho D_{z}\left( T_{b_{1},...,b_{q}}^{a_{1},...,a_{p}}s_{a_{1}}\otimes
...\otimes s_{a_{p}}\otimes s^{b_{1}}\otimes ...\otimes s^{b_{q}}\right)
\vspace*{1mm} \\
\qquad =z^{\gamma }\left( \rho _{\gamma }^{k}\displaystyle\frac{\partial
T_{b_{1},...,b_{q}}^{a_{1},...,a_{p}}}{\partial x^{k}}+\rho \Gamma _{a\gamma
}^{a_{1}}T_{b_{1},...,b_{q}}^{a,a_{2}...,a_{p}}\right. \vspace*{1mm} \\
\qquad +\rho \Gamma _{a\gamma
}^{a_{2}}T_{b_{1},...,b_{q}}^{a_{1},a...,a_{p}}+...+\rho \Gamma _{a\gamma
}^{a_{p}}T_{b_{1},...,b_{q}}^{a_{1},a_{2}...,a}-...\vspace*{1mm} \\
\qquad -\rho \Gamma _{b_{1}\gamma
}^{b}T_{b,b_{2},...,b_{q}}^{a_{1},a_{2}...,a_{p}}-\rho \Gamma _{b_{2}\gamma
}^{b}T_{b_{1},b,...,b_{q}}^{a_{1},a_{2}...,a_{p}}-...\vspace*{1mm} \\
\qquad \left. -\rho \Gamma _{b_{q}\gamma
}^{b}T_{b_{1},b_{2},...,b}^{a_{1},a_{2}...,a_{p}}\right) s_{a_{1}}\otimes
...\otimes s_{a_{p}}\otimes s^{b_{1}}\otimes ...\otimes s^{b_{q}}\vspace*{1mm%
} \\
\qquad \overset{put}{=}z^{\gamma }T_{b_{1},...,b_{q}\mid \gamma
}^{a_{1},...,a_{p}}s_{a_{1}}\otimes ...\otimes s_{a_{p}}\otimes
s^{b_{1}}\otimes ...\otimes s^{b_{q}}.%
\end{array}%
\leqno(4.1.15^{\prime })
\end{equation*}

In the classical case, $\left( \rho ,\eta ,h\right) =\left(
Id_{TM},Id_{M},Id_{M}\right) ,$\ the relations $\left( 4.1.15^{\prime
}\right) $\ become
\begin{equation*}
\begin{array}{l}
D_{X}\left( T_{j_{1}...j_{q}}^{i_{1}...i_{p}}\partial _{i_{1}}\otimes
...\otimes \partial _{i_{p}}\otimes dx^{j_{1}}\otimes ...\otimes
dx^{j_{q}}\right) \vspace*{1mm} \\
\qquad =X^{k}\left( \frac{\partial T_{j_{1}...j_{q}}^{i_{1}...i_{p}}}{%
\partial x^{k}}+\Gamma
_{ik}^{i_{1}}T_{j_{1}...j_{q}}^{ii_{2}...i_{p}}+...+\Gamma
_{ik}^{i_{p}}T_{j_{1}...j_{q}}^{i_{1}...i_{p-1}i}\right. \vspace*{1mm}%
\vspace*{1mm} \\
\qquad -\Gamma _{j_{1}k}^{j}T_{jj_{2}...j_{q}}^{i_{1}...i_{p}}-...\vspace*{%
1mm}\left. -\Gamma _{j_{q}k}^{j}T_{j_{1}...j_{q-1}j}^{i_{1}...i_{p}}\right)
\partial _{i_{1}}\otimes ...\otimes \partial _{i_{p}}\otimes
dx^{j_{1}}\otimes ...\otimes dx^{j_{q}}\vspace*{1mm} \\
\qquad \overset{put}{=}X^{k}T_{j_{1}...j_{q}\mid k}^{i_{1}...i_{p}}\partial
_{i_{1}}\otimes ...\otimes \partial _{i_{p}}\otimes dx^{j_{1}}\otimes
...\otimes dx^{j_{q}}.%
\end{array}%
\leqno(4.1.15^{\prime \prime })
\end{equation*}

\section{Torsion and curvature. Formulas of Ricci and Bianchi type}

We apply our theory for the diagram:
\begin{equation*}
\begin{array}{c}
\xymatrix{E\ar[d]_\pi&\left( F,\left[ , \right] _{F,h},\left( \rho
,Id_{N}\right) \right)\ar[d]^\nu \\ M\ar[r]^h&N}%
\end{array}%
,\leqno(5.1)
\end{equation*}%
where $\left( E,\pi ,M\right) \in \left\vert \mathbf{B}^{\mathbf{v}%
}\right\vert $ and $\left( \left( F,\nu ,N\right) ,\left[ ,\right]
_{F,h},\left( \rho ,Id_{N}\right) \right) \in \left\vert \mathbf{GLA}%
\right\vert .$

Let $\rho \Gamma $ be a linear $\rho $-connection for the vector bundle $%
\left( E,\pi ,M\right) $ by components $\rho \Gamma _{b\alpha }^{a}.$

Using the components of the linear $\rho $-connection $\rho \Gamma $, then
we obtain a linear $\rho $-connection $\rho \dot{\Gamma}$ for the vector
bundle $\left( E,\pi ,M\right) $ given by the diagram:
\begin{equation*}
\begin{array}{ccl}
~\ \ \ E &  & \left( h^{\ast }F,\left[ ,\right] _{h^{\ast }F},\left( \overset%
{h^{\ast }F}{\rho },Id_{M}\right) \right) \\
\pi \downarrow &  & ~\ \ \downarrow h^{\ast }\nu \\
~\ \ \ M & ^{\underrightarrow{~\ \ \ \ Id_{M}~\ \ }} & ~\ \ \ M%
\end{array}%
\leqno(5.2)
\end{equation*}

If $\left( E,\pi ,M\right) =\left( F,\nu ,N\right) ,$ then, using the
components of the linear $\rho $-connection $\rho \Gamma ,$ we can consider
a linear $\rho $-connection $\rho \ddot{\Gamma}$ for the vector bundle $%
\left( h^{\ast }E,h^{\ast }\pi ,M\right) $ given by the diagram:
\begin{equation*}
\begin{array}{ccl}
~\ \ ~\ \ \ \ h^{\ast }E &  & \left( h^{\ast }E,\left[ ,\right] _{h^{\ast
}E},\left( \overset{h^{\ast }E}{\rho },Id_{M}\right) \right) \\
h^{\ast }\pi \downarrow &  & ~\ \ \downarrow h^{\ast }\pi \\
~\ \ \ \ \ M & ^{\underrightarrow{~\ \ \ \ Id_{M}~\ \ }} & ~\ \ M%
\end{array}%
\leqno(5.3)
\end{equation*}%
\emph{Definition 5.1 }If $\left( E,\pi ,M\right) =\left( F,\nu ,N\right) $,
then the application
\begin{equation*}
\begin{array}{ccc}
\Gamma \left( h^{\ast }E,h^{\ast }\pi ,M\right) ^{2} & ^{\underrightarrow{\
\left( \rho ,h\right) \mathbb{T}\ }} & \Gamma \left( h^{\ast }E,h^{\ast }\pi
,M\right) \\
\left( U,V\right) & \longrightarrow & \rho \mathbb{T}\left( U,V\right)%
\end{array}%
\end{equation*}%
defined by:
\begin{equation*}
\begin{array}{c}
\left( \rho ,h\right) \mathbb{T}\left( U,V\right) =\rho \ddot{D}_{U}V-\rho
\ddot{D}_{V}U-\left[ U,V\right] _{h^{\ast }E},\,%
\end{array}%
\leqno(5.4)
\end{equation*}%
for any $U,V\in \Gamma \left( h^{\ast }E,h^{\ast }\pi ,M\right) ,$ will be
called $\left( \rho ,h\right) $\emph{-torsion associated to the linear }$%
\rho $\emph{-connection }$\rho \Gamma .$

In the particular case of Lie algebroids,\ $h=Id_{M}$, we obtain the
application
\begin{equation*}
\begin{array}{ccc}
\Gamma \left( E,\pi ,M\right) ^{2} & ^{\underrightarrow{\ \rho \mathbb{T}\ }}
& \Gamma \left( E,\pi ,M\right) \\
\left( u,v\right) & \longrightarrow & \rho \mathbb{T}\left( u,v\right)%
\end{array}%
\end{equation*}%
defined by:
\begin{equation*}
\begin{array}{c}
\rho \mathbb{T}\left( u,v\right) =\rho D_{u}v-\rho D_{v}u-\left[ u,v\right]
_{E},\,%
\end{array}%
\leqno(5.4^{\prime })
\end{equation*}%
for any $u,v\in \Gamma \left( E,\pi ,M\right) ,$\ which will be called the $%
\rho $\emph{-torsion associated to the linear }$\rho $\emph{-connection }$%
\rho \Gamma .$

In the classical case, $\rho =Id_{TM}$, we obtain the torsion $\mathbb{T}$
associated to the linear connection~$\Gamma .$

\textbf{Proposition 5.1 }\emph{The }$\left( \rho ,h\right) $\emph{-torsion }$%
\left( \rho ,h\right) \mathbb{T}$\emph{\ associated to the linear }$\rho $%
\emph{-connection }$\rho \Gamma $\emph{\ is }$\mathbb{R}$\emph{-bilinear and
antisymmetric.}

\emph{If }%
\begin{equation*}
\left( \rho ,h\right) \mathbb{T}\left( S_{a},S_{b}\right) \overset{put}{=}%
\left( \rho ,h\right) \mathbb{T}_{~ab}^{c}S_{c}
\end{equation*}%
\emph{\ then }%
\begin{equation*}
\begin{array}{c}
\left( \rho ,h\right) \mathbb{T}_{~ab}^{c}=\rho \Gamma _{ab}^{c}-\rho \Gamma
_{ba}^{c}-L_{ab}^{c}\circ h.%
\end{array}%
\leqno(5.5)
\end{equation*}

\emph{In the particular case of Lie algebroids, }$h=Id_{M},$\emph{\ we have }%
$\rho \mathbb{T}\left( s_{a},s_{b}\right) \overset{put}{=}\rho \mathbb{T}%
_{ab}^{c}s_{c}$\emph{\ and }%
\begin{equation*}
\begin{array}{c}
\rho \mathbb{T}_{~ab}^{c}=\rho \Gamma _{ab}^{c}-\rho \Gamma
_{ba}^{c}-L_{ab}^{c}.%
\end{array}%
\leqno(5.5^{\prime })
\end{equation*}

\emph{In the classical case, }$\left( \rho ,h\right) =\left(
Id_{TM},Id_{M}\right) ,$\emph{\ the equality }$\left( 5.5^{\prime }\right) $%
\emph{\ becomes:}
\begin{equation*}
\begin{array}{c}
\mathbb{T}_{~jk}^{i}=\Gamma _{jk}^{i}-\Gamma _{kj}^{i}.%
\end{array}%
\leqno(5.5^{\prime \prime })
\end{equation*}

\textbf{Definition 5.2 }The application
\begin{equation*}
\begin{array}{ccl}
(\Gamma \left( h^{\ast }F,h^{\ast }\nu ,M\right) ^{2}{\times }\Gamma (E,\pi
,M) & ^{\underrightarrow{\ \left( \rho ,h\right) \mathbb{R}\ }} & \Gamma
(E,\pi ,M) \\
((Z,V),u) & \longrightarrow & \left( \rho ,h\right) \mathbb{R}(Z,V)u%
\end{array}%
\end{equation*}%
defined by
\begin{equation*}
\left( \rho ,h\right) \mathbb{R}\left( Z,V\right) u=\rho \dot{D}_{Z}\left(
\rho \dot{D}_{V}u\right) -\rho \dot{D}_{V}\left( \rho \dot{D}_{Z}u\right)
-\rho \dot{D}_{\left[ Z,V\right] _{h^{\ast }F}}u,\,\leqno(5.6)
\end{equation*}%
for any $Z,V\in \Gamma \left( h^{\ast }F,h^{\ast }\nu ,M\right) $ and$~u\in
\Gamma \left( E,\pi ,M\right) ,$ will be called $\left( \rho ,h\right) $%
\emph{-curvature associated to the linear }$\rho $\emph{-connection }$\rho
\Gamma .$

In the particular case of Lie algebroids, $h=Id_{M}$, we obtain the
application
\begin{equation*}
\begin{array}{ccl}
\Gamma \left( F,\nu ,M\right) ^{2}{\times }\Gamma (E,\pi ,M) & ^{%
\underrightarrow{\ \rho \mathbb{R}\ }} & \Gamma (E,\pi ,M) \\
((z,v),u) & \longrightarrow & \rho \mathbb{R}(z,v)u%
\end{array}%
\end{equation*}%
defined by
\begin{equation*}
\rho \mathbb{R}\left( z,v\right) u=\rho D_{z}\left( \rho D_{v}u\right) -\rho
D_{v}\left( \rho D_{z}u\right) -\rho D_{\left[ z,v\right] _{F}}u,\,\leqno%
(5.6^{\prime })
\end{equation*}%
for any $z,v\in \Gamma \left( F,\nu ,M\right) $ and$~u\in \Gamma \left(
E,\pi ,M\right) ,$ which will be called $\rho $\emph{-curvature associated
to the linear }$\rho $\emph{-connection }$\rho \Gamma .$

In the classical case, $\left( \rho ,h\right) =\left( Id_{TM},Id_{M}\right)
, $ we obtain the curvature $\mathbb{R}$ associated to the linear connection
$\Gamma .$

\textbf{Proposition 5.2 }\emph{The }$\left( \rho ,h\right) $\emph{-curvature
}$\left( \rho ,h\right) \mathbb{R}$\emph{\ associated to the linear }$\rho $%
\emph{-connection }$\rho \Gamma $\emph{, is }$\mathbb{R}$\emph{-linear in
each argument and antisymmetric in the first two arguments.}

\emph{If }%
\begin{equation*}
\left( \rho ,h\right) \mathbb{R}\left( T_{\beta },T_{\alpha }\right) s_{b}%
\overset{put}{=}\left( \rho ,h\right) \mathbb{R}_{b~\alpha \beta }^{a}s_{a},
\end{equation*}%
\emph{then}%
\begin{equation*}
\begin{array}[b]{cl}
\left( \rho ,h\right) \mathbb{R}_{b~\alpha \beta }^{a} & =\rho _{\beta
}^{j}\circ h\frac{\partial \rho \Gamma _{b\alpha }^{a}}{\partial x^{j}}+\rho
\Gamma _{e\beta }^{a}\rho \Gamma _{b\alpha }^{e}-\rho _{\alpha }^{i}\circ h%
\frac{\partial \rho \Gamma _{b\beta }^{a}}{\partial x^{i}} \\
& -\rho \Gamma _{e\alpha }^{a}\rho \Gamma _{b\beta }^{e}+\rho \Gamma
_{b\gamma }^{a}L_{\alpha \beta }^{\gamma }\circ h.%
\end{array}%
\leqno(5.7)
\end{equation*}%
\emph{\ }

\emph{In the particular case of Lie algebroids, }$h=Id_{M},$\emph{\ we
obtain }$\rho \mathbb{R}\left( t_{\beta },t_{\alpha }\right) s_{b}\overset{%
put}{=}\rho \mathbb{R}_{b~\alpha \beta }^{a}s_{a}$\emph{, and }%
\begin{equation*}
\begin{array}[b]{c}
\rho \mathbb{R}_{b~\alpha \beta }^{a}=\rho _{\beta }^{j}\frac{\partial \rho
\Gamma _{b\alpha }^{a}}{\partial x^{j}}+\rho \Gamma _{e\beta }^{a}\rho
\Gamma _{b\alpha }^{e}-\rho _{\alpha }^{i}\frac{\partial \rho \Gamma
_{b\beta }^{a}}{\partial x^{i}}-\rho \Gamma _{e\alpha }^{a}\rho \Gamma
_{b\beta }^{e}+\rho \Gamma _{b\gamma }^{a}L_{\alpha \beta }^{\gamma }.%
\end{array}%
\leqno(5.7^{\prime })
\end{equation*}

\emph{In the classical case, }$\left( \rho ,h\right) =\left(
Id_{TM},Id_{M}\right) ,$\emph{\ we obtain }$\mathbb{R}\left( \partial
_{k},\partial _{h}\right) s_{b}\overset{put}{=}\mathbb{R}_{b~hk}^{a}s_{a}$%
\emph{, and the equality }$\left( 5.7^{\prime }\right) $\emph{\ becomes:}
\begin{equation*}
\begin{array}[b]{c}
\mathbb{R}_{b~hk}^{a}=\frac{\partial \Gamma _{bh}^{a}}{\partial x^{k}}%
+\Gamma _{ek}^{a}\Gamma _{bh}^{e}-\frac{\partial \Gamma _{bk}^{a}}{\partial
x^{h}}-\Gamma _{eh}^{a}\Gamma _{bk}^{e}.%
\end{array}%
\leqno(5.7^{\prime \prime })
\end{equation*}

\textbf{Theorem 5.1 }\emph{For any } $u^{a}s_{a}\in \Gamma \left( E,\pi
,M\right) $ \emph{we shall use the notation}
\begin{equation*}
\begin{array}[b]{c}
u_{~|\alpha \beta }^{a}=\rho _{\beta }^{j}\circ h\frac{\partial }{\partial
x^{j}}\left( u_{~|\alpha }^{a_{1}}\right) +\rho \Gamma _{b_{\beta
}}^{a_{1}}u_{~|\alpha }^{b},%
\end{array}%
\leqno(5.8)
\end{equation*}%
\emph{and we verify the formulas: }
\begin{equation*}
u_{~|\alpha \beta }^{a_{1}}-u_{~|\beta \alpha }^{a_{1}}=u^{a}\left( \rho
,h\right) \mathbb{R}_{a~\alpha \beta }^{a_{1}}-u_{~|\gamma
}^{a_{1}}L_{\alpha \beta }^{\gamma }\circ h.\leqno(5.9)
\end{equation*}

\emph{After some calculations, we obtain}
\begin{equation*}
\left( \rho ,h\right) \mathbb{R}_{a~\alpha \beta }^{a_{1}}=u_{a}\left(
u_{~|\alpha \beta }^{a_{1}}-u_{~|\beta \alpha }^{a_{1}}+u_{~|\gamma
}^{a_{1}}L_{\alpha \beta }^{\gamma }\circ h\right) ,\leqno(5.10)
\end{equation*}%
\emph{where }$u_{a}s^{a}\in \Gamma \left( \overset{\ast }{E},\overset{\ast }{%
\pi },M\right) $\emph{\ such that }$u_{a}u^{b}=\delta _{a}^{b}.$

\emph{In the particular case of Lie algebroids, }$h=Id_{M}$\emph{, the
relations }$\left( 5.10\right) $\emph{\ become}%
\begin{equation*}
\rho \mathbb{R}_{a~\alpha \beta }^{a_{1}}=u_{a}\left( u_{~|\alpha \beta
}^{a_{1}}-u_{~|\beta \alpha }^{a_{1}}+u_{~|\gamma }^{a_{1}}L_{\alpha \beta
}^{\gamma }\right) .\leqno(5.10^{\prime })
\end{equation*}

\emph{In the classical case,} $\left( \rho ,h\right) =\left(
Id_{TM},Id_{M}\right) ,$\emph{\ the relations }$\left( 5.10^{\prime }\right)
$\emph{\ become}%
\begin{equation*}
\mathbb{R}_{a~ij}^{a_{1}}=u_{a}\left(
u_{~|ij}^{a_{1}}-u_{~|ji}^{a_{1}}\right) .\leqno\left( 5.10^{\prime \prime
}\right)
\end{equation*}

\emph{Proof.} Since
\begin{equation*}
\begin{split}
u_{~|\alpha \beta }^{a_{1}}& =\rho _{\beta }^{j}\circ h\left( \frac{\partial
}{\partial x^{j}}\left( \rho _{\alpha }^{i}\circ h\frac{\partial u^{a_{1}}}{%
\partial x^{i}}+\rho \Gamma _{a\alpha }^{a_{1}}u^{a}\right) \right) \\
& +\rho \Gamma _{b\beta }^{a_{1}}\left( \rho _{\alpha }^{i}\circ h\frac{%
\partial u^{b}}{\partial x^{i}}+\rho \Gamma _{a\alpha }^{b}u^{a}\right) \\
& =\rho _{\beta }^{j}\circ h\frac{\partial \rho _{\alpha }^{i}\circ h}{%
\partial x^{j}}\frac{\partial u^{a_{1}}}{\partial x^{i}}+\rho _{\beta
}^{j}\circ h\rho _{\alpha }^{i}\circ h\frac{\partial }{\partial x^{j}}\left(
\frac{\partial u^{a_{1}}}{\partial x^{i}}\right) \\
& +\rho _{\beta }^{j}\circ h\frac{\partial \rho \Gamma _{a\alpha }^{a_{1}}}{%
\partial x^{j}}u^{a}+\rho _{\beta }^{j}\circ h\rho \Gamma _{a\alpha }^{a_{1}}%
\frac{\partial u^{a}}{\partial x^{j}} \\
& +\rho _{\alpha }^{i}\circ h\rho \Gamma _{b\beta }^{a_{1}}\frac{\partial
u^{b}}{\partial x^{i}}+\rho \Gamma _{b\beta }^{a_{1}}\rho \Gamma _{a\alpha
}^{b}u^{a}
\end{split}%
\end{equation*}%
and%
\begin{equation*}
\begin{split}
u_{~|\beta \alpha }^{a_{1}}& =\rho _{\alpha }^{i}\circ h\left( \frac{%
\partial }{\partial x^{i}}\left( \rho _{\beta }^{j}\circ h\frac{\partial
u^{a_{1}}}{\partial x^{j}}+\rho \Gamma _{a\beta }^{a_{1}}u^{a}\right) \right)
\\
& +\rho \Gamma _{b\alpha }^{a_{1}}\left( \rho _{\beta }^{j}\circ h\frac{%
\partial u^{b}}{\partial x^{j}}+\rho \Gamma _{a\beta }^{b}u^{a}\right) \\
& =\rho _{\alpha }^{i}\circ h\frac{\partial \rho _{\beta }^{j}\circ h}{%
\partial x^{i}}\frac{\partial u^{a_{1}}}{\partial x^{j}}+\rho _{\beta
}^{j}\circ h\rho _{\alpha }^{i}\circ h\frac{\partial }{\partial x^{i}}\left(
\frac{\partial u^{a_{1}}}{\partial x^{j}}\right) \\
& +\rho _{\alpha }^{i}\circ h\frac{\partial \rho \Gamma _{a\beta }^{a_{1}}}{%
\partial x^{i}}u^{a}+\rho _{\alpha }^{i}\circ h\rho \Gamma _{a\beta }^{a_{1}}%
\frac{\partial u^{a}}{\partial x^{i}} \\
& +\rho _{\beta }^{j}\circ h\rho \Gamma _{b\alpha }^{a_{1}}\frac{\partial
u^{b}}{\partial x^{j}}+\rho \Gamma _{b\alpha }^{a_{1}}\rho \Gamma _{a\beta
}^{b}u^{a},
\end{split}%
\end{equation*}%
it results that
\begin{equation*}
\begin{split}
u_{~|\alpha \beta }^{a_{1}}-u_{~|\beta \alpha }^{a_{1}}& =\rho _{\beta
}^{j}\circ h\frac{\partial \rho _{\alpha }^{i}\circ h}{\partial x^{j}}\frac{%
\partial u^{a_{1}}}{\partial x^{i}}-\rho _{\alpha }^{i}\circ h\frac{\partial
\rho _{\beta }^{j}\circ h}{\partial x^{i}}\frac{\partial u^{a_{1}}}{\partial
x^{j}} \\
& +\left( \rho _{\beta }^{j}\circ h\rho _{\alpha }^{i}\circ h\frac{\partial
^{2}u^{a_{1}}}{\partial x^{i}\partial x^{j}}-\rho _{\beta }^{j}\circ h\rho
_{\alpha }^{i}\circ h\frac{\partial ^{2}u^{a_{1}}}{\partial x^{j}\partial
x^{i}}\right) \\
& +\left( \rho _{\beta }^{j}\circ h\frac{\partial \rho \Gamma _{a\alpha
}^{a_{1}}}{\partial x^{j}}u^{a}-\rho _{\alpha }^{i}\circ h\frac{\partial
\rho \Gamma _{a\beta }^{a_{1}}}{\partial x^{i}}u^{a}\right) \\
& +\left( \rho _{\beta }^{j}\circ h\rho \Gamma _{a\alpha }^{a_{1}}\frac{%
\partial u^{a}}{\partial x^{j}}-\rho _{\beta }^{j}\circ h\rho \Gamma
_{b\alpha }^{a_{1}}\frac{\partial u^{b}}{\partial x^{j}}\right) \\
& +\left( \rho _{\alpha }^{i}\circ h\rho \Gamma _{b\beta }^{a_{1}}\frac{%
\partial u^{b}}{\partial x^{i}}-\rho _{\alpha }^{i}\circ h\rho \Gamma
_{a\beta }^{a_{1}}\frac{\partial u^{a}}{\partial x^{i}}\right) \\
& +\rho \Gamma _{b\beta }^{a_{1}}\rho \Gamma _{a\alpha }^{b}u^{a}-\rho
\Gamma _{b\alpha }^{a_{1}}\rho \Gamma _{a\beta }^{b}u^{a}.
\end{split}%
\end{equation*}

After some calculations, we obtain:
\begin{equation*}
\begin{split}
u_{~|\alpha \beta }^{a_{1}}-u_{~|\beta \alpha }^{a_{1}}& =L_{\beta \alpha
}^{\gamma }\circ h\rho _{\gamma }^{k}\circ h\frac{\partial u^{a_{1}}}{%
\partial x^{k}} \\
& +\left( \rho _{\beta }^{j}\circ h\frac{\partial \rho \Gamma _{a\alpha
}^{a_{1}}}{\partial x^{j}}u^{a}-\rho _{\alpha }^{i}\circ h\frac{\partial
\rho \Gamma _{a\beta }^{a_{1}}}{\partial x^{i}}u^{a}\right) \\
& +\rho \Gamma _{b\beta }^{a_{1}}\rho \Gamma _{a\alpha }^{b}u^{a}-\rho
\Gamma _{b\alpha }^{a_{1}}\rho \Gamma _{a\beta }^{b}u^{a}.
\end{split}%
\end{equation*}

Since
\begin{eqnarray*}
u^{a}\left( \rho ,h\right) \mathbb{R}_{a~\alpha \beta }^{a_{1}}
&=&u^{a}\left( \rho _{\beta }^{j}\circ h\frac{\partial \rho \Gamma _{a\alpha
}^{a_{1}}}{\partial x^{j}}+\rho \Gamma _{e\beta }^{a_{1}}\rho \Gamma
_{a\alpha }^{e}-\rho _{\alpha }^{i}\circ h\frac{\partial \rho \Gamma
_{a\beta }^{a_{1}}}{\partial x^{i}}\right. \\
&&\left. -\rho \Gamma _{e\alpha }^{a_{1}}\rho \Gamma _{a\beta }^{e}-\rho
\Gamma _{a\gamma }^{a_{1}}L_{\beta \alpha }^{\gamma }\circ h\right) .
\end{eqnarray*}%
and
\begin{equation*}
u_{~|\gamma }^{a_{1}}L_{\alpha \beta }^{\gamma }\circ h=\left( \rho _{\gamma
}^{k}\circ h\frac{\partial u^{a_{1}}}{\partial x^{k}}+\rho \Gamma _{a\gamma
}^{a_{1}}u^{a}\right) L_{\alpha \beta }^{\gamma }\circ h
\end{equation*}%
it results that
\begin{equation*}
\begin{split}
u^{a}\left( \rho ,h\right) \mathbb{R}_{a~\alpha \beta }^{a_{1}}-u_{~|\gamma
}^{a_{1}}L_{\alpha \beta }^{\gamma }\circ h& =-L_{\alpha \beta }^{\gamma
}\circ h\rho _{\gamma }^{k}\circ h\frac{\partial u^{a_{1}}}{\partial x^{k}}
\\
& +\left( \rho _{\beta }^{j}\circ h\frac{\partial \rho \Gamma _{a\alpha
}^{a_{1}}}{\partial x^{j}}u^{a}-\rho _{\alpha }^{i}\circ h\frac{\partial
\rho \Gamma _{a\beta }^{a_{1}}}{\partial x^{i}}u^{a}\right) \\
& +\rho \Gamma _{b\beta }^{a_{1}}\rho \Gamma _{a\alpha }^{b}u^{a}-\rho
\Gamma _{b\alpha }^{a_{1}}\rho \Gamma _{a\beta }^{b}u^{a}.
\end{split}%
\end{equation*}

\begin{flushright}
\emph{q.e.d.}
\end{flushright}

\textbf{Lemma 5.1 }\emph{If }$\left( E,\pi ,M\right) =\left( F,\nu ,N\right)
$\emph{, then, for any} $u^{a}s_{a}\in \Gamma \left( E,\pi ,M\right) ,$\emph{%
\ we have that } $u_{~|c}^{a},~a,c\in \overline{1,n}$ are the \emph{%
components of a tensor field of }$\left( 1,1\right) $\emph{\ type.}

\emph{Proof.} Let $U$ and $U^{\prime }$ be two vector local $\left(
m+n\right) $-charts such that $U\cap U^{\prime }\neq \phi .$

Since $u^{a^{\prime }}\left( x\right) =M_{a}^{a^{\prime }}\left( x\right)
u^{a}\left( x\right) ,~$for any $x\in U\cap U^{\prime },$ it results that
\begin{equation}
\begin{array}[b]{c}
\rho _{c^{\prime }}^{k^{\prime }}\circ h\left( x\right) \frac{\partial
u^{a^{\prime }}\left( x\right) }{\partial x^{k^{\prime }}}=\rho _{c^{\prime
}}^{k^{\prime }}\circ h\left( x\right) \frac{\partial }{\partial
x^{k^{\prime }}}\left( M_{a}^{a^{\prime }}\left( x\right) \right)
u^{a}\left( x\right) +M_{a}^{a^{\prime }}\left( x\right) \rho _{c^{\prime
}}^{k^{\prime }}\circ h\left( x\right) \frac{\partial u^{a}\left( x\right) }{%
\partial x^{k^{\prime }}}.%
\end{array}%
\end{equation}

Since , for any $x\in U\cap U^{\prime }$, we have
\begin{equation}
\begin{array}[b]{c}
\rho \Gamma _{b^{\prime }c^{\prime }}^{a^{\prime }}(x)=M_{a}^{a^{\prime
}}(x)\left( \rho _{c}^{k}\circ h(x)\frac{\partial }{\partial x^{k}}%
(M_{b^{\prime }}^{a}(x))+\rho \Gamma _{bc}^{a}(x)M_{b^{\prime
}}^{b}(x)\right) M_{c^{\prime }}^{c}(x),%
\end{array}%
\end{equation}%
and
\begin{equation}
\begin{array}[b]{c}
0=\frac{\partial }{\partial x^{k^{\prime }}}\left( M_{a}^{a^{\prime }}\left(
x\right) M_{b^{\prime }}^{a}\left( x\right) \right) =\frac{\partial }{%
\partial x^{k^{\prime }}}\left( M_{a}^{a^{\prime }}\left( x\right) \right)
M_{b^{\prime }}^{a}\left( x\right) +M_{a}^{a^{\prime }}\left( x\right) \frac{%
\partial }{\partial x^{k^{\prime }}}\left( M_{b^{\prime }}^{a}\left(
x\right) \right)%
\end{array}%
\end{equation}%
it results that
\begin{equation}
\begin{split}
\rho \Gamma _{b^{\prime }c^{\prime }}^{a^{\prime }}\left( x\right)
u^{b^{\prime }}\left( x\right) & =-\rho _{c^{\prime }}^{k^{\prime }}\circ
h\left( x\right) \frac{\partial }{\partial x^{k^{\prime }}}\left(
M_{a}^{a^{\prime }}\left( x\right) \right) u^{a}\left( x\right) \\
& +M_{a}^{a^{\prime }}\left( x\right) \rho \Gamma _{bc}^{a}\left( x\right)
u^{b}\left( x\right) M_{c^{\prime }}^{c}\left( x\right) .
\end{split}%
\end{equation}

Summing the equalities $\left( 1\right) $ and $\left( 4\right) $, it results
the conclusion of lemma.\hfill \emph{q.e.d.}

\textbf{Theorem 5.2 }\emph{If }$\left( E,\pi ,M\right) =\left( F,\nu
,N\right) $\emph{, then, for any }
\begin{equation*}
\begin{array}{c}
u^{a}s_{a}\in \Gamma \left( E,\pi ,M\right) ,%
\end{array}%
\end{equation*}%
\emph{\ we shall use the notation }%
\begin{equation*}
\begin{array}{c}
u_{~|a|b}^{a_{1}}=u_{~|ab}^{a_{1}}-\rho \Gamma _{ab}^{d}u_{~|d}^{a_{1}}%
\end{array}%
\leqno\left( 5.11\right)
\end{equation*}%
\emph{and we verify the formulas of Ricci type }%
\begin{equation*}
\begin{array}{c}
u_{~|a|b}^{a_{1}}-u_{~|b|a}^{a_{1}}+\left( \rho ,h\right) \mathbb{T}%
_{~ab}^{d}u_{~|d}^{a_{1}}=u^{d}\left( \rho ,h\right) \mathbb{R}%
_{d~ab}^{a_{1}}-u_{~|c}^{a_{1}}L_{ab}^{c}\circ h%
\end{array}%
\leqno\left( 5.12\right)
\end{equation*}

\emph{In the particular case of Lie algebroids, }$h=Id_{M}$\emph{, the
relations }$\left( 5.12\right) $\emph{\ become}
\begin{equation*}
\begin{array}{c}
u_{~|a|b}^{a_{1}}-u_{~|b|a}^{a_{1}}+\rho \mathbb{T}%
_{~ab}^{d}u_{~|d}^{a_{1}}=u^{d}\rho \mathbb{R}%
_{d~ab}^{a_{1}}-u_{~|c}^{a_{1}}L_{ab}^{c}%
\end{array}%
\leqno\left( 5.12^{\prime }\right)
\end{equation*}

\emph{In the classical case,}$\left( \rho ,h\right) =\left(
Id_{TM},Id_{M}\right) ,$\emph{\ the relations }$\left( 5.12^{\prime }\right)
$\emph{\ become }%
\begin{equation*}
\begin{array}{c}
u_{~|i|j}^{i_{1}}-u_{~|i|j}^{i_{1}}+\mathbb{T}_{~ij}^{k}u_{~|k}^{i_{1}}=u^{k}%
\mathbb{R}_{k~ij}^{i_{1}}%
\end{array}%
\leqno\left( 5.12^{\prime \prime }\right)
\end{equation*}

\textbf{Theorem 5.3 } \emph{For any }$u_{a}s^{a}\in \Gamma \left( \overset{%
\ast }{E},\overset{\ast }{\pi },M\right) $ \emph{we shall use the notation}
\begin{equation*}
\begin{array}{c}
u_{b_{1}\mid \alpha \beta }=\rho _{\beta }^{j}\circ h\frac{\partial }{%
\partial x^{j}}\left( u_{b_{1}\mid \alpha }\right) -\rho \Gamma _{b_{1}\beta
}^{b}u_{b\mid \alpha }%
\end{array}%
\leqno\left( 5.13\right)
\end{equation*}%
\emph{and we verify the formulas: }%
\begin{equation*}
\begin{array}{c}
u_{b_{1}\mid \alpha \beta }-u_{b_{1}\mid \beta \alpha }=-u_{b}\left( \rho
,h\right) \mathbb{R}_{b_{1}~\alpha \beta }^{b}-u_{b_{1}\mid \gamma
}L_{\alpha \beta }^{\gamma }\circ h%
\end{array}%
\leqno\left( 5.14\right)
\end{equation*}

\emph{After some calculations, we obtain}%
\begin{equation*}
\begin{array}{c}
\left( \rho ,h\right) \mathbb{R}_{b_{1}~\alpha \beta }^{b}=u^{b}\left(
-u_{b_{1}\mid \alpha \beta }+u_{b_{1}\mid \beta \alpha }-u_{b_{1}\mid \gamma
}L_{\alpha \beta }^{\gamma }\circ h\right) ,%
\end{array}%
\leqno\left( 5.15\right)
\end{equation*}%
\emph{where} $u^{a}s_{a}\in \Gamma \left( E,\pi ,M\right) $\emph{\ such that
}$u_{a}u^{b}=\delta _{a}^{b}.$

\emph{In the particular case of Lie algebroids, }$h=Id_{M}$\emph{, the
relations }$\left( 5.15\right) $\emph{\ become}%
\begin{equation*}
\begin{array}{c}
\rho \mathbb{R}_{b_{1}~\alpha \beta }^{b}=u^{b}\left( -u_{b_{1}\mid \alpha
\beta }+u_{b_{1}\mid \beta \alpha }-u_{b_{1}\mid \gamma }L_{\alpha \beta
}^{\gamma }\right) .%
\end{array}%
\leqno\left( 5.15^{\prime }\right)
\end{equation*}

\emph{In the classical case, }$\left( \rho ,h\right) =\left(
Id_{TM},Id_{M}\right) ,$ \emph{the relations }$\left( 5.15^{\prime }\right) $%
\emph{\ become}%
\begin{equation*}
\begin{array}{c}
\mathbb{R}_{b_{1}~ij}^{b}=u^{b}\left( -u_{b_{1}\mid ij}+u_{b_{1}\mid
ji}\right) .%
\end{array}%
\leqno(5.15^{\prime \prime })
\end{equation*}

\emph{Proof.} Since
\begin{equation*}
\begin{split}
u_{b_{1}\mid \alpha \beta }& =\rho _{\beta }^{j}\circ h\left( \frac{\partial
}{\partial x^{j}}\left( \rho _{\alpha }^{i}\circ h\frac{\partial u_{b_{1}}}{%
\partial x^{i}}-\rho \Gamma _{b_{1}\alpha }^{b}u_{b}\right) \right) \\
& -\rho \Gamma _{b_{1}\beta }^{b}\left( \rho _{\alpha }^{i}\circ h\frac{%
\partial u_{b}}{\partial x^{i}}-\rho \Gamma _{b\alpha }^{a}u_{a}\right) \\
& =\rho _{\beta }^{j}\circ h\frac{\partial \rho _{\alpha }^{i}\circ h}{%
\partial x^{j}}\frac{\partial u_{b_{1}}}{\partial x^{i}}+\rho _{\beta
}^{j}\circ h\rho _{\alpha }^{i}\circ h\frac{\partial }{\partial x^{j}}\left(
\frac{\partial u_{b_{1}}}{\partial x^{i}}\right) \\
& -\rho _{\beta }^{j}\circ h\frac{\partial \rho \Gamma _{b_{1}\alpha }^{b}}{%
\partial x^{j}}u_{b}-\rho _{\beta }^{j}\circ h\rho \Gamma _{b_{1}\alpha }^{b}%
\frac{\partial u_{b}}{\partial x^{j}} \\
& -\rho _{\alpha }^{i}\circ h\rho \Gamma _{b_{1}\beta }^{b}\frac{\partial
u_{b}}{\partial x^{i}}+\rho \Gamma _{b_{1}\beta }^{b}\rho \Gamma _{b\alpha
}^{a}u_{a}
\end{split}%
\end{equation*}%
and%
\begin{equation*}
\begin{split}
u_{b_{1}\mid \beta \alpha }& =\rho _{\alpha }^{i}\circ h\left( \frac{%
\partial }{\partial x^{i}}\left( \rho _{\beta }^{j}\circ h\frac{\partial
u_{b_{1}}}{\partial x^{j}}-\rho \Gamma _{b_{1}\beta }^{b}u_{b}\right) \right)
\\
& -\rho \Gamma _{b_{1}\alpha }^{b}\left( \rho _{\beta }^{j}\circ h\frac{%
\partial u_{b}}{\partial x^{j}}-\rho \Gamma _{b\beta }^{a}u_{a}\right) \\
& =\rho _{\alpha }^{i}\circ h\frac{\partial \rho _{\beta }^{j}\circ h}{%
\partial x^{i}}\frac{\partial u_{b_{1}}}{\partial x^{i}}+\rho _{\beta
}^{j}\circ h\rho _{\alpha }^{i}\circ h\frac{\partial }{\partial x^{i}}\left(
\frac{\partial u_{b_{1}}}{\partial x^{j}}\right) \\
& -\rho _{\alpha }^{i}\circ h\frac{\partial \rho \Gamma _{b_{1}\beta }^{b}}{%
\partial x^{i}}u_{b}-\rho _{\alpha }^{i}\circ h\rho \Gamma _{b_{1}\beta }^{b}%
\frac{\partial u_{b}}{\partial x^{i}} \\
& -\rho _{\beta }^{j}\circ h\rho \Gamma _{b_{1}\alpha }^{b}\frac{\partial
u_{b}}{\partial x^{i}}+\rho \Gamma _{b_{1}\alpha }^{b}\rho \Gamma _{b\beta
}^{a}u_{a}
\end{split}%
\end{equation*}%
it results that%
\begin{equation*}
\begin{split}
\ u_{b_{1}\mid \alpha \beta }-u_{b_{1}\mid \beta \alpha }& =\rho _{\beta
}^{j}\circ h\frac{\partial \rho _{\alpha }^{i}\circ h}{\partial x^{j}}\frac{%
\partial u_{b_{1}}}{\partial x^{i}}-\rho _{\alpha }^{i}\circ h\frac{\partial
\rho _{\beta }^{j}\circ h}{\partial x^{i}}\frac{\partial u_{b_{1}}}{\partial
x^{j}} \\
& +\rho _{\beta }^{j}\circ h\rho _{\alpha }^{i}\circ h\frac{\partial }{%
\partial x^{j}}\left( \frac{\partial u_{b_{1}}}{\partial x^{i}}\right) -\rho
_{\beta }^{j}\circ h\rho _{\alpha }^{i}\circ h\frac{\partial }{\partial x^{i}%
}\left( \frac{\partial u_{b_{1}}}{\partial x^{j}}\right) \\
& +\rho _{\alpha }^{i}\circ h\frac{\partial \rho \Gamma _{b_{1}\beta }^{b}}{%
\partial x^{i}}u_{b}-\rho _{\beta }^{j}\circ h\frac{\partial \rho \Gamma
_{b_{1}\alpha }^{b}}{\partial x^{j}}u_{b} \\
& +\rho _{\beta }^{j}\circ h\rho \Gamma _{b_{1}\alpha }^{b}\frac{\partial
u_{b}}{\partial x^{j}}-\rho _{\beta }^{j}\circ h\rho \Gamma _{b_{1}\alpha
}^{b}\frac{\partial u_{b}}{\partial x^{j}} \\
& +\rho _{\alpha }^{i}\circ h\rho \Gamma _{b_{1}\alpha }^{b}\frac{\partial
u_{b}}{\partial x^{i}}-\rho _{\alpha }^{i}\circ h\rho \Gamma _{b_{1}\alpha
}^{b}\frac{\partial u_{b}}{\partial x^{i}} \\
& +\rho \Gamma _{b_{1}\beta }^{b}\rho \Gamma _{b\alpha }^{a}u_{a}-\rho
\Gamma _{b_{1}\alpha }^{b}\rho \Gamma _{b\beta }^{a}u_{a}.
\end{split}%
\end{equation*}

After some calculations, we obtain:
\begin{equation*}
\begin{array}{cl}
u_{b_{1}\mid \alpha \beta }-u_{b_{1}\mid \beta \alpha } & =L_{\beta \alpha
}^{\gamma }\circ h\rho _{\gamma }^{k}\circ h\frac{\partial u_{b_{1}}}{%
\partial x^{k}}\vspace*{1mm} \\
& +\left( \rho _{\alpha }^{i}\circ h\frac{\partial \rho \Gamma _{b_{1}\beta
}^{b}}{\partial x^{i}}u_{b}-\rho _{\beta }^{j}\circ h\frac{\partial \rho
\Gamma _{b_{1}\alpha }^{b}}{\partial x^{j}}u_{b}\right) \vspace*{1mm} \\
& +\rho \Gamma _{b_{1}\beta }^{b}\rho \Gamma _{b\alpha }^{a}u_{a}-\rho
\Gamma _{b_{1}\alpha }^{b}\rho \Gamma _{b\beta }^{a}u_{a}.%
\end{array}%
\end{equation*}

Since
\begin{equation*}
\begin{split}
u_{b}\left( \rho ,h\right) \mathbb{R}_{b_{1}\alpha \beta }^{b}& =u_{b}\left(
\rho _{\beta }^{j}\circ h\frac{\partial \rho \Gamma _{b_{1}\alpha }^{b}}{%
\partial x^{j}}+\rho \Gamma _{e\beta }^{b}\rho \Gamma _{b_{1}\alpha
}^{e}\right. \\
& -\rho _{\alpha }^{i}\circ h\frac{\partial \rho \Gamma _{b_{1}\beta }^{b}}{%
\partial x^{i}}\left. -\rho \Gamma _{e\alpha }^{b}\rho \Gamma _{b_{1}\beta
}^{e}-\rho \Gamma _{b_{1}\gamma }^{b}L_{\beta \alpha }^{\gamma }\circ
h\right)
\end{split}%
\end{equation*}%
and
\begin{equation*}
u_{b_{1}\mid \gamma }L_{\alpha \beta }^{\gamma }\circ h=\left( \rho _{\gamma
}^{k}\circ h\frac{\partial u_{b_{1}}}{\partial x^{k}}-\rho \Gamma
_{b_{1}\gamma }^{b}u_{b}\right) L_{\alpha \beta }^{\gamma }\circ h
\end{equation*}%
it results that 
\begin{equation*}
\begin{split}
-u_{b}\left( \rho ,h\right) \mathbb{R}_{b_{1},\alpha \beta
}^{b}-u_{b_{1}\mid \gamma }L_{\alpha \beta }^{\gamma }\circ h& =-L_{\alpha
\beta }^{\gamma }\circ h\rho _{\gamma }^{k}\circ h\frac{\partial u_{b_{1}}}{%
\partial x^{k}} \\
& +\left( \rho _{\alpha }^{i}\circ h\frac{\partial \rho \Gamma _{b_{1}\beta
}^{b}}{\partial x^{i}}u_{b}-\rho _{\beta }^{j}\circ h\frac{\partial \rho
\Gamma _{b_{1}\alpha }^{b}}{\partial x^{j}}u_{b}\right) \\
& +\rho \Gamma _{b_{1}\beta }^{b}\rho \Gamma _{b\alpha }^{a}u_{a}-\rho
\Gamma _{b_{1}\alpha }^{b}\rho \Gamma _{b\beta }^{a}u_{a}.
\end{split}%
\end{equation*}

\begin{flushright}
\emph{q.e.d.}
\end{flushright}

\textbf{Lemma 5.2 }\emph{If }$\left( E,\pi ,M\right) =\left( F,\nu ,N\right)
$\emph{, then, for any }
\begin{equation*}
\begin{array}{c}
u_{b}s^{b}\in \Gamma \left( \overset{\ast }{E},\overset{\ast }{\pi }%
,M\right) ,%
\end{array}%
\end{equation*}%
\emph{\ we have that }$u_{b~|c},~b,c\in \overline{1,n}$ are the \emph{%
components of a tensor field of }$\left( 0,2\right) $\emph{\ type.}

\emph{Proof.} Let $U$ and $U^{\prime }$ be two vector local $\left(
m+n\right) $-charts such that $U\cap U^{\prime }\neq \phi .$

Since $u_{b^{\prime }}\left( x\right) =M_{b^{\prime }}^{b}\left( x\right)
u_{b}\left( x\right) ,~$for any $x\in U\cap U^{\prime },$ it results that
\begin{equation*}
\begin{array}{cl}
\displaystyle\rho _{c^{\prime }}^{k^{\prime }}\circ h\left( x\right) \frac{%
\partial u_{b^{\prime }}\left( x\right) }{\partial x^{k^{\prime }}} & =%
\displaystyle\rho _{c^{\prime }}^{k^{\prime }}\circ h\left( x\right) \frac{%
\partial }{\partial x^{k^{\prime }}}\left( M_{b^{\prime }}^{b}\left(
x\right) \right) u_{b}\left( x\right) \vspace*{1mm} \\
& \displaystyle+M_{b^{\prime }}^{b}\left( x\right) \rho _{c^{\prime
}}^{k^{\prime }}\circ h\left( x\right) \frac{\partial u_{b}\left( x\right) }{%
\partial x^{k^{\prime }}}.%
\end{array}%
\leqno(1)
\end{equation*}

Since, for any $x\in U\cap U^{\prime }$, we have
\begin{equation*}
\begin{array}{cl}
\displaystyle\rho \Gamma _{b^{\prime }c^{\prime }}^{a^{\prime }}\left(
x\right) & \displaystyle=M_{a}^{a^{\prime }}\left( x\right) \left( \rho
_{c}^{k}\circ h\left( x\right) \frac{\partial }{\partial x^{k}}\left(
M_{b^{\prime }}^{a}\left( x\right) \right) \right. \vspace*{1mm} \\
& \displaystyle\left. +\rho \Gamma _{bc}^{a}\left( x\right) M_{b^{\prime
}}^{b}\left( x\right) \right) M_{c^{\prime }}^{c}\left( x\right) ,%
\end{array}%
\leqno(2)
\end{equation*}%
and
\begin{equation*}
\begin{array}{cl}
0 & \displaystyle=\frac{\partial }{\partial x^{k^{\prime }}}\left(
M_{a}^{a^{\prime }}\left( x\right) M_{b^{\prime }}^{a}\left( x\right)
\right) \vspace*{1mm} \\
& \displaystyle=\frac{\partial }{\partial x^{k^{\prime }}}\left(
M_{a}^{a^{\prime }}\left( x\right) \right) M_{b^{\prime }}^{a}\left(
x\right) +M_{a}^{a^{\prime }}\left( x\right) \frac{\partial }{\partial
x^{k^{\prime }}}\left( M_{b^{\prime }}^{a}\left( x\right) \right)%
\end{array}%
\leqno(3)
\end{equation*}%
it results that
\begin{equation*}
\begin{array}{c}
\displaystyle\rho \Gamma _{b^{\prime }c^{\prime }}^{a^{\prime }}\left(
x\right) u_{a^{\prime }}\left( x\right) =-\rho _{c^{\prime }}^{k^{\prime
}}\circ h\left( x\right) \frac{\partial }{\partial x^{k^{\prime }}}\left(
M_{b^{\prime }}^{b}\left( x\right) \right) u_{b}\left( x\right) \vspace*{1mm}
\\
\displaystyle+M_{b^{\prime }}^{b}\left( x\right) \rho \Gamma _{bc}^{a}\left(
x\right) u_{a}\left( x\right) M_{c^{\prime }}^{c}\left( x\right) .%
\end{array}%
\leqno(4)
\end{equation*}

Summing the equalities $\left( 1\right) $ and $\left( 4\right) $, it results
the conclusion of lemma.\hfill \emph{q.e.d.}

\textbf{Theorem 5.4 }\emph{If }$\left( E,\pi ,M\right) =\left( F,\nu
,N\right) $\emph{, then, for any }
\begin{equation*}
\begin{array}{c}
u_{b}s^{b}\in \Gamma \left( \overset{\ast }{E},\overset{\ast }{\pi }%
,M\right) ,%
\end{array}%
\end{equation*}%
\emph{\ we shall use the notation }%
\begin{equation*}
\begin{array}{c}
u_{b_{1}~|a|b}=u_{b_{1}~|ab}-\rho \Gamma _{ab}^{d}u_{b_{1}~|d}%
\end{array}%
\leqno(5.16)
\end{equation*}%
\emph{and we verify the formulas of Ricci type }%
\begin{equation*}
\begin{array}{c}
u_{b_{1}~|a|b}-u_{b_{1}~|b|a}+\left( \rho ,h\right) \mathbb{T}%
_{~ab}^{d}u_{b_{1}~|d}=-u_{d}\left( \rho ,h\right) \mathbb{R}%
_{b_{1}~ab}^{d}-u_{b_{1}~|d}L_{ab}^{d}\circ h%
\end{array}%
\leqno(5.17)
\end{equation*}

\emph{In the particular case of Lie algebroids, }$h=Id_{M}$\emph{, the
relations }$\left( 5.17\right) $\emph{\ become}%
\begin{equation*}
\begin{array}{c}
u_{b_{1}~|a|b}-u_{b_{1}~|b|a}+\rho \mathbb{T}_{~ab}^{d}u_{b_{1}~|d}=-u_{d}%
\rho \mathbb{R}_{b_{1}~ab}^{d}-u_{b_{1}~|d}L_{ab}^{d}.%
\end{array}%
\leqno(5.17^{\prime })
\end{equation*}

\emph{In the classical case, }$\left( \rho ,h\right) =\left(
Id_{TM},Id_{M}\right) ,$ \emph{the relations }$\left( 5.17^{\prime }\right) $%
\emph{\ become}%
\begin{equation*}
\begin{array}{c}
u_{j_{1}~|i|j}-u_{j_{1}~|j|i}+\mathbb{T}_{~ij}^{h}u_{j_{1}~|h}=u_{h}\mathbb{R%
}_{j_{1}~ij}^{h}.%
\end{array}%
\leqno(5.17^{\prime \prime })
\end{equation*}

\textbf{Theorem 5.5 }\emph{For any tensor field }%
\begin{equation*}
\begin{array}{c}
T_{b_{1}...b_{q}}^{a_{1}...a_{p}}s_{a_{1}}\otimes ...\otimes
s_{a_{p}}\otimes s^{b_{1}}\otimes ...\otimes s^{b_{q}}%
\end{array}%
,
\end{equation*}%
\emph{we verify the equality: }%
\begin{equation*}
\begin{array}{l}
T_{b_{1}...b_{q}\mid \alpha \beta }^{a_{1}...a_{p}}-T_{b_{1}...b_{q}\mid
\beta \alpha }^{a_{1}...a_{p}}=T_{b_{1}...b_{q}}^{aa_{2}...a_{p}}\left( \rho
,h\right) \mathbb{R}_{a~\alpha \beta }^{a_{1}}+...\vspace*{1mm} \\
+T_{b_{1}...b_{q}}^{a_{1}...a_{p-1}a}\left( \rho ,h\right) \mathbb{R}%
_{a~\alpha \beta }^{a_{p}}-T_{bb_{2}...b_{q}}^{a_{1}...a_{p}}\left( \rho
,h\right) \mathbb{R}_{b_{1}~\alpha \beta }^{b}-...\vspace*{1mm} \\
-T_{b_{1}...b_{q-1}b}^{a_{1}...a_{p}}\left( \rho ,h\right) \mathbb{R}%
_{b_{q}~\alpha \beta }^{b}-T_{b_{1}...b_{q}\mid \gamma
}^{a_{1}...a_{p}}L_{\alpha \beta }^{\gamma }\circ h.%
\end{array}%
\leqno\left( 5.18\right)
\end{equation*}

\emph{In the particular case of Lie algebroids, }$h=Id_{M}$\emph{, the
relations }$\left( 5.18\right) $\emph{\ become}%
\begin{equation*}
\begin{array}{l}
T_{b_{1}...b_{q}\mid \alpha \beta }^{a_{1}...a_{p}}-T_{b_{1}...b_{q}\mid
\beta \alpha }^{a_{1}...a_{p}}=T_{b_{1}...b_{q}}^{aa_{2}...a_{p}}\rho
\mathbb{R}_{a~\alpha \beta }^{a_{1}}+...\vspace*{1mm} \\
+T_{b_{1}...b_{q}}^{a_{1}...a_{p-1}a}\rho \mathbb{R}_{a~\alpha \beta
}^{a_{p}}-T_{bb_{2}...b_{q}}^{a_{1}...a_{p}}\rho \mathbb{R}_{b_{1}~\alpha
\beta }^{b}-...\vspace*{1mm} \\
-T_{b_{1}...b_{q-1}b}^{a_{1}...a_{p}}\rho \mathbb{R}_{b_{q}~\alpha \beta
}^{b}-T_{b_{1}...b_{q}\mid \gamma }^{a_{1}...a_{p}}L_{\alpha \beta }^{\gamma
}.%
\end{array}%
\leqno(5.18^{\prime })
\end{equation*}

\emph{In the classical case,} $\left( \rho ,h\right) =\left(
Id_{TM},Id_{M}\right) ,$ \emph{the relations }$\left( 5.18^{\prime }\right) $%
\emph{\ become}%
\begin{equation*}
\begin{array}{l}
T_{j_{1}...j_{q}\mid hk}^{i_{1}...i_{p}}-T_{j_{1}...j_{q}\mid
kh}^{i_{1}...i_{p}}=T_{j_{1}...j_{q}}^{ii_{2}...i_{p}}\mathbb{R}%
_{i~hk}^{i_{1}}+...\vspace*{1mm} \\
+T_{j_{1}...j_{q}}^{i_{1}...i_{p-1}i}\mathbb{R}%
_{i~hk}^{i_{p}}-T_{jj_{2}...j_{q}}^{i_{1}...i_{p}}\mathbb{R}%
_{j_{1}~hk}^{j}-...\vspace*{1mm}-T_{j_{1}...j_{q-1}j}^{i_{1}...i_{p}}\mathbb{%
R}_{j_{q}~hk}^{j}.%
\end{array}%
\leqno(5.18^{\prime \prime })
\end{equation*}

\textbf{Theorem 5.6 }\emph{If }$\left( E,\pi ,M\right) =\left( F,\nu
,N\right) $\emph{, then we obtain the following formulas of Ricci type:}%
\begin{equation*}
\begin{array}{l}
T_{b_{1}...b_{q}~|b|c}^{a_{1}...a_{p}}-T_{b_{1}...b_{q}\mid
c|b}^{a_{1}...a_{p}}+\left( \rho ,h\right) \mathbb{T}%
_{~bc}^{d}T_{b_{1}...b_{q}\mid d}^{a_{1}...a_{p}}\vspace*{1mm} \\
=T_{b_{1}...b_{q}}^{aa_{2}...a_{p}}\left( \rho ,h\right) \mathbb{R}%
_{a~bc}^{a_{1}}+...+T_{b_{1}...b_{q}}^{a_{1}...a_{p-1}a}\left( \rho
,h\right) \mathbb{R}_{a~bc}^{a_{p}}\vspace*{1mm} \\
-T_{bb_{2}...b_{q}}^{a_{1}...a_{p}}\left( \rho ,h\right) \mathbb{R}%
_{b_{1}~bc}^{b}-...-T_{b_{1}...b_{q-1}b}^{a_{1}...a_{p}}\left( \rho
,h\right) \mathbb{R}_{b_{q}~bc}^{b}-T_{b_{1}...b_{q}\mid
d}^{a_{1}...a_{p}}L_{bc}^{d}\circ h.%
\end{array}%
\leqno(5.19)
\end{equation*}

\emph{In the particular case of Lie algebroids, }$h=Id_{M}$\emph{, the
relations }$\left( 5.19\right) $\emph{\ become}%
\begin{equation*}
\begin{array}{l}
T_{b_{1}...b_{q}~|b|c}^{a_{1}...a_{p}}-T_{b_{1}...b_{q}\mid
c|b}^{a_{1}...a_{p}}+\rho \mathbb{T}_{~bc}^{d}T_{b_{1}...b_{q}\mid
d}^{a_{1}...a_{p}}\vspace*{1mm} \\
=T_{b_{1}...b_{q}}^{aa_{2}...a_{p}}\rho \mathbb{R}%
_{a~bc}^{a_{1}}+...+T_{b_{1}...b_{q}}^{a_{1}...a_{p-1}a}\rho \mathbb{R}%
_{a~bc}^{a_{p}}\vspace*{1mm} \\
-T_{bb_{2}...b_{q}}^{a_{1}...a_{p}}\rho \mathbb{R}%
_{b_{1}~bc}^{b}-...-T_{b_{1}...b_{q-1}b}^{a_{1}...a_{p}}\rho \mathbb{R}%
_{b_{q}~bc}^{b}-T_{b_{1}...b_{q}\mid d}^{a_{1}...a_{p}}L_{bc}^{d}.%
\end{array}%
\leqno(5.19^{\prime })
\end{equation*}

\emph{In the classical case,} $\left( \rho ,h\right) =\left(
Id_{TM},Id_{M}\right) ,$ \emph{the relations }$\left( 5.19^{\prime }\right) $%
\emph{\ become}%
\begin{equation*}
\begin{array}{l}
T_{j_{1}...j_{q}~|h|k}^{i_{1}...i_{p}}-T_{j_{1}...j_{q}\mid
k|h}^{i_{1}...i_{p}}+\mathbb{T}_{~hk}^{m}T_{j_{1}...j_{q}\mid
m}^{i_{1}...i_{p}}\vspace*{1mm} \\
=T_{j_{1}...j_{q}}^{ii_{2}...i_{p}}\mathbb{R}%
_{i~hk}^{i_{1}}+...+T_{j_{1}...j_{q}}^{i_{1}...i_{p-1}i}\mathbb{R}%
_{i~hk}^{i_{p}}\vspace*{1mm}-T_{jj_{2}...j_{q}}^{i_{1}...i_{p}}\mathbb{R}%
_{j_{1}~hk}^{j}-...-T_{j_{1}...j_{q-1}j}^{i_{1}...i_{p}}\mathbb{R}%
_{j_{q}~hk}^{j}.%
\end{array}%
\leqno(5.19^{\prime \prime })
\end{equation*}

We observe that if the structure functions of generalized Lie algebroid
\begin{equation*}
\left( \left( F,\nu ,M\right) ,\left[ ,\right] _{F,h},\left( \rho
,Id_{M}\right) \right) ,
\end{equation*}%
the $\left( \rho ,h\right) $-torsion associated to linear $\rho $-connection
$\rho \Gamma $ and the $\left( \rho ,h\right) $-curvature asso\-cia\-ted to
linear $\rho $-connection $\rho \Gamma $ are null, then we have the
equality:
\begin{equation*}
\begin{array}{c}
T_{b_{1}...b_{q}\mid b|c}^{a_{1}...a_{p}}=T_{b_{1}...b_{q}\mid
c|b}^{a_{1}...a_{p}},%
\end{array}%
\leqno(5.20)
\end{equation*}%
which generalizes the Schwartz equality.

\textbf{Theorem 5.7}{\ }\emph{If }$(E,\pi ,M){=}(F,\nu ,N)$\emph{, then the
following relations hold good}
\begin{equation*}
\begin{array}{r}
\underset{cyclic\left( U_{1},U_{2},U_{3}\right) }{\dsum }\left\{ \left( \rho
\ddot{D}_{U_{1}}\left( \rho ,h\right) \mathbb{T}\right) \left(
U_{2},U_{3}\right) -\left( \rho ,h\right) \mathbb{R}\left(
U_{1},U_{2}\right) U_{3}\right. \vspace*{1mm} \\
\left. +\left( \rho ,h\right) \mathbb{T}\left( \left( \rho ,h\right) \mathbb{%
\ T}\left( U_{1},U_{2}\right) ,U_{3}\right) \right\} =0,%
\end{array}%
\leqno(\tilde{B}_{1})
\end{equation*}%
\emph{and}
\begin{equation*}
\underset{cyclic\left( U_{1},U_{2},U_{3},U\right) }{\dsum }\!\!\left\{
\left( \rho \ddot{D}_{U_{1}}\left( \rho ,h\right) \mathbb{R}\right) \left(
U_{2},U_{3}\right) U{+}\left( \rho ,h\right) \mathbb{R}\left( \left( \rho
,h\right) \mathbb{T}\left( U_{1},U_{2}\right) ,U_{3}\right) U\right\} {=}0.%
\leqno(\tilde{B}_{2})
\end{equation*}%
\emph{respectively. This identities will be called the first respectively
the second identity of Bianchi type.}

\emph{In the particular case of Lie algebroids, }$h=Id_{M},$\emph{\ the
identities }$(\tilde{B}_{1})$\emph{\ and }$(\tilde{B}_{2})$\emph{\ become }%
\begin{equation*}
\underset{cyclic\left( u_{1},u_{2},u_{3}\right) }{\dsum }\left\{ \left( \rho
D_{u_{1}}\rho \mathbb{T}\right) \left( u_{2},u_{3}\right) -\rho \mathbb{R}%
\left( u_{1},u_{2}\right) u_{3}+\rho \mathbb{T}\left( \rho \mathbb{\ T}%
\left( u_{1},u_{2}\right) ,u_{3}\right) \right\} =0,\leqno(\tilde{B}%
_{1}^{\prime })
\end{equation*}%
\begin{equation*}
\begin{array}{c}
\underset{cyclic\left( u_{1},u_{2},u_{3},u\right) }{\dsum }\left\{ \left(
\rho D_{u_{1}}\rho \mathbb{R}\right) \left( u_{2},u_{3}\right) u+\rho
\mathbb{R}\left( \rho \mathbb{T}\left( u_{1},u_{2}\right) ,u_{3}\right)
u\right\} =0.%
\end{array}%
\leqno(\tilde{B}_{2}^{\prime })
\end{equation*}

\emph{In the classical case,} $\left( \rho ,h\right) =\left(
Id_{TM},Id_{M}\right) ,$\emph{\ the identities }$(\tilde{B}_{1}^{\prime })$%
\emph{\ and }$(\tilde{B}_{2}^{\prime })$\emph{\ become }%
\begin{equation*}
\underset{cyclic\left( X_{1},X_{2},X_{3}\right) }{\dsum }\left\{ \left(
D_{X_{1}}\mathbb{T}\right) \left( X_{2},X_{3}\right) -\mathbb{R}\left(
X_{1},X_{2}\right) X_{3}+\mathbb{T}\left( \mathbb{\ T}\left(
X_{1},X_{2}\right) ,X_{3}\right) \right\} =0,\leqno(\tilde{B}_{1}^{\prime
\prime })
\end{equation*}%
\begin{equation*}
\begin{array}{c}
\underset{cyclic\left( X_{1},X_{2},X_{3},X\right) }{\dsum }\left\{ \left(
D_{X_{1}}\mathbb{R}\right) \left( X_{2},X_{3}\right) X+\mathbb{R}\left(
\mathbb{T}\left( X_{1},X_{2}\right) ,X_{3}\right) X\right\} =0.%
\end{array}%
\leqno(\tilde{B}_{2}^{\prime \prime })
\end{equation*}

\emph{Proof: }Using the equality
\begin{equation*}
\begin{array}[b]{l}
\left( \rho \ddot{D}_{U_{1}}\left( \rho ,h\right) \mathbb{T}\right) \left(
U_{2},U_{3}\right) =\rho \ddot{D}_{U_{1}}\left( \left( \rho ,h\right)
\mathbb{T}\left( U_{2},U_{3}\right) \right) \\
-\left( \rho ,h\right) \mathbb{T}\left( \rho \ddot{D}_{U_{1}}U_{2},U_{3}%
\right) -\left( \rho ,h\right) \mathbb{T}\left( U_{2},\rho \ddot{D}%
_{U_{1}}U_{3}\right)%
\end{array}%
\end{equation*}%
and the Jacobi identity we obtain the first identity of Bianchi type.

Using the equality%
\begin{equation*}
\begin{array}[b]{l}
\left( \rho \ddot{D}_{U_{1}}\left( \rho ,h\right) \mathbb{R}\right) \left(
U_{2},U_{3}\right) U=\rho \ddot{D}_{U_{1}}\left( \left( \rho ,h\right)
\mathbb{R}\left( U_{2},U_{3}\right) U\right) \\
-\left( \rho ,h\right) \mathbb{R}\left( \rho \ddot{D}_{U_{1}}U_{2},U_{3}%
\right) U-\left( \rho ,h\right) \mathbb{R}\left( U_{2},\rho \ddot{D}%
_{U_{1}}U_{3}\right) U-\left( \rho ,h\right) \mathbb{R}\left(
U_{2},U_{3}\right) \rho \ddot{D}_{U_{1}}U%
\end{array}%
\end{equation*}%
and the Jacobi identity we obtain the second identity of Bianchi type.\hfill
\emph{q.e.d.}

\emph{Remark 5.1} On components, the identities of Bianchi type become:
\begin{equation*}
\begin{array}{r}
\underset{cyclic\left( a_{1},a_{2},a_{3}\right) }{\sum }\left\{ \left( \rho
,h\right) \mathbb{T}_{~\ a_{2}a_{3}{\mid a_{1}}}^{b}+\left( \rho ,h\right)
\mathbb{T}_{~ga_{3}}^{b}\cdot \left( \rho ,h\right) \mathbb{T}%
_{~a_{1}a_{2}}^{g}\right\} \vspace*{1mm} \\
=\underset{cyclic\left( a_{1},a_{2},a_{3}\right) }{\sum }\left( \rho
,h\right) \mathbb{R}_{a_{3}~a_{1}a_{2}}^{a}%
\end{array}%
\end{equation*}%
and
\begin{equation*}
\underset{cyclic\left( a,a_{1},a_{2},a_{3}\right) }{\sum }\left\{ \left(
\rho ,h\right) \mathbb{R}_{a~a_{2}a_{3}{\mid a_{1}}}^{b}+\left( \rho
,h\right) \mathbb{R}_{a~ga_{3}}^{b}\cdot \left( \rho ,h\right) \mathbb{T}%
_{~a_{2}a_{1}}^{g}\right\} =0.
\end{equation*}

If the $\left( \rho ,h\right) $-torsion is null, then the identities of
Bianchi type become:
\begin{equation*}
\underset{cyclic\left( a_{1}~a_{2},a_{3}\right) }{\sum }\left( \rho
,h\right) \mathbb{R}_{a_{3}~a_{1}a_{2}}^{b}=0
\end{equation*}%
and
\begin{equation*}
\underset{cyclic\left( a,a_{1},a_{2},a_{3}\right) }{\sum }\left( \rho
,h\right) \mathbb{R}_{a~a_{2}a_{3}\mid _{a_{1}}}^{b}=0.
\end{equation*}

\section{(Pseudo)metrizable vector bundles. Formulas of Levi-Civita type}

We will apply our theory for the diagram:
\begin{equation*}
\begin{array}{c}
\xymatrix{E\ar[d]_\pi&\left( F,\left[ , \right] _{F,h},\left( \rho
,Id_{N}\right) \right)\ar[d]^\nu \\ M\ar[r]^h&N}%
\end{array}%
,\leqno(6.1)
\end{equation*}%
where $\left( E,\pi ,M\right) \in \left\vert \mathbf{B}^{\mathbf{v}%
}\right\vert $ and $\left( \left( F,\nu ,N\right) ,\left[ ,\right]
_{F,h},\left( \rho ,Id_{N}\right) \right) $ is a generalized Lie algebroid.

\textbf{Definition 6.1} We will say that the \emph{vector bundle }$\left(
E,\pi ,M\right) $\emph{\ is endowed with a pseudometrical structure}\textit{%
\ }if it exists $g=g_{ab}s^{a}\otimes s^{b}\in \mathcal{T}_{2}^{0}\left(
E,\pi ,M\right) $ such that for each $x\in M,$ the matrix $\left\Vert
g_{ab}\left( x\right) \right\Vert $ is nondegenerate and symmetric.

Moreover, if for each $x\in M$ the matrix $\left\Vert g_{ab}\left( x\right)
\right\Vert $ has constant signature, then we will say that \emph{the vector
bundle }$\left( E,\pi ,M\right) $\emph{\ is endowed with a metrical
structure.}

If $g=g_{ab}s^{a}\otimes s^{b}\in \mathcal{T}_{2}^{0}\left( E,\pi ,M\right) $
is a (pseudo) metrical structure, then, for any $a,b\in \overline{1,r}$ and
for any vector local $\left( m+r\right) $-chart $\left( U,s_{U}\right) $ of $%
\left( E,\pi ,M\right) $, we consider the real functions
\begin{equation*}
\begin{array}{ccc}
U & ^{\underrightarrow{~\ \ \tilde{g}^{ba}~\ \ }} & \mathbb{R}%
\end{array}%
\end{equation*}%
such that $\left\Vert \tilde{g}^{ba}\left( x\right) \right\Vert =\left\Vert
g_{ab}\left( x\right) \right\Vert ^{-1},~$for any $\forall x\in U.$

\textbf{Definition 6.2} We admit that $\left( E,\pi ,M\right) $ is a vector
bundle endowed with a (pseudo)metrical structure $g$ and with a linear $\rho
$-connection $\rho \Gamma .$

We will say that the \emph{linear }$\rho $\emph{-connection }$\rho \Gamma $%
\emph{\ is compatible with the (pseudo)metrical structure }$g$ if
\begin{equation*}
\begin{array}{c}
\rho D_{z}g=0,\ \forall z\in \Gamma \left( F,\nu ,N\right) .%
\end{array}%
\leqno(6.2)
\end{equation*}

\textbf{Definition 6.3} We will say that the vector bundle $(E,\pi ,M)$ is $%
\rho $\emph{-(pseudo)metrizable,} if it exists a (pseudo)metrical structure $%
g\in \mathcal{T}_{2}^{0}\left( E,\pi ,M\right) $ and a linear $\rho $%
-connection $\rho \Gamma $ for $\left( E,\pi ,M\right) $ compatible with $g.$
The $id_{TM}$-(pseudo)metri\-zable vector bundles will be called \textit{(}%
\textit{pseudo})\textit{metrizable vector bundles.}

In particular, if $\left( TM,\tau _{M},M\right) $ is a (pseudo)metri\-zable
vector bundle, then we will say that $\left( TM,\tau _{M},M\right) $ \textit{%
is a \emph{(pseudo)Riemannian space}, }and the manifold $M$ will be called (%
\textit{pseudo})\textit{Riemannian manifold.}

The linear connection of a (pseudo)Riemannian space will be called (\textit{%
pseudo})\textit{Rieman\-nian linear connection.}

\textbf{Theorem 6.1 }\emph{If }$(E,\pi ,M)=(F,\nu ,N)$\emph{\ and }$g\in
\mathcal{T}_{2}^{0}\left( h^{\ast }E,h^{\ast }\pi ,M\right) $\emph{\ is a
(pseudo)metri\-cal structure, then the local real functions}%
\begin{equation*}
\begin{array}{cl}
\rho \Gamma _{bc}^{a} & \displaystyle=\frac{1}{2}\tilde{g}^{ad}\left( \rho
_{c}^{k}\circ h\frac{\partial g_{bd}}{\partial x^{k}}+\ \rho _{b}^{j}\circ h%
\frac{\partial g_{dc}}{\partial x^{j}}-\rho _{d}^{l}\circ h\frac{\partial
g_{bc}}{\partial x^{l}}\right. \vspace*{1mm} \\
& \displaystyle\left. -\left( L_{bc}^{e}\circ h\right) g_{ed}-\left(
L_{bd}^{e}\circ h\right) g_{ec}+\left( L_{dc}^{e}\circ h\right)
g_{eb}\right) ,%
\end{array}%
\leqno(6.3)
\end{equation*}%
\emph{are the components of a linear }$\rho $\emph{-connection }$\rho \Gamma
$\emph{\ for the vector bundle }$\left( h^{\ast }E,h^{\ast }\pi ,M\right) $%
\emph{\ such that }$\left( \rho ,h\right) \mathbb{T}=0$\emph{\ and the
vector bundle }$\left( h^{\ast }E,h^{\ast }\pi ,M\right) $\emph{\ becomes} $%
\rho $\emph{-(pseudo)metrizable.}\textit{\ }This linear $\rho $-connection $%
\rho \Gamma $ will be called the \emph{linear }$\rho $\emph{-connection of
Levi-Civita type.}

\emph{In the particular case of Lie algebroids, }$h=Id_{M},$\emph{\ the
relations }$\left( 6.3\right) $\emph{\ become}%
\begin{equation*}
\begin{array}{cl}
\rho \Gamma _{bc}^{a} & \displaystyle=\frac{1}{2}\tilde{g}^{ad}\left( \rho
_{c}^{k}\frac{\partial g_{bd}}{\partial x^{k}}+\ \rho _{b}^{j}\frac{\partial
g_{dc}}{\partial x^{j}}-\rho _{d}^{l}\frac{\partial g_{bc}}{\partial x^{l}}%
-L_{bc}^{e}g_{ed}-L_{bd}^{e}g_{ec}+L_{dc}^{e}g_{eb}\right) .\vspace*{1mm}%
\end{array}%
\leqno(6.3^{\prime })
\end{equation*}

\emph{In the classical case, }$\left( \rho ,h\right) =\left(
Id_{TM},Id_{M}\right) ,$\emph{\ the relations }$\left( 6.3^{\prime }\right) $%
\emph{\ become}%
\begin{equation*}
\begin{array}{cl}
\Gamma _{jk}^{i} & \displaystyle=\frac{1}{2}\tilde{g}^{ih}\left( \frac{%
\partial g_{jh}}{\partial x^{k}}+\ \frac{\partial g_{hk}}{\partial x^{j}}-%
\frac{\partial g_{jk}}{\partial x^{h}}\right) .%
\end{array}%
\leqno(6.3^{\prime \prime })
\end{equation*}

\bigskip \noindent \textit{Proof.} Since
\begin{equation*}
\begin{array}{cl}
\left( \rho \ddot{D}_{U}g\right) V\otimes Z & =\Gamma \left( \overset{%
h^{\ast }E}{\rho },Id_{M}\right) \left( U\right) \left( \left( g\left(
V\otimes Z\right) \right) -g\left( \left( \rho \ddot{D}_{U}V\right) \otimes
Z\right) \right. \\
& -g\left( V\otimes \left( \rho \ddot{D}_{U}Z\right) \right) ,\vspace*{1mm}%
~\forall U,V,Z\in \Gamma \left( h^{\ast }E,h^{\ast }\pi ,M\right) .%
\end{array}%
\end{equation*}%
it results that, for any $U,V,Z\in \Gamma \left( h^{\ast }E,h^{\ast }\pi
,M\right) ,$ we obtain the equalities:
\begin{eqnarray*}
(1) &&\Gamma \left( \overset{h^{\ast }E}{\rho },Id_{M}\right) \left(
U\right) \left( g\left( V\otimes Z\right) \right) =g\left( \left( \rho \ddot{%
D}_{U}V\right) \otimes Z\right) +g\left( V\otimes \left( \rho \ddot{D}%
_{U}Z\right) \right) , \\
(2) &&\Gamma \left( \overset{h^{\ast }E}{\rho },Id_{M}\right) \left(
Z\right) \left( g\left( U\otimes V\right) \right) =g\left( \left( \rho \ddot{%
D}_{Z}U\right) \otimes V\right) +g\left( U\otimes \left( \rho \ddot{D}%
_{Z}V\right) \right) , \\
(3) &&\Gamma \left( \overset{h^{\ast }E}{\rho },Id_{M}\right) \left(
V\right) \left( g\left( Z\otimes U\right) \right) =g\left( \left( \rho \ddot{%
D}_{V}Z\right) \otimes U\right) +g\left( Z\otimes \left( \rho \ddot{D}%
_{V}U\right) \right) .
\end{eqnarray*}

We observe that $\left( 1\right) +\left( 3\right) -\left( 2\right) $ is
equivalent with the equality:
\begin{equation*}
\begin{array}{l}
g\left( \left( \rho \ddot{D}_{U}V+\rho \ddot{D}_{V}U\right) \otimes Z\right)
+g\left( \left( \rho \ddot{D}_{V}Z-\rho \ddot{D}_{Z}V\right) \otimes
U\right) \vspace*{1mm} \\
\qquad +g\left( \left( \rho \ddot{D}_{U}Z-\rho \ddot{D}_{Z}U\right) \otimes
V\right) =\Gamma \left( \overset{h^{\ast }E}{\rho },Id_{M}\right) \left(
U\right) \left( g\left( V\otimes Z\right) \right) \vspace*{1mm} \\
\qquad +\Gamma \left( \overset{h^{\ast }E}{\rho },Id_{M}\right) \left(
V\right) \left( g\left( Z\otimes U\right) \right) -\ \Gamma \left( \overset{%
h^{\ast }E}{\rho },Id_{M}\right) \left( Z\right) \left( g\left( U\otimes
V\right) \right) .%
\end{array}%
\end{equation*}

Using the condition $\left( \rho ,h\right) \mathbb{T}=0,$ which is
equivalent with the equality%
\begin{equation*}
\rho \ddot{D}_{U}V-\rho \ddot{D}_{V}U-\left[ U,V\right] _{h^{\ast }E}=0,
\end{equation*}%
we obtain the equality%
\begin{equation*}
\begin{array}{l}
2g\left( \left( \rho \ddot{D}_{U}V\right) \otimes Z\right) +g\left( \left[
V,U\right] _{h^{\ast }E}\otimes Z\right) +g\left( \left[ V,Z\right]
_{h^{\ast }E}\otimes U\right) +\vspace*{1mm}g\left( \left[ U,Z\right]
_{h^{\ast }E}\otimes V\right) \\
=\Gamma \left( \overset{h^{\ast }E}{\rho },Id_{M}\right) \left( U\right)
\left( g\left( V\otimes Z\right) \right) \vspace*{1mm}+\Gamma \left( \overset%
{h^{\ast }E}{\rho },Id_{M}\right) \left( V\right) \left( g\left( Z\otimes
U\right) \right) \\
-\Gamma \left( \overset{h^{\ast }E}{\rho },Id_{M}\right) \left( Z\right)
\left( g\left( U\otimes V\right) \right) ,~\forall U,V,Z\in \Gamma \left(
h^{\ast }E,h^{\ast }\pi ,M\right) .%
\end{array}%
\end{equation*}

This equality is equivalent with the following equality:
\begin{equation*}
\begin{array}{cl}
2g\left( \left( \rho \ddot{D}_{U}V\right) \otimes Z\right) & =\Gamma \left(
\overset{h^{\ast }E}{\rho },Id_{M}\right) \left( U\right) \cdot \left(
g\left( V\otimes Z\right) \right) \vspace*{1mm} \\
& +\Gamma \left( \overset{h^{\ast }E}{\rho },Id_{M}\right) \left( V\right)
\left( g\left( Z\otimes U\right) \right) -\Gamma \left( \overset{h^{\ast }E}{%
\rho },Id_{M}\right) \left( Z\right) \left( g\left( U\otimes V\right)
\right) \vspace*{1mm} \\
& +g\left( \left[ U,V\right] _{h^{\ast }E}\otimes Z\right) -g\left( \left[
V,Z\right] _{h^{\ast }E}\otimes U\right) +g\left( \left[ Z,U\right]
_{h^{\ast }E}\otimes V\right) \vspace*{1mm}%
\end{array}%
\end{equation*}%
for any $U,V,Z\in \Gamma \left( h^{\ast }E,h^{\ast }\pi ,M\right) .$

If $U=S_{c},V=S_{b}$ and $Z=S_{d},$ then we obtain the equality:%
\begin{equation*}
\begin{array}[b]{cl}
2g\left( \left( \rho \Gamma _{bc}^{a}S_{a}\right) \otimes S_{d}\right) &
=\rho _{c}^{k}\circ h\frac{\partial g\left( S_{b}\otimes S_{d}\right) }{%
\partial xk}+\ \rho _{b}^{j}\circ h\frac{\partial g\left( S_{d}\otimes
S_{c}\right) }{\partial x^{j}}-\rho _{d}^{l}\circ h\frac{\partial g\left(
S_{b}\otimes S_{c}\right) }{\partial x^{l}}\vspace*{1mm} \\
& +g\left( \left( L_{cb}^{e}\circ h\right) S_{e}\otimes S_{d}\right)
-g\left( \left( L_{bd}^{e}\circ h\right) S_{e}\otimes S_{c}\right) +g\left(
\left( L_{dc}^{e}\circ h\right) S_{e}\otimes S_{b}\right) ,%
\end{array}%
\end{equation*}
which is equivalent with:%
\begin{equation*}
\begin{array}[b]{cl}
2g_{da}\rho \Gamma _{bc}^{a} & =\rho _{c}^{k}\circ h\frac{\partial g_{bd}}{%
\partial xk}+\ \rho _{b}^{j}\circ h\frac{\partial g_{dc}}{\partial x^{j}}%
-\rho _{d}^{l}\circ h\frac{\partial g_{bc}}{\partial x^{l}} \\
& -\left( L_{bc}^{e}\circ h\right) g_{ed}-\left( L_{bd}^{e}\circ h\right)
g_{ec}+\left( L_{dc}^{e}\circ h\right) g_{eb}%
\end{array}%
\end{equation*}

Finally, we obtain:
\begin{equation*}
\begin{array}{cl}
\rho \Gamma _{bc}^{a} & =\frac{1}{2}\tilde{g}^{ad}\left( \rho _{c}^{k}\circ h%
\frac{\partial g_{bd}}{\partial xk}+\ \rho _{b}^{j}\circ h\frac{\partial
g_{dc}}{\partial x^{j}}-\rho _{d}^{l}\circ h\frac{\partial g_{bc}}{\partial
x^{l}}\right. \vspace*{1mm} \\
& \left. -\left( L_{bc}^{e}\circ h\right) g_{ed}-\left( L_{bd}^{e}\circ
h\right) g_{ec}+\left( L_{dc}^{e}\circ h\right) g_{eb}\right) ,%
\end{array}%
\end{equation*}%
where $\left\Vert \tilde{g}^{ad}\left( x\right) \right\Vert =\left\Vert
g_{da}\left( x\right) \right\Vert ^{-1},$ for any $x\in M.$ \hfill \emph{%
q.e.d.}

\textbf{Theorem 6.2. }\emph{If }$(E,\pi ,M)=(F,\nu ,N),$\emph{\ }$g\in
\mathcal{T}_{2}^{0}\left( h^{\ast }E,h^{\ast }\pi ,M\right) $\emph{\ is a
pseudo(metrical) structure and }$\mathbb{T}\in \mathcal{T}_{2}^{1}\left(
h^{\ast }E,h^{\ast }\pi ,M\right) $\emph{\ such that its components are skew
symmetric in the lover indices, then the local real functions }%
\begin{equation*}
\begin{array}{c}
\rho \mathring{\Gamma}_{bc}^{a}=\rho \Gamma _{bc}^{a}+\displaystyle\frac{1}{2%
}\tilde{g}^{ad}\left( g_{de}\mathbb{T}_{bc}^{e}-g_{be}\mathbb{T}%
_{dc}^{e}+g_{ec}\mathbb{T}_{bd}^{e}\right) ,%
\end{array}%
\leqno(6.4)
\end{equation*}%
\emph{are the components of a linear} $\rho $\emph{-connection compatible
with the (pseudo) metrical structure }$g,$ \emph{\ where }$\rho \Gamma
_{bc}^{a}$\emph{\ are the components of linear} $\rho $\emph{-connection of
Levi-Civita type }$\left( 6.3\right) $\emph{. Therefore, the vector bundle }$%
\left( h^{\ast }E,h^{\ast }\pi ,M\right) $\emph{\ becomes} $\rho $\emph{%
-(pseudo)metrizable and the tensor field }$\mathbb{T}$\emph{\ is the }$%
\left( \rho ,h\right) $\emph{-torsion tensor field.}

\emph{In the particular case of Lie algebroids, }$h=Id_{M},$\emph{\ }$g\in
\mathcal{T}_{2}^{0}\left( E,\pi ,M\right) $\emph{\ is a pseudo(metrical)
structure and }$T\in \mathcal{T}_{2}^{1}\left( E,\pi ,M\right) $\emph{\ such
that its components are skew symmetric in the lover indices, then the local
real functions }%
\begin{equation*}
\begin{array}{c}
\rho \mathring{\Gamma}_{bc}^{a}=\rho \Gamma _{bc}^{a}+\displaystyle\frac{1}{2%
}\tilde{g}^{ad}\left( g_{de}\mathbb{T}_{bc}^{e}-g_{be}\mathbb{T}%
_{dc}^{e}+g_{ec}\mathbb{T}_{bd}^{e}\right) ,%
\end{array}%
\leqno(6.4^{\prime })
\end{equation*}%
\emph{are the components of a linear} $\rho $\emph{-connection compatible
with the (pseudo)metrical structure }$g,$ \emph{where }$\rho \Gamma
_{bc}^{a} $\emph{\ are the components of linear} $\rho $\emph{-connection of
Levi-Civita type }$\left( 6.3^{\prime }\right) $\emph{. }

\emph{In the classical case, }$\left( \rho ,h\right) =\left(
Id_{TM},Id_{M}\right) ,$\emph{\ }$g\in \mathcal{T}_{2}^{0}\left( TM,\tau
_{M},M\right) $\emph{\ is a pseudo(metrical) structure and }$T\in \mathcal{T}%
_{2}^{1}\left( TM,\tau _{M},M\right) $\emph{\ such that its components are
skew symmetric in the lover indices, then the local real functions }%
\begin{equation*}
\begin{array}{c}
\mathring{\Gamma}_{jk}^{i}=\Gamma _{jk}^{i}+\displaystyle\frac{1}{2}\tilde{g}%
^{ih}\left( g_{he}\mathbb{T}_{jk}^{e}-g_{je}\mathbb{T}_{hk}^{e}+g_{ek}%
\mathbb{T}_{jh}^{e}\right) ,%
\end{array}%
\leqno(6.4^{\prime \prime })
\end{equation*}%
\emph{are the components of a linear} \emph{connection compatible with the
(pseudo)metrical structure }$g,$ \emph{where }$\Gamma _{jk}^{i}$\emph{\ are
the components of linear} \emph{connection of Levi-Civita type }$\left(
6.3^{\prime \prime }\right) $\emph{. }

\textbf{Theorem 6.3 }\emph{If }$(E,\pi ,M)=(F,\nu ,M),$\emph{\ }$g\in
\mathcal{T}_{2}^{0}\left( h^{\ast }E,h^{\ast }\pi ,M\right) $\emph{\ is a
pseudo (metrical) structure and }$\rho \mathring{\Gamma}$\emph{\ is the
linear }$\rho $\emph{-connection }$\left( 6.4\right) $\emph{\ for the vector
bundle }$\left( h^{\ast }E,h^{\ast }\pi ,M\right) $\emph{, then the local
real functions }%
\begin{equation*}
\begin{array}{c}
\rho \tilde{\Gamma}_{b\alpha }^{a}=\rho \mathring{\Gamma}_{b\alpha }^{a}+%
\frac{1}{2}\tilde{g}^{ac}g_{cb\overset{\circ }{\mid }\alpha }%
\end{array}%
\leqno(6.5)
\end{equation*}%
\emph{are the components of a linear }$\rho $\emph{-connection such that the
vector bundle }$\left( h^{\ast }E,h^{\ast }\pi ,M\right) $\emph{\ becomes} $%
\rho $\emph{-(pseudo)metrizable.}

\emph{In the particular case of Lie algebroids, }$h=Id_{M},$ $g\in \mathcal{T%
}_{2}^{0}\left( E,\pi ,M\right) $\emph{\ is a pseudo(metrical) structure and
}$\rho \mathring{\Gamma}$\emph{\ is the linear }$\rho $\emph{-connection }$%
\left( 6.4^{\prime }\right) $\emph{\ for the vector bundle }$\left( E,\pi
,M\right) $\emph{, then the local real functions }%
\begin{equation*}
\begin{array}{c}
\rho \tilde{\Gamma}_{b\alpha }^{a}=\rho \mathring{\Gamma}_{b\alpha }^{a}+%
\frac{1}{2}\tilde{g}^{ac}g_{cb\overset{\circ }{\mid }\alpha }%
\end{array}%
\leqno(6.5^{\prime })
\end{equation*}%
\emph{are the components of a linear }$\rho $\emph{-connection such that the
vector bundle }$\left( E,\pi ,M\right) $\emph{\ becomes} $\rho $\emph{%
-(pseudo)metrizable.}

\emph{In the classical case, }$\left( \rho ,h\right) =\left(
Id_{TM},Id_{M}\right) ,$\emph{\ }$g\in \mathcal{T}_{2}^{0}\left( TM,\tau
_{M},M\right) $\emph{\ is a pseudo(metrical) structure and }$\rho \mathring{%
\Gamma}$\emph{\ is the linear }$\rho $\emph{-connection }$\left( 6.4^{\prime
}\right) $\emph{\ for the vector bundle }$\left( TM,\tau _{M},M\right) $%
\emph{, then the local real functions }%
\begin{equation*}
\begin{array}{c}
\tilde{\Gamma}_{jk}^{i}=\mathring{\Gamma}_{jk}^{i}+\frac{1}{2}\tilde{g}%
^{ih}g_{hj\overset{\circ }{\mid }k}%
\end{array}%
\leqno(6.5^{\prime \prime })
\end{equation*}%
\emph{are the components of a linear connection such that the vector bundle }%
$\left( TM,\tau _{M},M\right) $\emph{\ becomes} \emph{(pseudo)metrizable.}

\textbf{Theorem 6.4 }\emph{If }$g\in \mathcal{T}_{2}^{0}\left( h^{\ast
}E,h^{\ast }\pi ,M\right) $\emph{\ is a pseudo (metrical) structure, }$\rho
\tilde{\Gamma}$\emph{\ is the linear }$\rho $\emph{-connection }$\left(
6.5\right) $\emph{\ for the vector bundle }$\left( h^{\ast }E,h^{\ast }\pi
,M\right) $\emph{, }$T=T_{c\alpha }^{d}S_{d}\otimes S^{c}\otimes t^{\alpha }$%
\emph{, and }$O_{bd}^{ca}=\frac{1}{2}\delta _{b}^{c}\delta _{d}^{a}-g_{bd}%
\tilde{g}^{ca}$ \emph{is the Obata operator, then the local real functions}
\begin{equation*}
\begin{array}{c}
\rho \hat{\Gamma}_{b\alpha }^{a}=\rho \tilde{\Gamma}_{b\alpha }^{a}+\frac{1}{%
2}O_{bd}^{ca}T_{c\alpha }^{d},%
\end{array}%
\leqno(6.6)
\end{equation*}%
\emph{are the components of a linear }$\rho $\emph{-connection such that the
vector bundle }$\left( h^{\ast }E,h^{\ast }\pi ,M\right) $\emph{\ becomes} $%
\rho $\emph{-(pseudo)metrizable. }

\emph{In the particular case of Lie algebroids, }$h=Id_{M},$\emph{\ }$g\in
\mathcal{T}_{2}^{0}\left( E,\pi ,M\right) $\emph{\ is a pseudo (metrical)
structure, }$\rho \tilde{\Gamma}$\emph{\ is the linear }$\rho $\emph{%
-connection }$\left( 6.5^{\prime }\right) $\emph{\ for the vector bundle }$%
\left( E,\pi ,M\right) $\emph{, }$T=T_{c\alpha }^{d}s_{d}\otimes
s^{c}\otimes t^{\alpha }$\emph{\ and and }$O_{bd}^{ca}=\frac{1}{2}\delta
_{b}^{c}\delta _{d}^{a}-g_{bd}\tilde{g}^{ca}$ \emph{is the Obata operator,
then the local real functions}
\begin{equation*}
\begin{array}{c}
\rho \hat{\Gamma}_{b\alpha }^{a}=\rho \tilde{\Gamma}_{b\alpha }^{a}+\frac{1}{%
2}O_{bd}^{ca}T_{c\alpha }^{d},%
\end{array}%
\leqno(6.6^{\prime })
\end{equation*}%
\emph{are the components of a linear }$\rho $\emph{-connection such that the
vector bundle }$\left( E,\pi ,M\right) $\emph{\ becomes} $\rho $\emph{%
-(pseudo)metrizable. }

\emph{In the classical case, }$\left( \rho ,h\right) =\left(
Id_{TM},Id_{M}\right) ,$\emph{\ }$g\in \mathcal{T}_{2}^{0}\left( TM,\tau
_{M},M\right) $ \emph{is a pseudo(metrical) structure, }$\tilde{\Gamma}$%
\emph{\ is the linear connection }$\left( 6.5^{\prime \prime }\right) $\emph{%
\ for the vector bundle }$\left( TM,\tau _{M},M\right) $\emph{, }$%
T=T_{hk}^{l}\partial _{l}\otimes dx^{h}\otimes dx^{k}$\emph{\ and and }$%
O_{jl}^{hi}=\frac{1}{2}\delta _{j}^{h}\delta _{l}^{i}-g_{jl}\tilde{g}^{hi}$
\emph{is the Obata operator, then the local real functions}
\begin{equation*}
\begin{array}{c}
\hat{\Gamma}_{jk}^{i}=\tilde{\Gamma}_{jk}^{i}+\frac{1}{2}%
O_{jl}^{hi}T_{hk}^{l},%
\end{array}%
\leqno(6.6^{\prime })
\end{equation*}%
\emph{are the components of a linear connection such that the vector bundle }%
$\left( TM,\tau _{M},M\right) $\emph{\ becomes} \emph{(pseudo)metrizable. }



\bigskip

\bigskip

\bigskip

\bigskip

\bigskip

\bigskip

\bigskip




\hfill

\end{document}